\documentclass[leqno,12pt]{article}  
\usepackage{amsmath}
\usepackage{amssymb}

\usepackage[german,english]{babel}
\usepackage{amsthm}

\title{Families of two dimensional modular $(\varphi,\Gamma)$-modules}

\author{\textsc{Elmar Grosse-Kl\"onne}}
\date{}

\theoremstyle{plain} 
\newtheorem{satz}{Theorem}[section]  
\newtheorem{kor}[satz]{Corollary}  
\newtheorem{lem}[satz]{Lemma}  
\newtheorem{pro}[satz]{Proposition}  
 
\newcommand{\q}{\mbox{\rm Frac}}  

\theoremstyle{remark}

\theoremstyle{definition}

\newcommand{\F}{\ensuremath{\mathbb{F}}}

\newcommand{\w}{\omega}

\newcommand{\0}{\ensuremath{\overrightarrow{0}}}

\begin{document}

\maketitle

%


\begin{abstract} Let $F/{\mathbb Q}_p$ be a finite unramified extension, let $k$ be a finite extension of the residue field of $F$. We explicitly construct integral structures for all rank two \'{e}tale Lubin-Tate $(\varphi,{\mathcal O}_F^{\times})$-modules over $k$. We construct algebraic families of such integral structures and show that these comprehensively reflect the degeneration behaviour of $(\varphi,{\mathcal O}_F^{\times})$-modules. These results reveal new combinatorial structures of the moduli stack of $(\varphi,{\mathcal O}_F^{\times})$-modules, and allow us, in particular, to rederive the fact that the Serre weights assigned to a two dimensional ${\rm Gal}(\overline{F}/F)$-representation over $k$ can be read off from the geometry of the stack.  
\end{abstract}
\tableofcontents

\vspace{1cm}

Let $F$ be an unramified field extension of ${\mathbb Q}_p$ with residue field ${\mathbb F}_q$, where $q=p^f<\infty$. Let $k$ be a finite field extension of ${\mathbb F}_q$. Write $\Gamma={\mathcal O}_F^{\times}$. Write ${\mathcal G}_{F}$, resp. ${\mathcal I}_{F}$, for the absolute Galois group of $F$, resp. its inertia subgroup. In this introduction, if not specified otherwise, $(\varphi,\Gamma)$-modules are \'{e}tale Lubin-Tate $(\varphi,\Gamma)$-modules of dimension $2$ over the Laurent series field $k((t))$. Their isomorphism classes are in bijection with the isomorphism classes of two dimensional continuous ${\mathcal G}_{F}$-representations over $k$.

{\bf (1)} In this paper we aim to

(a) describe explicitly all $(\varphi,\Gamma)$-modules (in fact, even $k[[t]]$-integral structures in their $k$-duals),

(b) describe explicitly certain classes of non-trivial (i.e. not simply being a variation of free parameters) families of such $(\varphi,\Gamma)$-modules (in fact, even $k[[t]]$-integral structures in their $k$-duals),

(c) identify underlying combinatorial structures of the families in (b). One way to describe these structures is in terms of the labelling by Serre weights of the irreducible components of the stack of $(\varphi,\Gamma)$-modules ${\mathcal X}_{2,{\rm red}}$ constructed by Matthew Emerton and Toby Gee \cite{eg}.

For $p>2$ we in particular rederive the main result of \cite{cegs} for unramified $F/{\mathbb Q}_p$ \footnote{The main result of \cite{cegs} allows finite extensions $F/{\mathbb Q}_p$ which are not necessarily unramified.}: the Serre weights attached to a two dimensional continuous ${\mathcal G}_{F}$-representation over $k$ (see \cite{bdj}) can be read off from the degeneration behaviour of families of $(\varphi,\Gamma)$-modules, or equivalently, from the intersection behaviour of the components of ${\mathcal X}_{2,{\rm red}}$. But in fact, we

(d) provide a stratification of ${\mathcal X}_{2,{\rm red}}$ which is finer than the one implied by the decomposition into irreducible components (i.e. the irreducible components are stratified further). The strata and their relative position on ${\mathcal X}_{2,{\rm red}}$ are given and analysed in terms of purely combinatorial data.

{\bf (2)} Let us give more details on (a), (b), (c), (d).

(a) Consider the ring $k[[t]][\varphi]$ with commutation relation $\varphi t=t^q\varphi$. Let $M$ be a free rank two $k[[t]][\varphi]$-module, with a specified basis. To describe a $(\varphi,\Gamma)$-module, we

${\bullet}$ specify an action of $\Gamma$ on $M$ by specifying its effect on the given basis and by using the Lubin-Tate action of $\Gamma$ on $k[[t]][\varphi]$, and

${\bullet}$ specify four particular elements inside $M$.

We then look at the quotient $\Delta$ of $M$ by the $k[[t]][\varphi]$-sub module generated by these four elements. For suitable such specifications we show that the action of $\Gamma$ on $M$ passes to an action on $\Delta$. Next $\Delta^*={\rm Hom}_k(\Delta,k)$ is naturally a $k[[t]]$-module (we remark in passing that $\Delta$ is torsion as a $k[[t]]$-module), and $\Delta^*((t))=\Delta^*\otimes_{k[[t]]}k((t))$ is a two-dimensional $k((t))$-vector space, which then naturally becomes a $(\varphi,\Gamma)$-module.\footnote{Here is a rank-one model for this construction: Consider the quotient $\Delta$ of $k[[t]][\varphi]$, viewed as a free rank-one module over itself, by the $k[[t]][\varphi]$-submodule generated by the two elements $t^{q-1}\varphi-1$ and $t$; then $\Delta^*$ is a free rank-one module over $k[[t]]$.}

We show that each $(\varphi,\Gamma)$-module arises from this construction.

As $\varphi$ acts on $\Delta$, its dual acts as a $k$-linear operator on the $k[[t]]$-lattice $\Delta^*$ in $\Delta^*((t))$. It is a $\psi$-operator with respect to $\varphi$. For $F={\mathbb Q}_p$ such $(\psi,\Gamma)$-stable lattices play a crucial role in Colmez's assignment of ${\rm GL}_2({\mathbb Q}_p)$-representations to $(\varphi,\Gamma)$-modules and vice versa, and one might ask if for general $F$ they can play a similar role. Here however we do not even mention the $\psi$-operator on $\Delta^*$, instead we directly pass from $\varphi$ on $\Delta$ to $\varphi$ on $\Delta^*((t))$.

Denote by $\omega:{\mathcal G}_F\to\Gamma\to{\mathbb F}_q^{\times}\subset k^{\times}$ the Lubin Tate
character associated with the uniformizer $p$ of ${\mathcal O}_F$ and fix a pair $\ell,u$ in $[0,q-2]\cong {\mathbb Z}/(q-1){\mathbb Z}$. Consider ${\mathcal G}_F$-representations $W$ which are {\it reducible} and sit, after restriction to ${\mathcal I}_F$, in an exact sequence$$0\longrightarrow \omega^{-u}\longrightarrow W\longrightarrow \omega^{\ell-u}\longrightarrow0.$$We may ask for all $(\varphi,\Gamma)$-modules $\Delta^*((t))$ as constructed above for which the corresponding ${\mathcal G}_F$-representation is of this sort. The pair $\ell,u$ specifies a $\Gamma$-action on a $k[[t]][\varphi]$-basis of $M$ (as above). Lifting the ${\mathcal I}_F$-actions on the two subquotients $\omega^{-u}$, $\omega^{\ell-u}$ of $W$ to ${\mathcal G}_F$-actions comes down to specifying two (one for each subquotient) elements in $k^{\times}$ (telling us how a Frobenius element of ${\mathcal G}_F$ acts); further, the extension class between the two subquotients is then given by $f$ (resp. $f+1$ in exceptional cases) arbitrary parameters in $k$. And indeed, the said four elements of $M$ which define $\Delta$ are then given in terms of these input data from $(k^{\times})^2\times k^f$. In this way we even exhibit a {\it canonical} basis of the ${\rm Ext}^1$-space, which allows us to regard these extension classes as elements in the $k$-vector space on the basis set ${\mathcal F}:=[0,f-1]\cong{\mathbb Z}/f{\mathbb Z}$ (and one more basis element in the exceptional cases). One of the main themes of this paper is then the obvious observation that subsets ${\mathcal I}\subset{\mathcal F}$ therefore single out specific subclasses of $(\varphi,\Gamma)$-modules, namely those for which the said parameters belong to the subspace spanned by the elements of ${\mathcal I}$. Let us say that these $(\varphi,\Gamma)$-modules (over $k$, and similarly over finite extensions of $k$) are of type $(\ell,u,{\mathcal I})$.\footnote{Notice that a $(\varphi,\Gamma)$-module of type $(\ell,u,{\mathcal I})$ is also of type $(\ell,u,{\mathcal I}')$ for each ${\mathcal I}'\supset{\mathcal I}$.} We then call an algebraic family of $(\varphi,\Gamma)$-modules an $(\ell,u,{\mathcal I})$-family if all its fibres at closed points are of type $(\ell,u,{\mathcal I})$.

(b) It follows that, when $\ell,u$ are fixed, $\Delta$ by construction varies algebraically over a parameter space ${\mathbb G}_m^2\times{\mathbb A}^f$ over $k$, i.e. ($\Delta$ and hence) $\Delta^*$ may be regarded as a flat family of $(\varphi,\Gamma)$-modules over ${\mathbb G}_m^2\times{\mathbb A}^f$. \footnote{In the exceptional cases there is one more parameter (in $k$) which we may safely neglect for the present purposes.} For such a family the restriction to ${\mathcal I}_F$ of the semisimplified associated ${\mathcal G}_F$-representations always remains the {\it same} (determined by $\ell$, $u$). But one may aim to extend such a family over larger parameter spaces $\overline{{\mathbb G}_m^2\times{\mathbb A}^f}$, containing ${\mathbb G}_m^2\times{\mathbb A}^f$ as an open subset, and with {\it new} ${\mathcal I}_F$-restricted semisimplified ${\mathcal G}_F$-representations (reducible or irreducible) appearing 'at the boundary'. We (somewhat implicity\footnote{In fact, we make no attempts to give closed formulae for the (maximal) $\overline{{\mathbb G}_m^2\times{\mathbb A}^f}$; they remain entirely in the background.}) describe such $\overline{{\mathbb G}_m^2\times{\mathbb A}^f}$-families by describing their restrictions to ${\rm Spec}(k[\tau])$ along certain well chosen morphisms ${\rm Spec}(k[\tau])\to\overline{{\mathbb G}_m^2\times{\mathbb A}^f}$ for which ${\rm Spec}(k[\tau^{\pm}])$ maps to ${\mathbb G}_m^2\times{\mathbb A}^f$. (We suggest to think of this as 'going to infinity' along selected directions in ${\mathbb A}^f$.) In other words, in the notations of sections \ref{okt3} and \ref{zoomgergelyfengwei}, we describe algebraic one-parameter families over $k[\tau]$, with $\tau\ne0$ being the locus where the ${\mathcal I}_F$-restriction of the semisimplified ${\mathcal G}_F$-representations is given by $\ell$, $u$ as before, and with $\tau=0$ being the locus which yields 'new' ${\mathcal I}_F$-restricted semisimplified ${\mathcal G}_F$-representations.

The construction of these families follows the same outline as before: four elements in a free $k[\tau][[t]][\varphi]$-module of rank two are specified, dividing out the $k[\tau][[t]][\varphi]$-sub module generated by them yields a $k[\tau][[t]][\varphi]$-module $\Delta$, and hence a $(\varphi,\Gamma)$-module ${\rm Hom}_{k[\tau]}(k[\tau][[t]],k[\tau])\otimes_{k[\tau][[t]]}k[\tau]((t))$ over $k[\tau]$. Notice that we really not just build families of $(\varphi,\Gamma)$-modules, but even families of integral structures for them. Here it seems to be of critical importance to work with 'optimal' such integral structures (in notations of the main text: to choose the parameters $n_x$, $n_y$ minimally), and our feeling is that any attempt to build interesting families (i.e. degenerations) of $(\varphi,\Gamma)$-modules will need to overcome computational obstructions which in one way or the other reflect our analysis of 'optimal' integral structures.

(c), (d) What we explained so far suggests to define, for any given $\widetilde{\ell},\widetilde{u},\overline{\ell},\overline{u}$ in $[0,q-2]$ and subsets $\widetilde{\mathcal I}$, $\overline{\mathcal I}$ of ${\mathcal F}$, the relation\begin{gather}(\widetilde{\ell},\widetilde{u},\widetilde{\mathcal I})\rhd(\overline{\ell},\overline{u},\overline{\mathcal I}),\label{intro0}\end{gather}by demanding
 that there is an algebraic family of $(\varphi,\Gamma)$-module over $k[\tau]$, as in (b) above, which over $k[\tau^{\pm}]$ is an $(\widetilde{\ell},\widetilde{u},\widetilde{\mathcal I})$-family but which at $\tau=0$ degenerates, as in (b) above, into an $(\overline{\ell},\overline{u},\overline{\mathcal I})$-family. Which $\widetilde{\ell},\widetilde{u}, \overline{\ell},\overline{u},\widetilde{\mathcal I}$, $\overline{\mathcal I}$ satisfy the relation (\ref{intro0}) ? Our conventions for $\ell$, $u$ are such that, if $\widetilde{\ell}$, $\overline{\ell}$ are given, then $\widetilde{u}$ (resp. $\overline{u}$) is uniquely determined by $\overline{u}$ (resp. $\widetilde{u}$), and the latter can be freely chosen if there is {\it some} relation (\ref{intro0}) at all. Moreover, a relation (\ref{intro0}) automatically specifies a subset $J\subset{\mathcal F}$ which is 'interval-like', i.e. is the image of an interval in ${\mathbb Z}$ under  ${\mathbb Z}\to {\mathbb Z}/f{\mathbb Z}\cong {\mathcal F}$. These two considerations transform the study of the relation (\ref{intro0}) into that of a similar relation\begin{gather}(\widetilde{\ell},\widetilde{\mathcal I})\succ_J(\overline{\ell},\overline{\mathcal I})\label{intro1}\end{gather}for $\widetilde{\ell}, \overline{\ell}$ in $[0,q-2]$ and subsets $\widetilde{\mathcal I}$, $\overline{\mathcal I}$, $J$ of ${\mathcal F}$, with $J$ being interval-like. 

We invest quite some effort into the study of the  relation (\ref{intro1}). First, we give a characterization of it which no longer involves $(\varphi,\Gamma)$-modules but is purely combinatorial. More precisely, we define a map $\nu_{\widetilde{\ell}}^{\overline{\ell}}:{\mathcal F}\to{\mathcal F}$ (depending only on $\widetilde{\ell}$, ${\overline{\ell}}$) and show that (\ref{intro1}) is equivalent with the conjunction of conditions\begin{gather}\widetilde{\ell}\succ_J\overline{\ell},\quad\quad\quad J^-\cup J^+\subset \widetilde{\mathcal I},\quad\quad\quad\overline{\mathcal I}\subset  \nu_{\widetilde{\ell}}^{\overline{\ell}}(\widetilde{\mathcal I}\cap J^{c,1})\label{intro2}\end{gather}with certain subsets $J^+, J^-, J^{c,1}$ of ${\mathcal F}$ derived from $J$, and a certain numerical relation $\widetilde{\ell}\succ_J\overline{\ell}$. But these conditions make sense without requiring $J$ to be interval-like, thus define relations (\ref{intro1}) for {\it any} $J$. In fact, one of our main combinatorial results is that (\ref{intro1}) for arbitrary $J$'s is precisely the transitive closure of (\ref{intro1}) for interval-like $J$'s: Theorem \ref{tobbu}. Accordingly, this yields a characterization of the transitive closure of the relation (\ref{intro0}). 

All this has parallels (relations similar to (\ref{intro0}) and (\ref{intro1}), and a combinatorial characterization similar to formula (\ref{intro2})) for degenerations of families of reducible $(\varphi,\Gamma)$-modules into irreducible $(\varphi,\Gamma)$-modules, with Theorem \ref{ostersa} complementing Theorem \ref{tobbu}.

Of course, the organizing principle and object in the background of our attempts to understand degenerations of $(\varphi,\Gamma)$-modules is the moduli stack ${\mathcal X}_{2,{\rm red}}$ of rank two $(\varphi,\Gamma)$-modules constructed by Matthew Emerton and Toby Gee \cite{eg}, or rather its isomorphic Lubin-Tate adaption due to Dat Pham \cite{pham}, even though as an {\it algebro geometric} object it makes only implicit appearance in our work. Each of the irreducible components of ${\mathcal X}_{2,{\rm red}}$ comes with a pair $\ell$, $u$ as above, in such a way that {\it generically} the $(\varphi,\Gamma)$-modules sitting on that component are $(\ell,u,{\mathcal F})$-families (in our terminology from above). This allows us to label that component by a Serre weight (irreducible $\overline{\mathbb F}_{p}[{\rm GL}_2({\mathbb F}_{q})]$-module) $V(\ell,u)$ constructed from $\ell$, $u$. (The Steinberg components of ${\mathcal X}_{2,{\rm red}}$ receive an additional label.) By the main result of \cite{cegs}\footnote{The main result of \cite{cegs} even applies to finite extensions $F/{\mathbb Q}_p$ which are not necessarily unramified.}, for a given two dimensional ${\mathcal G}_F$-representation $W$ over $k$, the set of Serre weights attached to it by the (proven) Serre weight conjecture is precisely the set of labels of the ${\mathcal X}_{2,{\rm red}}$-components on which the $(\varphi,\Gamma)$-module of $W$ lies. The proof in \cite{cegs} is indirect, employing auxiliary moduli stacks and the formalism of the Breuil-M\'{e}zard conjecture. Here we suggest an alternative proof which instead makes critical use of the combinatorics of the relations (\ref{intro0}) and (\ref{intro1}) (and their analogs for degeneration into irreducibles). It is based on the fact that each triple $(\ell,u,{\mathcal I})$ defines a subset of the set of $k$-points (or better: $\overline{\mathbb F}_p$-points) of ${\mathcal X}_{2,{\rm red}}$, namely, those $(\varphi,\Gamma)$-modules which are of type $(\ell,u,{\mathcal I})$ but not of type $(\ell,u,{\mathcal I}')$ for any ${\mathcal I}'\subsetneq{\mathcal I}$. This (together with the irreducible points) yields a stratification of ${\mathcal X}_{2,{\rm red}}$. This stratification is finer than the one implied by the decomposition into irreducible components. Namely, asking for the degeneration of $(\widetilde{\ell},\widetilde{u},\widetilde{\mathcal I})$-families into $(\overline{\ell},\overline{u},\overline{\mathcal I})$-families, for given $\widetilde{\ell},\widetilde{u},\overline{\ell},\overline{u},\widetilde{\mathcal I}$, $\overline{\mathcal I}$, (and similarly, asking for the degeneration of $(\widetilde{\ell},\widetilde{u},\widetilde{\mathcal I})$-families into a given irreducible $(\varphi,\Gamma)$-module) is finer than asking for the intersection behaviour of the ${\mathcal X}_{2,{\rm red}}$-components --- the point is that additionally $\widetilde{\mathcal I}$, $\overline{\mathcal I}$ (and not just $\widetilde{\ell},\widetilde{u},\overline{\ell},\overline{u}$) are taken into account.

At no point do we make genericity assumptions on the ${\mathcal G}_F$-representations resp. $(\varphi,\Gamma)$-modules; genericity assumptions would very drastically cut down the length of this paper. More technically, any $\ell\in[0,q-2]$ as above gives rise to a set decomposition ${\mathcal F}={\mathcal D}(\ell)\coprod{\mathcal E}(\ell)$. Among the four $\Delta$-defining elements in $M$ mentioned above there is a distinguished one (the element $\underline{R}$ in Proposition \ref{comme}) defined as a certain sum, with (besides other summands) a summand for each element of ${\mathcal F}$. The summands corresponding to elements in ${\mathcal D}(\ell)$ are much more benign than those corresponding to elements in ${\mathcal E}(\ell)$. For generic $\ell$ we have ${\mathcal F}={\mathcal D}(\ell)$. Similarly, for generic $\widetilde{\ell}$, $\overline{\ell}$ the above map $\nu_{\widetilde{\ell}}^{\overline{\ell}}$ is the identity map on ${\mathcal F}$, and the combinatorics of the relations (\ref{intro0}) and (\ref{intro1}) becomes almost trivial. 

{\bf (3)} In the introduction to \cite{cegs} we read "it seems hard to understand in any explicit
way which Galois representations arise as the limits of a family of extensions of given
characters, and the description of the sets $W(r)$ is very complicated." We view the present work as an attempt to respond to this challenge. We are (mildly) confident that our methods generalize to higher dimensional ${\mathcal G}_F$-representations. If so, this would constitute an approach towards the Serre weight conjecture in the spirit of the above somewhat skeptical citation.

As has been explained, our combinatorial analysis of the limiting behaviour and the Serre weight assignment for ${\mathcal G}_F$-representations relies on working with a {\it canonical} $k$-basis of the ${\rm Ext}^1$-space between two given ${\mathcal G}_F$-characters. This feature is not new; the most notable paper here is \cite{dediro}\footnote{As an aside, the calculations of section 8 in \cite{dediro} were a very helpful testing ground for ours.} in which such a basis is defined and analysed in purely Galois theoretic terms. For a comparison with ours see the discussion of formula (\ref{yufasha}). Yet, our findings which relate such canonical bases for {\it different} pairs of ${\mathcal G}_F$-characters seem to be one of the main advances of the present work.

{\bf (4)} In section \ref{ottolebt} we work out the relation (\ref{intro0}) in the case $f=2$, formulating it without reference to $(\varphi,\Gamma)$-modules. We also work out  for $f=2$ some of the combinatorial maps occuring later in the general context. In section \ref{secone} we begin with an investigation of the coefficients of Lubin-Tate-multiplication power series, we recall basics on semisimple $(\varphi,\Gamma)$-modules and classify those of rank $\le2$. In section \ref{sectwo} we explain how to convert $k[[t]][\varphi]$-modules $\Delta$ with $\Gamma$-action into \'{e}tale $(\varphi,\Gamma)$-modules $\Delta^*((t))$, and we work this out for rank $\le 2$. Subsection \ref{viervor} contains the construction explained in (a) above, in the case where the restrictions to ${\mathcal I}_F$ of the two one dimensional subquotients of the ${\mathcal G}_F$-representation differ; in our notations this is the case $\ell\ne 0$. As usual, it comes with an additional subdiscussion if the said difference is the cyclotomic character. Subsection \ref{secthreeone} contains the case $\ell=0$. Subsections \ref{viervor} and \ref{secthreeone} in particular yield our explicit description of all rank two $(\varphi,\Gamma)$-modules.\footnote{Recently, another explicit description of the $(\varphi,\Gamma)$-modules in question was independently given in \cite{yiwa} (with Lemma 2.1 of \cite{yiwa} being the parallel to our Lemma \ref{doju}). However, integral structures do not seem to be described in \cite{yiwa}.} Subsection \ref{stackstrataset} defines the set of strata of ${\mathcal X}_{2,{\rm red}}$ alluded to above. In section \ref{okt3} it is explained, as outlined in (b) above, how a family of reducible $(\varphi,\Gamma)$-modules can degenerate into an irreducible one. In the analogous section \ref{zoomgergelyfengwei}\footnote{From this section onwards we assume $p>2$. It might well be that everything generalizes without too much effort to the case $p=2$.} the degeneration into reducible $(\varphi,\Gamma)$-modules with new ${\mathcal I}_F$-semisimplification (i.e. with new pair $\ell$, $u$) is explained. Section \ref{combi} explores the relation (\ref{intro1}), i.e. the conditions (\ref{intro2}), as well as its analog relevant for degeneration into irreducible $(\varphi,\Gamma)$-modules. In section \ref{colook}, to a two dimensional ${\mathcal G}_F$-representation $W$ we assign a set ${\mathbb W}_{\rhd}(W)$ of Serre weights in terms of the relation (\ref{intro0}). The results of section \ref{combi} allow us to identify it with the set ${\mathbb W}_{\rm geom}(W)$ of Serre weights attached to $W$ via the stack ${\mathcal X}_{2,{\rm red}}$, as well as with the more traditionally defined set ${\mathbb W}_{BDJ}(W)$ (which, of course, is the same). For this a comparsion with the definitions of the paper \cite{dediro} is relevant. In section \ref{giovanni} we explain that, by construction, the elements of ${\mathbb W}_{\rhd}(W)$ carry canonically defined {\it multiplicities}; thus, ${\mathbb W}_{\rhd}(W)$ is canonically even a {\it multiset}. 

{\bf (5)} Besides calling for generalizations to higher dimensions, this work should have further applications in the present two dimensional case. For example, it seems to easily imply the main results of \cite{kan} (for $F/{\mathbb Q}_p$ unramified) and suggests generalizations of them. More importantly, we hope it to be relevant for the Breuil-M\'{e}zard conjecture.\\

{\bf Acknowledgments:} I am extremely grateful to Toby Gee for his patience, both in personal conversation as well as in email exchange, with my numerous (and, I am afraid, repetitive) questions on the moduli stack of $(\varphi,\Gamma)$-modules, and on which sort of degeneration phenomena precisely are to be expected. Huge thanks go to Stefano Morra for an invitation to the Laboratoire Analyse G\'{e}om\'{e}trie Alg\`{e}bre (Paris) as a CNRS invited professor, for subsequent discussions on the Serre weight conjecture and remarks on the present paper. First insights of this work were conceived at Paris. I thank Florian Herzig for detailed and helpful suggestions for improvements. I also thank Vytautas Paskunas. I am very grateful to the referee whose comments helped a lot in improving the overall presentation.\\

{\bf Notations:} Fix $f\in{\mathbb Z}_{>0}$. We write $${{{\mathcal F}}}=[0,f-1]$$and denote by $$\Pi:{\mathbb Z}\longrightarrow{{{\mathcal F}}}$$the map with $\Pi(n)-n\in f{\mathbb Z}$ for all $n\in{\mathbb Z}$. Let $F/{\mathbb Q}_p$ be the unramified finite field extension of degree $f$. Write $k_F={\mathbb F}_q$ for the
residue field of $F$, with $q=p^f$ elements. Let $k$ be a finite extension of $k_F$. Write $\Gamma={\mathcal O}_F^{\times}$. Every element ${\gamma}\in\Gamma$ factors uniquely as 
\begin{gather}\gamma=\gamma_1\gamma_2\quad\quad\mbox{ with }\quad[{\gamma}_1]=a_{{\gamma}_1,1}t\quad\mbox{ and }\quad a_{{\gamma}_2,1}=1.\label{floos}\end{gather}Namely, $\gamma_1$ is the Teichmuller lifting of the image of $\gamma$ in ${\mathbb F}_q^{\times}$. Put $${\xi}=q-1.$$Let $\pi$ be a uniformizer in $F$ (after Lemma \ref{doju} we take $\pi=p$). Let $\Phi(t)$ denote the power series in a
coordinate $t$ for the Lubin-Tate formal group for $F$ with respect to
$\pi$. We may and do assume $\Phi(t)=\pi t+t^q$. For $\gamma\in{\mathcal O}_F$ let $[\gamma]\in{\mathcal O}_F[t]]$ denote the formal
  power series describing the action of $\gamma$ in the Lubin-Tate formal group law
  associated with $\Phi$; we write $[\gamma]$ also for the induced element in $k[[t]]$. Via $\gamma*t=[\gamma]$ we have a (continuous) action of
$\Gamma$ on $k((t))$. Similarly, we have a (continuous) $k$-algebra endomorphism $\varphi$ of $k((t))$ satisfying $\varphi(t)=\Phi(t)=t^q$.

Let ${\mathcal G}_{F}={\rm Gal}(\overline{F}/F)$ be the absolute Galois group of $F$, let ${\mathcal I}_{F}$ denote its inertia subgroup. Denote by $\omega:{\mathcal G}_F\to\Gamma$ the Lubin Tate
  character associated with $\pi$.
By composing with $\Gamma\to{\mathbb F}_q^{\times}\subset k^{\times}$ we read
$\omega$ also as a character ${\mathcal G}_F\to k^{\times}$. For $\lambda\in k^{\times}$ denote by $\mu_{\lambda}:{\mathcal G}_F\to k^{\times}$ the unique character with $\mu_{\lambda}|_{{\mathcal I}_F}=1$ and sending a fixed (arithmetic) Frobenius element of ${\mathcal G}_F$ to $\lambda^{-1}$.

\section{The case $f=2$}

\label{ottolebt}

This section is meant to further illustrate the explanations from the preceding introduction in the case $f=2$ (and hence ${\mathcal F}=\{0,1\}$) by describing the construction of $\Delta$ as well as the relation $(\widetilde{\ell},\widetilde{u},\widetilde{\mathcal I})\rhd(\overline{\ell},\overline{u},\overline{\mathcal I})$, without reference to $(\varphi,\Gamma)$-modules. It also writes out in detail some of the combinatorial maps to be defined in the main text for general $f$.

For this section only, we distinguish three types of $0<\ell<q-1$. We write $\ell=m_0+m_1p$ with $0\le m_i\le p-1$ and say that $\ell$ is of type (a) if $0\le m_i<p-1$ for $i=0,1$, of type (b) if $0\le m_0<p-1$ and $m_1=p-1$, and of type (c) if $0\le m_1<p-1$ and $m_0=p-1$.\footnote{If we define $1\le a_i\le p$ with [$a_0\ne p$ or $a_1\ne p$] by asking $\ell\equiv -a_0-a_1p$ modulo $\xi{\mathbb Z}$, then $\ell$ is of type (a) if $1\le a_i<p$ for $i=0,1$, of type (b) if $(a_0,a_1)=(p,p-1)$ or [$a_1=p$ and $a_0\le p-2$], and of type (c) if $(a_0,a_1)=(p-1,p)$ or [$a_0=p$ and $a_1\le p-2$].}\\

\underline{(1) Construction of $\Delta$:}

Let $\ell=m_0+m_1p$ be as above. Let $u\in[0,q-2]$ and $\alpha, \beta\in k^{\times}$, let ${\bf b}_0, {\bf b}_1\in k$.  Consider the free $k[[t]][\varphi]$-module $M$ of rank two with basis $x$, $y$, and the $k[[t]][\varphi]$-submodule $\nabla$ of $M$ generated by the four elements\footnote{For a brief motivation concerning the number and the shape of these four generators cf. the comments below pertaining to $\Gamma$-stability, as well as the Remark following Lemma \ref{fassa}.}$$ty,\quad\quad\quad t^{n_x}x,\quad\quad\quad t^{\xi}\varphi y-\beta y,\quad\quad\quad R=t^{\xi}\varphi x-\alpha x+\alpha\beta^{1-N}({\bf b}_0t^{r_0}+{\bf b}_1t^{r_1})\varphi^Ny$$(the normalization factor $\alpha\beta^{1-N}$ is not important) where $n_x, N\in{\mathbb N}$ are sufficiently large, and where$$r_0=q^N-1+\ell-\xi,\quad\quad\quad\quad  r_1=q^N-1+p(m_1+m_0p-\xi).$$We let $\Gamma$ act semilinearly on $M$ by demanding$$\gamma*x=a_{\gamma,1}^{u+1}\frac{t}{[\gamma]}x,\quad\quad\quad\quad\gamma*y=a_{\gamma,1}^{u-\ell+1}\frac{t}{[\gamma]}y$$for $\gamma\in\Gamma$. We claim that if $\ell$ is of type (a) then $\nabla$ is stable under the action of $\Gamma$.

Let us explain this. It is enough to see that $\gamma\in\Gamma$ takes the four given generators of $\nabla$ again into $\nabla$. For the three generators $ty$, $t^{n_x}x$ and $t^{\xi}\varphi y-\beta y$ this is straightforward. It remains to see $\gamma*R\in\nabla$. If $\gamma$ is a Teichm\"uller lift (i.e. $[\gamma]=a_{\gamma,1}t$), then this is again straightforward (as a consequence of $r_i\equiv\ell$ modulo $\xi{\mathbb Z}$ for $i=1,2$). Thus we may assume that $\gamma$ is such that $a_{\gamma,1}=1$ (since a general element in $\Gamma$ is a product of such an element and a Teich\"uller lift). Since $R\in \nabla$ and $\frac{t}{[\gamma]}\in k[[t]]$ it will be enough to see $\gamma*R-\frac{t}{[\gamma]}R\in\nabla$. One easily checks$$\gamma*R-\frac{t}{[\gamma]}R\equiv\alpha\beta^{1-N}({\bf b}_0([\gamma]^{r_0}-\frac{t}{[\gamma]}t^{r_0})+{\bf b}_1([\gamma]^{r_1}-\frac{t}{[\gamma]})t^{r_1})\varphi^Ny$$modulo $\nabla$, so it will be enough to show $[\gamma]^{r_i}-\frac{t}{[\gamma]}t^{r_i}\varphi^Ny\in\nabla$ for $i=1,2$. For this it will be enough to show $$\frac{[\gamma]^{r_i+1}}{t}\varphi^Ny\equiv t^{r_i}\varphi^Ny\quad\quad{\rm modulo }\quad tk[[t]]y$$ (as $ty\in\nabla$). This requires the vanishing of the coefficients of various $t$-power monomials in the power series $[\gamma]^{r_i+1}$. These vanishings follow from two observations. First, from $[\gamma]\in t+t^qk[[t]]$ (i.e. the coefficients of $t^i$ in $[\gamma]$ vanish for $1<i<q$, Lemma \ref{doju}), and secondly, from\begin{gather}{r_i+1\choose m}\equiv 0\quad\mbox{ mod }p{\mathbb Z}\quad\quad\mbox{ for }m\in[1,p^i-1].\label{revf2}\end{gather}This latter vanishing is precisely due to the type (a) assumption.

Thus, $\nabla$ being stable under the action of $\Gamma$, it follows that $k[[t]][\varphi]$ and $\Gamma$ act on $\Delta=M/\nabla$. It turns out that $\Delta^*\otimes_{k[[t]]}k((t))$ thus becomes an \'{e}tale $(\varphi,\Gamma)$-module ${\bf D}$ of rank two.

In terminology from the introduction, ${\bf D}$ is an $(\ell,u,{\mathcal F})$-family; in terminology from subsection \ref{stackstrataset} it belongs to ${\mathcal V}(\ell,u,{\mathcal F})$. But for our purposes it is crucial to be more specific. Namely, the canonical basis of the ${\rm Ext}^1$-space mentioned in the introduction is precisely the one for which ${\bf b}_0, {\bf b}_1\in k$ are the two coordinates. Thus, if ${\bf b}_0=0$ (resp. if ${\bf b}_1=0$, resp. if ${\bf b}_0={\bf b}_1=0$), then ${\bf D}$ is even an $(\ell,u,\{1\})$-family (resp. an $(\ell,u,\{0\})$-family, resp. an $(\ell,u,\emptyset)$-family); in terminology from subsection \ref{stackstrataset} it belongs to ${\mathcal V}(\ell,u,\{1\})$ (resp. to ${\mathcal V}(\ell,u,\{0\})$, resp. to ${\mathcal V}(\ell,u,\emptyset)$). 

The ${\mathcal G}_F$-representations with \'{e}tale $(\varphi,\Gamma)$-module as just constructed, i.e. with $\ell$ of type (a), are classified in \cite{dediro} as being the "generic" ones.\\

For types (b) and (c) the congruence (\ref{revf2}) fails for either $r_0$ or $r_1$ (and no other number could replace that problematic $r_i$). However, in the cases (b) and (c) the essentially same construction works with a modified definition of $R$; namely, the problematic summand $t^{r_0}$ resp. $t^{r_1}$ in $R$ must be replaced by a somewhat more complex polynomial (and a few more congruences similar to (\ref{revf2}) are then needed).\footnote{Later, when $f\ge2$, we will collect the indices $i$ (in ${\mathcal F}$) for which $r_i$ satisfies formula (\ref{revf2}) into the index set ${\mathcal D}(\ell)$, and the remaining indices into the set ${\mathcal E}(\ell)$. Thus, the elements in ${\mathcal D}(\ell)$ correspond to "nice" summands in $R$, whereas those in ${\mathcal E}(\ell)$ correspond to more intricate ones. In generic cases ${\mathcal E}(\ell)$ is empty.}

Unless $\ell=(p-2)(p+1)$ (the cyclotomic case, which requires an additional construction as explained in the main text), each reducible \'{e}tale $(\varphi,\Gamma)$-module of rank two arises in this way.\\

\underline{(2) The relation $(\widetilde{\ell},\widetilde{u},\widetilde{\mathcal I})\rhd(\overline{\ell},\overline{u},\overline{\mathcal I})$:}

Assume $p>2$ for simplicity. Let $\widetilde{\ell},\overline{\ell}, \widetilde{u},\overline{u}\in[0,q-2]$ and $\widetilde{\mathcal I},\overline{\mathcal I}\subset{\mathcal F}$ be given. The relation $(\widetilde{\ell},\widetilde{u},\widetilde{\mathcal I})\rhd(\overline{\ell},\overline{u},\overline{\mathcal I})$ holds true in each of the following instances:

(I) $(\widetilde{\ell},\widetilde{u})=(\overline{\ell},\overline{u})$ and $\overline{\mathcal I}\subset \widetilde{\mathcal I}$. This corresponds to the (trivial from the definitions) fact that a ${\mathcal G}_{F}$-representation of type $(\widetilde{\ell},\widetilde{u},\overline{\mathcal I})$ is also of type $(\widetilde{\ell},\widetilde{u},\widetilde{\mathcal I})$ if $\overline{\mathcal I}\subset \widetilde{\mathcal I}$.  

(II) $\widetilde{\ell}\equiv\widetilde{u}-\overline{u}\equiv- \overline{\ell}$ modulo $\xi{\mathbb Z}$ and $\overline{\mathcal I}=\widetilde{\mathcal I}=\emptyset$. This corresponds to swapping the two direct summands in a reducible split two dimensional ${\mathcal G}_{F}$-representation.

(III) $\overline{\mathcal I}=\emptyset$ and $\overline{\ell}=0$, and either [$\widetilde{\ell}+2=\xi$ and $\overline{u}\equiv 1+\widetilde{u}$ and $\widetilde{\mathcal I}=\{0\}$] or [$\widetilde{\ell}+2p=\xi$ and $\overline{u}\equiv p+\widetilde{u}$ and $\widetilde{\mathcal I}=\{1\}$]. This corresponds to degenerations of families of reducible ${\mathcal G}_{F}$-representations almost all of whose members (as well as all its character twists) are ramified, into split reducible ${\mathcal G}_{F}$-representations which, after character twist, are unramified.

(IV) Assume $\widetilde{\ell}\ne 0$ and define $1\le \widetilde{a_i}\le p$ with [$\widetilde{a_0}\ne p$ or $\widetilde{a_1}\ne p$] by asking $\widetilde{\ell}\equiv -\widetilde{a_0}-\widetilde{a_1}p$ modulo $\xi{\mathbb Z}$. We require that $0<\overline{\ell}<q-1$ is given by one of the following:$$\overline{\ell}\equiv \widetilde{a_0}-\widetilde{a_1}p\quad{\rm modulo}\quad \xi{\mathbb Z}\quad\quad\quad(\mbox{case } (1)),$$$$\overline{\ell}\equiv -\widetilde{a_0}+\widetilde{a_1}p\quad{\rm modulo}\quad \xi{\mathbb Z}\quad\quad\quad(\mbox{case } (2)).$$We require $\widetilde{\mathcal I}={\mathcal F}$ and\begin{gather}\overline{\mathcal I}=\{0\}\quad\quad\mbox{ if }\overline{\ell}\mbox{ is of type }(c)\mbox{ or }[\overline{\ell}\mbox{ is of type }(a)\mbox{ and we are in case } (1)],\notag\\\overline{\mathcal I}=\{1\}\quad\quad\mbox{ if }\overline{\ell}\mbox{ is of type }(b)\mbox{ or }[\overline{\ell}\mbox{ is of type }(a)\mbox{ and we are in case } (2)].\notag\end{gather}We require $2(\widetilde{u}-\overline{u})\equiv\widetilde{\ell}-\overline{\ell}$ modulo $\xi{\mathbb Z}$ and$$\widetilde{u}-\overline{u}-\widetilde{\ell}\in\xi{\mathbb Z}+p[0,p+1]\quad\quad\quad \mbox{ in case } (1),$$$$\widetilde{u}-\overline{u}-\widetilde{\ell}\in\xi{\mathbb Z}+[0,p+1]\quad\quad\quad \mbox{ in case } (2).$$This case (IV) corresponds to degenerations of families of reducible ${\mathcal G}_{F}$-representations into ${\mathcal G}_{F}$-representations with {\it new} semisimplified ${\mathcal I}_{F}$-restriction, but not yet covered by (III). (Here all ${\mathcal G}_{F}$-representations (as well as all its character twists) in the family, including the special fibre, are ramified.) \\

{\bf Remarks:} For $0<\ell<q-1$ let $\sigma^{\ell}:{\mathbb Z}\to\{1,2\}$ and $i_{\mathcal D}^{\ell}:{\mathbb Z}\to\{0,1\}={\mathcal F}$ as in the main text.

For $\ell$ of type (a) we find ${\mathcal D}(\ell)=\{0,1\}$ and ${\mathcal E}(\ell)=\emptyset$, as well as $i_{\mathcal D}^{\ell}(x)\equiv x$ modulo $2{\mathbb Z}$ for all $x$, and $\sigma^{\ell}(j)=1$ for all $j$. 

For $\ell$ of type (b) we find ${\mathcal D}(\ell)=\{1\}$ and ${\mathcal E}(\ell)=\{0\}$, as well as $i_{\mathcal D}^{\ell}(x)=1$ for all $x$, and $\sigma^{\ell}(j)=1$ if $j\in 1+2{\mathbb Z}$, but $\sigma^{\ell}(j)=2$ if $j\in 2{\mathbb Z}$.

For $\ell$ of type (c) we find ${\mathcal D}(\ell)=\{0\}$ and ${\mathcal E}(\ell)=\{1\}$, as well as $i_{\mathcal D}^{\ell}(x)=0$ for all $x$, and $\sigma^{\ell}(j)=1$ if $j\in 2{\mathbb Z}$, but $\sigma^{\ell}(j)=2$ if $j\in 1+2{\mathbb Z}$.

Now let again $\widetilde{\ell}$ and $\overline{\ell}$ be as in (IV) above. In case (1) we observe\begin{align}\overline{\ell}\equiv \widetilde{a_0}-\widetilde{a_1}p&\equiv -(p-\widetilde{a_0})-(\widetilde{a_1}-1)p\quad\quad&\mbox{ if }\widetilde{a_1}>1\mbox{ and }\widetilde{a_0}<p\notag\\{}&\equiv-(p-1)-(p-1)p\quad\quad&\mbox{ if }\widetilde{a_1}=1\mbox{ and }\widetilde{a_0}=p\notag\\{}&\equiv-(p-\widetilde{a_0}-1)-pp\quad\quad&\mbox{ if }\widetilde{a_1}=1\mbox{ and }\widetilde{a_0}<p-1\notag\\{}&\equiv-p-(p-1)p\quad\quad&\mbox{ if }\widetilde{a_1}=1\mbox{ and }\widetilde{a_0}=p-1\notag\\{}&\equiv-p-(\widetilde{a_1}-2)p\quad\quad&\mbox{ if }\widetilde{a_1}>2\mbox{ and }\widetilde{a_0}=p\notag\\{}&\equiv-(p-1)-pp\quad\quad&\mbox{ if }\widetilde{a_1}=2\mbox{ and }\widetilde{a_0}=p.\notag\end{align}It follows that, in the first two lines, $\overline{\ell}$ is of type (a). In the next two lines, $\overline{\ell}$ is of type (b). In the last two lines, $\overline{\ell}$ is of type (c). Similar formulae hold in case (2).

We derive: If $\widetilde{\ell}$ is of type (b) and $\overline{\ell}$ is of type (c) then we are in case (1), and in fact $(\widetilde{a_0},\widetilde{a_1})=(p,p-1)$. If $\widetilde{\ell}$ is of type (c) and $\overline{\ell}$ is of type (b) then we are in case (2), and in fact $(\widetilde{a_0},\widetilde{a_1})=(p-1,p)$. If $\widetilde{\ell}$ is of type (a) and $\overline{\ell}$ is of type (b) then we are in case (1). If $\widetilde{\ell}$ is of type (a) and $\overline{\ell}$ is of type (c) then we are in case (2). It is impossible that $\widetilde{\ell}$ and $\overline{\ell}$ are both of type (b), or both of type (c). The generic case is, of course, where $\widetilde{\ell}$ and $\overline{\ell}$  are both of type (a).

From all this, it can be read off that the map $\nu^{\overline{\ell}}_{\widetilde{\ell}}:{\mathcal F}=\{0,1\}\to {\mathcal F}=\{0,1\}$ of the main text is given by$$\nu^{\overline{\ell}}_{\widetilde{\ell}}={\rm id}\quad\quad\quad\quad\,\,\quad\mbox{ if }\overline{\ell}\mbox{ is of type }(a),$$$$\nu^{\overline{\ell}}_{\widetilde{\ell}}(0)=\nu^{\overline{\ell}}_{\widetilde{\ell}}(1)=0\quad\mbox{ if }\overline{\ell}\mbox{ is of type }(c),$$$$\nu^{\overline{\ell}}_{\widetilde{\ell}}(0)=\nu^{\overline{\ell}}_{\widetilde{\ell}}(1)=1\quad\mbox{ if }\overline{\ell}\mbox{ is of type }(b).$$In later notations, in case (1) we have $J=\{0\}$ and $i_1\equiv 1$ and $i_2\equiv 0$ modulo $f{\mathbb Z}$, whereas in case (2) we have $J=\{1\}$ and $i_1\equiv 0$ and $i_2\equiv 1$ modulo $f{\mathbb Z}$. We have $\overline{\mathcal I}=\{\nu^{\overline{\ell}}_{\widetilde{\ell}}(0)\}$ in case (1) and $\overline{\mathcal I}=\{\nu^{\overline{\ell}}_{\widetilde{\ell}}(1)\}$ in case (2).\\

{\bf Remark:} What is representative from this $f=2$-discussion for the case of general $f\ge2$ is the case above distinction (I)-(IV) (together with its qualitative meaning in terms of ${\mathcal G}_{F}$-representations). The all important maps $\nu^{\overline{\ell}}_{\widetilde{\ell}}$ hold the key for relating $\widetilde{\mathcal I}$ with $\overline{\mathcal I}$. They can be fairly subtle (and chaotic) in general. In the present $f=2$ case we used the type distinction (a), (b), (c) in order to work out $\nu^{\overline{\ell}}_{\widetilde{\ell}}$ explicitly. In the general case we do not formally proceed via such a type distinction, but encode everything into combinatorics of the sets ${\mathcal D}(\ell)$ and ${\mathcal E}(\ell)$.

\section{$(\varphi,\Gamma)$-modules}

\label{secone}

In the following Lemma \ref{doju} we prove a little bit more than we really need later on; what we really need are formulae (\ref{ottolaender}), (\ref{freivorquinqua1}) and (\ref{amrit00}).

\begin{lem}\label{doju} Assume $F\ne{\mathbb Q}_2$. Let ${\mathcal N}=1+{\mathbb Z}_{\ge0}{\xi}=\{1, q, 2q-1, 3q-2,\ldots\}$. For each $\gamma\in\Gamma$ we have$$[\gamma]=\sum_{i\in {\mathcal N}}a_{\gamma,i}t^{i}$$with certain $a_{\gamma,i}\in {\mathcal O}_F$ for $i\in {\mathcal N}$. We
    have\begin{align}a_{\gamma,1}&=\gamma,\label{ottolaender}\\a_{\gamma,q}&=\frac{{\gamma}^q-\gamma}{\pi^q-\pi}=\frac{a_{\gamma,1}^q-a_{\gamma,1}}{\pi^q-\pi},\label{ruofra}\\a_{\gamma,2q-1}&=q\frac{\pi^{\xi}-a_{\gamma,1}^{\xi}}{\pi-\pi^{2q-1}}a_{\gamma,q}.\label{freivorquinqua}\end{align}For $\gamma,\gamma'\in \Gamma$ we
    have \begin{gather}\frac{a_{\gamma\gamma',q}}{a_{\gamma\gamma',1}}=\frac{a_{\gamma,q}}{a_{\gamma,1}}+\frac{a_{\gamma',q}}{a_{\gamma',1}}.\label{freivorquinqua1}\end{gather}Let
$$[0,q]^{\mathcal N}=\{k_{\bullet}=(k_i)_{i\in {\mathcal N}}\,|\, k_i\in
    [0,q]\mbox{ and }k_i=0\mbox{ for }i>>0\}.$$For $n\in {\mathcal N}$
    with $n>q$ we have\begin{gather}a_{\gamma,n}=\frac{1}{\pi-\pi^n}(\sum_{m\in{\mathcal N}\atop
      q^{-1}n\le m<n}{m\choose
      {\xi}^{-1}(n-m)}\pi^{m-{\xi}^{-1}(n-m)}a_{\gamma,m}-\sum_{k_{\bullet}\in
      [0,q]^{\mathcal N}\atop\sum_ik_i=q,
        \sum_iik_i=n}\frac{q!}{\prod_i(k_i!)}\prod_ia_{\gamma,i}^{k_i}).\label{purgatoire1}\end{gather}If $F\ne{\mathbb Q}_p$ then\begin{gather}a_{\gamma,(rp^{f-1}+1){\xi}+1}=-{p\choose
        r}\frac{1}{\pi}a_{\gamma,1}(\frac{a_{\gamma,q}}{a_{\gamma,1}})^{rp^{f-1}}\quad\mbox{ in
      }k\label{purgatoire}\end{gather}for $0<r<p$, but $a_{\gamma,n}=0$ in $k$ for each other $n\in {\mathcal N}$ with $q^2>n>q$.         
\end{lem}

{\sc Proof:} Write $[\gamma]=\sum_{i\ge 1}a_{\gamma,i}t^{i}$ with certain $a_{\gamma,i}\in {\mathcal O}_F$ for $i\ge 1$. Given that $\Phi(t)=\pi t+t^q$, the $a_{\gamma,i}$ are characterized by
  the equality$$\pi\sum_{i\ge1}a_{\gamma,i}t^i+(\sum_{i\ge1}a_{\gamma,i}t^i)^q=\sum_{i\ge1}a_{\gamma,i}(\pi t+t^q)^i.$$First, this implies $a_{\gamma,i}=0$ whenever $i\notin {\mathcal N}$. The characterizing equation then reads$$\pi\sum_{i\in{\mathcal N}}a_{\gamma,i}t^i+(\sum_{i\in{\mathcal N}}a_{\gamma,i}t^i)^q=\sum_{i\in{\mathcal N}}a_{\gamma,i}(\pi t+t^q)^i.$$Equating the coefficients of $t^q$ gives $a_{\gamma,q}=\frac{a_{\gamma,1}^q-a_{\gamma,1}}{\pi^q-\pi}$. For $n\in {\mathcal N}$
    with $n>q$, equating the coefficients of $t^n$ gives $$\pi a_{\gamma,n}+\sum_{k_{\bullet}\in[0,q]^{\mathcal N}\atop\sum_ik_i=q,
        \sum_iik_i=n}\frac{q!}{\prod_i(k_i!)}\prod_ia_{\gamma,i}^{k_i}=\sum_{m\in{\mathcal N}\atop
      q^{-1}n\le m\le n}{m\choose
      {\xi}^{-1}(n-m)}\pi^{m-{\xi}^{-1}(n-m)}a_{\gamma,m}$$and hence formula
    (\ref{purgatoire1}). To derive formula (\ref{purgatoire}) from formula (\ref{purgatoire1}) notice
first that none of the summands $$(\pi-\pi^n)^{-1}{m\choose
  {\xi}^{-1}(n-m)}\pi^{m-{\xi}^{-1}(n-m)}a_{\gamma,m}$$ (for $m, n\in{\mathcal N}$
with $q^{-1}n\le m<n$ and  $q^2>n>q$) contributes to the class of $a_{\gamma,n}$
in $k$: Indeed, the conditions on $n,m$ exclude the case $n=qm$, and if
$n=qm-q+1$ then they imply $m=q$, in which case ${m\choose
  {\xi}^{-1}(n-m)}=q\in (\pi^2)$ by our assumption
$F\ne{\mathbb Q}_p$, hence $(\pi-\pi^n)^{-1}{m\choose
  {\xi}^{-1}(n-m)}\pi^{m-{\xi}^{-1}(n-m)}\in (\pi)$. Concerning the
summands\begin{gather}(\pi-\pi^n)^{-1}(-\frac{q!}{\prod_i(k_i!)}\prod_ia_{\gamma,i}^{k_i})\label{DovorDembe}\end{gather}notice
first that, as the fixed $k_{\bullet}\in [0,q]^{\mathcal N}$ satisfies
$\sum_ik_i=q$ and $\sum_iik_i=n$, and as $q^2>n>q$, there are $i_1\ne i_2$ with $k_{i_1}\ne 0$ and $k_{i_2}\ne 0$. Moreover, we then have
      $\frac{q!}{\prod_i(k_i!)}\notin p^2{\mathbb Z}$ if and only if $k_i=0$
      for all $i$ with $i\ne i_1$ and $i\ne i_2$, as well as $k_{i_1}=rp^{f-1}$ and
      $k_{i_2}=q-rp^{f-1}$ for some $0<r<p$ (use $\sum_ik_i=q$ again). Writing $n=1+a{\xi}$, $i_1=1+j_1{\xi}$,
      $i_2=1+j_2{\xi}$, the condition $\sum_iik_i=n$
      reads$$(1+j_1{\xi})rp^{f-1}+(1+j_2{\xi})(q-rp^{f-1})=1+a{\xi}$$which is
      equivalent with$$a=1+p^{f-1}(j_1r+j_2(p-r)).$$Since $q<n<q^2$ means $1<a\le q$
      we thus get $0<p^{f-1}(j_1r+j_2(p-r))\le q-1$, and hence
      $0< j_1r+j_2(p-r)<p$. This implies either $j_1=0$ or $j_2=0$. Swapping $i_1$,
      $i_2$ if necessary, we may assume $j_2=0$ (and hence $i_2=1$) and $1\le j_1<p$. At this
      point it is already clear that no summand (\ref{DovorDembe}) can occur
      if $q<n<1+p{\xi}$, thus $a_{\gamma,n}=0$ if $q<n<p{\xi}$. Applying
      this insight for $i_1$ instead of $n$ shows $a_{\gamma,i_1}=0$ if
      $1<j_1<p$. Thus, our summand (\ref{DovorDembe}) under discussion does
      not contribute to $a_{\gamma,n}=0$ if $1<j_1<p$, but possibly only if
      $j_1=1$ (and hence $i_1=q$). Finally, we observe$$\frac{q!}{(rp^{f-1})!(q-rp^{f-1})!}\equiv {p\choose
        r}\quad\mbox{ modulo }p^2{\mathbb Z}$$ as well as
    $a_{\gamma,1}^{q-rp^{f-1}}a_{\gamma,q}^{rp^{f-1}}=a_{\gamma,1}(\frac{a_{\gamma,q}}{a_{\gamma,1}})^{rp^{f-1}}$
in $k$. We have shown formula (\ref{purgatoire}).\hfill$\Box$\\

{\bf Remark:} If $F={\mathbb Q}_2$ then $\Phi(t)={{p}}t+t^q=2t+t^2$ defines the multiplicative group, and this violates
the conclusion of Lemma \ref{doju}.\\

From now on we assume $\pi=p$. For $\gamma\in\Gamma$ and $n\ge1$ let $a_{\gamma,n}\in k$ be as defined in Lemma \ref{doju}. It follows from Lemma \ref{doju} that if $a_{\gamma,1}=1$ then\footnote{Formula (\ref{amrit00}) was obtained independently of Lemma 2.1 in \cite{yiwa}.}\begin{gather}[\gamma]-t-a_{\gamma,q}t^q+a_{\gamma,q}^{p^{f-1}}t^{(p^{f-1}+1){\xi}+1}\in t^{(2p^{f-1}+1){\xi}+1}k[[t]]\label{amrit00}\end{gather}and in particular\begin{gather}[\gamma]-t-a_{\gamma,q}t^q\in t^{(p^{f-1}+1){\xi}+1}k[[t]].\label{amrit}\end{gather}Without further comment we often make use of formula (\ref{amrit}) by using its following implication: For $M,N\in{\mathbb Z}_{\ge0}$ we have\begin{gather}(\frac{[\gamma]}{t})^M\equiv1\mbox{ modulo }t^{N}k[[t]]\quad\mbox{ if }{M\choose m}\equiv 0\mbox{ modulo }p{\mathbb Z}\mbox{ for all }1\le m<\frac{N}{\xi}.\label{mariaeverkuendigung}\end{gather}In our discussion of formula (\ref{revf2}) we tried to explain why this formula is important for our purposes.\\

{\bf Definition:} An \'{e}tale $(\varphi,\Gamma)$-module over
$k((t))$ is a finite dimensional $k((t))$-vector space ${\bf D}$, together with a semilinear action $\Gamma\times{\bf D}\to {\bf D}$, $(\gamma,d)\mapsto \gamma*d$, of $\Gamma$ on ${\bf D}$, and a semilinear endomorphism $\varphi$ of ${\bf D}$ which commutes with the action by $\Gamma$ and is such that its $k((t))$-linearization$${\rm id}\otimes\varphi:k((t))\otimes_{\varphi,k((t))} {\bf
  D}\stackrel{\cong}{\longrightarrow}{\bf
  D}$$is bijective. \\

{\bf Remarks:} (1) Throughout this text, by a $(\varphi,\Gamma)$-module we understand an \'{e}tale $(\varphi,\Gamma)$-module over
$k((t))$, if not specified otherwise.\footnote{There will indeed show up \'{e}tale $(\varphi,\Gamma)$-modules over other base rings, but we will then say so.} Similarly, by a ${\mathcal G}_F$-representation we understand a continuous representation of ${\mathcal G}_F$ on a finite dimensional $k$-vector space. 

(2) We use the symbol $*$ to denote the $\Gamma$-action on ${\bf D}$ in order to distinguish it from the scalar action (via the natural projection  $\Gamma\to {\mathbb F}_q^{\times}\subset k$) of $\Gamma$ on the $k$-vector space underlying ${\bf D}$; this latter action we simply denote by $(\gamma,d)\mapsto \gamma d$.

\begin{satz}\label{sosego} (Fontaine, Kisin-Ren, Schneider \cite{peterlec}) There are functors ${\bf D}\mapsto W({\bf
  D})$ and $W\mapsto {\bf
  D}(W)$ setting up an equivalence between the category of $(\varphi,\Gamma)$-modules ${\bf D}$ and the category of ${\mathcal G}_F$-representations $W$; they preserve vector space dimensions.  
\end{satz}

\begin{lem}\label{allerseelmai} Let $\lambda\in k^{\times}$ and $a\in {\mathbb Z}$. Denote by $W$ the one dimensional ${\mathcal G}_F$-representation on which ${\mathcal G}_F$ acts through the character
$\omega^a\mu_{\lambda}:{\mathcal G}_F\to k^{\times}$. There is a $k((t))$-basis element $y$ of ${\bf D}(W)$ such that
$\varphi(y)=\lambda y$ and $\gamma* y=a_{\gamma,1}^ay$ for all $\gamma\in\Gamma$. Moreover, ${\bf D}(W)$ is (up to isomorphism) the unique \'{e}tale
$(\varphi,\Gamma)$-module ${\bf D}$ of $k((t))$-dimension one admitting a basis element $y$ with $\varphi(y)=\lambda y$ and $\gamma* y-a_{\gamma,1}^ay\in tk[[t]]y$ for all
$\gamma\in\Gamma$.
\end{lem}

{\sc Proof:} See \cite{ps}. The last statement follows from the first one since the $\omega^a\mu_{\lambda}$ (for all $\lambda$ and $a$) exhaust the continuous characters ${\mathcal G}_F\to k^{\times}$.\hfill$\Box$\\
  
\begin{lem}\label{endosterkreis} Let $h\in{\mathbb Z}_{>0}$ and $\alpha\in
  k^{\times}$. Let ${\bf D}$ be a two-dimensional \'{e}tale
  $(\varphi,\Gamma)$-module such that for some non-zero $x\in {\bf D}$ we have$$\varphi^2(x)=\alpha t^{-h{\xi}}x,\quad\quad \gamma*x=x$$for $\gamma\in\Gamma$ with $[\gamma]=a_{\gamma,1}t$.

  (a) If $h\notin(q+1){\mathbb Z}_{>0}$ then ${\bf D}$ is irreducible. The isomorphism class of the \'{e}tale $(\varphi,\Gamma)$-module ${\bf D}$ remains unchanged if $h$ is replaced by some $h'\in{\mathbb Z}_{>0}-(q+1){\mathbb Z}_{>0}$ with either $h\equiv h'$ modulo $(q^2-1){\mathbb Z}$ or with $qh\equiv h'$ modulo $(q^2-1){\mathbb Z}$.

  (b) If $h=(q+1)v$ for some $v\in{\mathbb Z}$ and if $\alpha$ has a square root in $k$, then the corresponding ${\mathcal G}_F$-representation is the direct sum of two copies of the same one-dimensional ${\mathcal G}_F$-representation which, when restricted to ${\mathcal I}_F$, is the character $\omega^v$.
  
\end{lem}
  
{\sc Proof:} For (a) see e.g. \cite{ps}. For (b) let $\beta\in k$ be such that $\beta^2=\alpha$ and look at the two submodules generated by $\beta t^{-v{\xi}}x\pm\varphi x$.\hfill$\Box$\\

\section{From $t$-torsion $k[[t]][\varphi,\Gamma]$-modules to $(\varphi,\Gamma)$-modules}

\label{sectwo}

We define the $k$-algebra $k[[t]][\varphi,\Gamma]$ with
commutation rules given by$$\gamma
\varphi=\varphi\gamma,\quad\quad \gamma
t=[\gamma]_{\Phi}(t)\gamma,\quad\quad\varphi t=\Phi(t) \varphi$$for $\gamma\in\Gamma$. Of course, $\Phi(t)=[\pi]_{\Phi}(t)=t^q$ in $k[[t]]$. In $k[[t]][\varphi,\Gamma]$ we have the subalgebras $k[[t]][\varphi]$ and $k[[t]][\Gamma]$. For a $k[[t]][\varphi,\Gamma]$-module $\Delta$ endow $\Delta^*={\rm Hom}_k(\Delta,k)$ with the structure of a $k[[t]][\Gamma]$-module by means
of $$(S\ell)(x)=\ell(Sx),\quad\quad (\gamma* \ell)(x)=\ell(\gamma^{-1}* x)$$
for $\ell\in\Delta^*$, $x\in\Delta$, $S\in k[[t]]$,
$\gamma\in\Gamma$. Consider the following ${k((t))}$-linear
maps:\begin{gather}\Delta^*\otimes_{k[[t]]}{k((t))}\longrightarrow
  (\Delta\otimes_{k[[t]],\varphi}k[[t]])^*\otimes_{k[[t]]}{k((t))},\label{pruefeaus0}\\\ell\otimes1\mapsto[x\otimes b\mapsto
    \ell(b\varphi(x))]\otimes1,\notag\\\Delta^*\otimes_{k[[t]],\varphi}{k((t))}\longrightarrow
  (\Delta\otimes_{k[[t]],\varphi}k[[t]])^*\otimes_{k[[t]]}{k((t))},\label{pruefeaus1}\\\ell\otimes1\mapsto[x\otimes
    t^{\kappa}\mapsto
    \left\{\begin{array}{l@{\quad:\quad}l}\ell(x)&\kappa=\xi\\0&0\le \kappa\le q-2\notag \end{array}\right\}]\otimes 1\end{gather}where
  $\ell\in\Delta^*$, $x\in\Delta$, $b\in k[[t]]$, and where we use the
  decomposition$$\Delta\otimes_{k[[t]],\varphi}k[[t]]=\bigoplus_{0\le\kappa<q}\Delta\otimes t^{\kappa}$$which results from
  the decomposition$$k[[t]]=\bigoplus_{0\le\kappa<q}t^{\kappa}\varphi k[[t]].$$

\begin{satz}\label{semferiend} Suppose that, as a $k[[t]]$-module, $\Delta^*$
  is free of finite rank $n$. The
  maps (\ref{pruefeaus0}) and (\ref{pruefeaus1}) are bijective. The
  composition$$\Delta^*\otimes_{{k[[t]]},\varphi}{k((t))}\longrightarrow
  \Delta^*\otimes_{{k[[t]]}}{k((t))}$$of (\ref{pruefeaus1}) with the inverse
  of (\ref{pruefeaus0}) thus
  induces\footnote{use the natural isomorphism $\Delta^*\otimes_{{k[[t]]}}{k((t))}\otimes_{{k((t))},\varphi}{k((t))}\cong
  \Delta^*\otimes_{{k[[t]]},\varphi}{k((t))}$} a $\varphi$-linear endomorphism$$\varphi:\Delta^*\otimes_{{k[[t]]}}{k((t))}\longrightarrow \Delta^*\otimes_{{k[[t]]}}{k((t))}.$$Together with the $\Gamma$-action obtained from that on $\Delta^*$ this defines on $\Delta^*\otimes_{{k[[t]]}}{k((t))}$ the structure of an \'{e}tale $(\varphi,\Gamma)$-module of ${k((t))}$-dimension $n$.
  \end{satz}

{\sc Proof:} $\Delta^*\otimes_{{k[[t]]}}{k((t))}$ is an $n$-dimensional ${k((t))}$-vector space. It follows that the sources and targets of the
$k((t))$-linear maps (\ref{pruefeaus1}) and (\ref{pruefeaus0}) all are $n$-dimensional ${k((t))}$-vector spaces. Since these maps are easily
seen to be injective, they are thus bijective.\hfill$\Box$\\

\begin{lem}\label{sur} If $\Delta$ is a $k[[t]]$-torsion module and $t:\Delta\to\Delta$ is surjective, with
  kernel of finite $k$-dimension $n$, then $\Delta^*$
  is free of rank $n$ as a $k[[t]]$-module. 
  \end{lem}

{\sc Proof:} As $t$ acts surjectively on $\Delta$ we see that $\Delta^*$ is
$t$-torsionfree. As ${\rm
    dim}_k{\rm
    ker}[t:\Delta\longrightarrow\Delta]=n<\infty$ it follows that $\Delta^*$ is
  free of rank $n$.\hfill$\Box$\\

{\bf Remark:} Consider the category of $k[[t]][\varphi,\Gamma]$-modules which, as $k[[t]]$-modules, have the following properties: They are torsion modules, with $t$ acting surjectively and with finite dimensional (over $k$) kernel. Theorem \ref{semferiend} together with Lemma \ref{sur} describe an exact functor from this category to the category of \'{e}tale $(\varphi,\Gamma)$-module. As in this paper we are interested in objects of rank two only, this functoriality will play no role.\\

In the remainder of this section we provide explicit constructions of $\Delta$'s as in Theorem \ref{semferiend} for which the resulting \'{e}tale $(\varphi,\Gamma)$-module is irreducible of rank one or two (Lemma \ref{modvacuumonedim0} and Lemma \ref{fassageo}). The main objective in section 4 will be to provide a similar construction in the case where the resulting \'{e}tale $(\varphi,\Gamma)$-module is reducible of rank two.
 
\begin{lem}\label{fassa} (a) For $\alpha\in k^{\times}$ and $\ell, n\in{\mathbb Z}_{>0}$ define the $k[[t]][\varphi]$-submodule $\nabla$ of $M=k[[t]][\varphi]\otimes_kke=k[[t]][\varphi]e$ to be generated by the elements $t^ne$ and $t^{\ell}\varphi e-\alpha e$ and put $\Delta=M/\nabla$. Identify $e$ with its image element in $\Delta$, and inside
  $\Delta$ consider $$C=\sum_{s>0}\sum_{0\le\theta<q^{s-1}\ell}t^{\theta}\varphi^ske.$$If $\ell\le n{\xi}$ then we have the $k$-vector space decomposition\begin{gather}\Delta=C\bigoplus\bigoplus_{i=0}^{n-1}k.t^ie.\label{wassermuehle0georg0}\end{gather}The action of $t$ on $\Delta$ is surjective. The $k[[t]]$-module $\Delta^*$ is free of rank one, a basis is given by the element $\lambda$ defined by $\lambda(t^{n-1}e)=1$ and $\lambda(C\oplus\oplus_{i=0}^{n-2}k.t^ie)=0$.

(b) Let $\alpha, \beta\in k^{\times}$ and $N,\ell, \mu, n_x, n_y, m\in{\mathbb Z}_{>0}$ and $F\in t^mk[[t]]$ be given. Write $V=kx\oplus ky$, define the $k[[t]][\varphi]$-submodule $\nabla$ of $M=k[[t]][\varphi]\otimes_kV=k[[t]][\varphi]V$ to be generated by$$t^{n_x}x,\quad\quad\quad  t^{n_y}y,\quad\quad\quad  t^{\mu}\varphi y-\beta y,\quad\quad\quad t^{\ell}\varphi x-\alpha x+F\varphi^Ny$$and put $\Delta=M/\nabla$. Identify $V$ with its image in $\Delta$, and inside
$\Delta$ consider$$C=\sum_{s>0}\sum_{0\le\theta<q^{s-1}\ell}t^{\theta}\varphi^skx+\sum_{s>0}\sum_{0\le\theta<q^{s-1}\mu}t^{\theta}\varphi^sky.$$If $\ell\le n_x{\xi}$ and $\mu\le n_y{\xi}$ and \begin{gather}m+qn_x\ge \ell+\mu(1+q+\ldots+q^{N-1})+n_y\label{univorso}\end{gather}then we have the $k$-vector space decomposition\begin{gather}\Delta= C\bigoplus (\bigoplus_{i=0}^{n_x-1}k.t^ix)\bigoplus (\bigoplus_{i=0}^{n_y-1}k.t^iy).\label{wassermuehle0georg1}\end{gather}The action of $t$ on $\Delta$ is surjective. The $k[[t]]$-module $\Delta^*$ is free of rank two, a basis is given by the elements $\lambda_x$, $\lambda_y$ defined by $\lambda_x(t^{n_x-1}x)=\lambda_y(t^{n_y-1}y)=1$ and\begin{gather}\lambda_x(C)=\lambda_x(\oplus_{i=0}^{n_x-2}k.t^ix)=\lambda_x(\oplus_{i=0}^{n_y-1}k.t^iy)=0\notag\\\lambda_y(C)=\lambda_y(\oplus_{i=0}^{n_x-1}k.t^ix)=\lambda_y(\oplus_{i=0}^{n_y-2}k.t^iy)=0.\notag\end{gather}
\end{lem}

{\sc Proof:} (a) The decomposition (\ref{wassermuehle0georg0}) follows from repeated application of the relation $\varphi t=t^q\varphi$ in $k[[t]][\varphi]$ and the defining relations for $\nabla$. (It also results from the subsequent discussion). For $h>0$ define in $\Delta$ the element$$t^{-h}e:=\alpha^{-h}t^{(1+q+q^2+\ldots+q^{h-1})\ell-h}\varphi^he$$and on the other hand consider (the classes of) the elements $e, te,\ldots, t^{n-1}e$ in $\Delta$. The assumption $\ell\le n{\xi}$ implies $t^{n-1}e\notin \nabla$, so that $t^{n-1}e$ does not vanish in $\Delta$. One checks that $\{t^ie\}_{i\in{\mathbb Z}_{\le n-2}}$ is a $k$-basis of $C\oplus\oplus_{i=0}^{n-2}k.t^ie$, and that $t(t^ie)=t^{i+1}e$ (all $i\in{\mathbb Z}_{\le n-2}$) and $t^{n}e=0$ (in $\Delta$). In particular, $t$ acts surjectively on $\Delta$, with one-dimensional kernel. Therefore $\lambda$ is a basis element for the $k[[t]]$-module $\Delta^*$, cf. Lemma \ref{sur}.

(b) Letting $\overline{x}$ denote the image of $x$ in $M/k[[t]][\varphi]y$, we have$$M/k[[t]][\varphi]y\cong k[[t]][\varphi]\otimes_kk.\overline{x}=k[[t]][\varphi]\overline{x}.$$We first claim\begin{gather}\nabla\cap k[[t]][\varphi]y=\nabla_y\label{cnord}\end{gather}where $\nabla_y$ is the $k[[t]][\varphi]$-submodule generated by $t^{\mu}\varphi y-\beta y$ and $t^{n_y}y$. The inclusion $k[[t]][\varphi]y\cap\nabla\supset\nabla_y$ is obvious. As $\nabla$ is generated modulo $\nabla_y$ by $R=t^{\ell}\varphi x-\alpha x+F\varphi^Ny$ and $t^{n_x}x$, and as $k[[t]][\varphi]y\cap k[[t]][\varphi]t^{n_x}x=0$, to see the inclusion $k[[t]][\varphi]y\cap\nabla\subset\nabla_y$ it will be enough to see that for all
$B\in k[[t]][\varphi]$ with $BR\in k[[t]][\varphi]y+k[[t]][\varphi]t^{n_x}x$
we in fact have $BR\in\nabla_y+k[[t]][\varphi]t^{n_x}x$. Now $BR\in k[[t]][\varphi]y+k[[t]][\varphi]t^{n_x}x$ means $$B (t^{\ell}\varphi
\overline{x}-\alpha \overline{x})\in k[[t]][\varphi]t^{n_x}\overline{x}$$inside
$k[[t]][\varphi]\overline{x}$, i.e. $$B (t^{\ell}\varphi
-\alpha )\in k[[t]][\varphi]t^{n_x}$$inside
$k[[t]][\varphi]$. We may write
$B=\sum_{n=0}^{n_0}B_n\varphi^n$ for some $n_0\in{\mathbb
  Z}_{\ge0}$, some $B_n\in k[[t]]$. From $B (t^{\ell}\varphi
-\alpha )\in k[[t]][\varphi]t^{n_x}$ we obtain
$$B_{n_0}t^{\ell q^{n_0}}\varphi^{n_0+1}\in
k[[t]]\varphi^{n_0+1}t^{n_x}=k[[t]]t^{n_xq^{n_0+1}}\varphi^{n_0+1}$$and this
implies \begin{gather}B_{n_0}\in t^{n_xq^{n_0+1}-\ell q^{n_0}}k[[t]].\label{steuer0}\end{gather}Formula (\ref{univorso}) implies that $t^{n_xq-\ell}F\varphi^N y$ belongs to $\nabla_y$, hence that $t^{n_xq^{n_0+1}-\ell q^{n_0}}\varphi^{n_0}F\varphi^N y$ belongs to $\nabla_y$. Combined with formula (\ref{steuer0}) this implies
$B_{n_0}\varphi^{n_0}F\varphi^N y\in \nabla_y$. On the other hand, formula (\ref{steuer0}) also implies $B_{n_0}\varphi^{n_0}(t^{\ell}\varphi x-\alpha x)\in k[[t]][\varphi]t^{n_x}x$. Thus, we may replace $B$ by $B-B_{n_0}\varphi^{n_0}$. This makes the $\varphi$-degree decrease. Continuing inductively, we finally arrive
at $B=0$. Formula (\ref{cnord}) is proven.

Given formula (\ref{cnord}), it follows that the submodule $k[[t]][\varphi]y$ and its quotient $\Delta/k[[t]][\varphi]y\cong k[[t]][\varphi]\overline{x}$ are both of the sort looked at in statement (a). Since (by statement (a)) $t$ acts surjectively on both of them, it acts surjectively on $\Delta$. The decomposition (\ref{wassermuehle0georg1})
follows from the decomposition (\ref{wassermuehle0georg0}), applied to both
$k[[t]][\varphi]y$ and $k[[t]][\varphi]\overline{x}$. Similarly, formula
(\ref{cnord}) implies that the claim on the $k[[t]]$-basis of $\Delta^*$
follows from that in statement (a).\footnote{A $k$-basis of $\Delta$ on which
  $k[[t]]$ acts similarly as in (a) can be constructed as follows. For integers $h\ge0$ put $V_h=kx_h\oplus ky_h$ and define the $k[[t]][\varphi]$-submodule $\nabla_h$ of $$M_h=k[[t]][\varphi]\otimes_kV_h=k[[t]][\varphi]V_h$$ to be generated by $t^{h+n_x}x_h$, $t^{h+n_y}y_h$ and $$t^{{\xi}h+\mu}\varphi y_h-\beta y_h,\quad t^{{\xi}h+\ell}\varphi x_h-\alpha x_h+t^{{\xi}h}F\varphi y_h.$$We then have a well defined $k[[t]][\varphi]$-linear isomorphism$$\sigma_h:\Delta\longrightarrow M_h/\nabla_h,\quad\mbox{ determined by }\quad x\mapsto t^hx_h,\quad y\mapsto t^hy_h,$$as follows from (a) (applied twice, using formula (\ref{cnord})). Define the element $t^{-h}x:=\sigma_h^{-1}(x_h)$ and $t^{-h}y:=\sigma_h^{-1}(y_h)$ in $\Delta$. Do this for all $h>0$, and on the other hand consider the elements $x, tx,\ldots, t^{n_x-1}x$ and $y, ty,\ldots, t^{n_y-1}y$. Then $$\{t^ix\}_{i\in{\mathbb Z}_{\le n_x-2}}\cup\{t^iy\}_{i\in{\mathbb Z}_{\le n_y-2}}$$is a $k$-basis of $C\oplus(\oplus_{i=0}^{n_x-2}k.t^ix)\oplus(\oplus_{i=0}^{n_y-2}k.t^iy)$, and $t(t^ix)=t^{i+1}x$ (all $i\in{\mathbb Z}_{\le n_x-2}$) and $t(t^iy)=t^{i+1}y$ (all $i\in{\mathbb Z}_{\le n_y-2}$) and $t^{n_x}x=0$ and $t^{n_y}y=0$ (in $\Delta$).}\hfill$\Box$\\

{\bf Remark:} At first the definition of $\nabla$ in Lemma \ref{fassa} (b) might look somewhat unmotivated or random. However, the intent is to define a reducible rank two object whose rank one subquotients are as in Lemma \ref{fassa} (a). Thus, for both of these subquotients one should employ the two definining relations in \ref{fassa} (a), i.e. one should employ defining relations $t^{n_y}y$, $t^{\mu}\varphi y-\beta y$ for the subobject, and $t^{n_x}x$ and $t^{\ell}\varphi x-\alpha x$ for the quotient object. But in order to allow non-splitting extensions of these we add to $t^{\ell}\varphi x-\alpha x$ an expression $F\varphi^Ny$; this added summand $F\varphi^Ny$ then dictates the extension type (and is in a sense the main term looked at in this paper).\\

{\bf Remark:} Let us look at various parameters appearing in Lemma \ref{fassa} (b). Obviously, the choice of $m$ does not affect $\Delta$ at all (as long as formula (\ref{univorso}) holds true). The numbers $n_x$, $n_y$ are auxiliary, in the following sense: For another such pair $n_x'$, $n_y'$, with corresponding $\Delta'$, we have a surjection $\Delta'\to \Delta$ with finite dimensional (over $k$) kernel, inducing an isomorphism $\Delta^*\otimes_{k[[t]]}k((t))\cong (\Delta')^*\otimes_{k[[t]]}k((t))$. The dependence of $\Delta$ on $N$ and $F$ is given according to the following observation: Putting $N'=N+1$ and $F'=\beta^{-1}t^{\mu q^N}F$ leads to an isomorphic $\Delta$. This is due to the congruence$$F'\varphi^{N'}y=\beta^{-1}F t^{\mu q^N} \varphi^{N+1}y= F \varphi^{N}(\beta^{-1}t^{\mu})\varphi y\equiv F \varphi^{N}y$$modulo the sub module generated by $t^{\mu}\varphi y-\beta y$ (which belongs to our generators of $\Delta$).

\begin{lem}\label{modvacuumonedim0} Fix $\eta\in
  k^{\times}$ and $w\in{\mathbb Z}$ and $N\in{\mathbb Z}_{>0}$. In
  $M=k[[t]][\varphi]\otimes_kkz=k[[t]][\varphi]z$ define $\nabla$ to be the
  $k[[t]][\varphi]$-submodule generated by $t^{N}z$ and
  $t^{\xi}\varphi z-\eta z$. The $k[[t]][\varphi]$-action on
  $\Delta=M/\nabla$ extends in a unique way to a
  $k[[t]][\varphi,\Gamma]$-action such that
  $\gamma*z=a_{\gamma,1}^{w+1}\frac{t}{[\gamma]}z$ for all
  $\gamma\in\Gamma$, and $\Delta^*\otimes_{{k[[t]]}}{k((t))}$ becomes a $(\varphi,\Gamma)$-module of rank one. It does not depend on $N$. We denote it by ${\bf E}(\eta,w)$. There is a $k((t))$-basis element $\lambda$ of ${\bf E}(\eta,w)$ such that\begin{gather}\varphi(\lambda)=\eta^{-1}\lambda,\label{altneu0}\\\gamma*\lambda=a_{\gamma,1}^{-w}\lambda\label{altneu1}\end{gather}for all $\gamma\in\Gamma$ with $[\gamma]=a_{\gamma,1}t$, and $\lambda$ is uniquely determined
up to a $k^{\times}$-scalar by formula (\ref{altneu0}). The ${\mathcal G}_F$-character corresponding to ${\bf E}(\eta,w)$ is $\omega^{-w}\mu_{\eta^{-1}}$.   
\end{lem}
 
{\sc Proof:} That $\Delta^*$ is free of rank one over $k[[t]]$ follows from Lemma \ref{fassa}, hence Theorem \ref{semferiend} applies. Identify $z$ with its image element in $\Delta$. Inside $\Delta$ define $$C=\sum_{s>0}\sum_{0\le\theta<q^{s-1}{\xi}}t^{\theta}\varphi^sz.$$We then have
the direct sum decomposition\begin{gather}\Delta=kz\oplus(C+\sum_{i=1}^{N-1}k.t^iz).\label{wassermuehle00}\end{gather}We use it to define
$\lambda\in\Delta^*$ by\footnote{If $N\ge 2$ then $\lambda$ is not a $k[[t]]$-basis
  of $\Delta^*$.} requiring $$\lambda(C+\sum_{i=1}^{N-1}k.t^iz)=0,\quad\quad \quad\lambda(z)=1.$$To prove formula (\ref{altneu0}) we need to see that the image of $\lambda$ under the map (\ref{pruefeaus1}) is the same as the image
    of $\eta^{-1}\lambda$ under the map (\ref{pruefeaus0}). Of course, we may
    just as well multiply with $t^{\xi}$ and prove that the image of $t^{\xi} \lambda$ under the map (\ref{pruefeaus1}) is the same as the image
    of $t^{\xi}\eta^{-1}\lambda$ under the map (\ref{pruefeaus0}). Both these image
    elements can be seen as $k$-linear forms on
    $\Delta\otimes_{k[[t]],\varphi}k[[t]]$, and to show that they are the
    same, we show that they take the same value at $x\otimes t^{\kappa}$ for
    all $x\in\Delta$, all $0\le\kappa\le \xi$. This comes down to proving\begin{align}(\eta^{-1}\lambda)(t^{\kappa+\xi}\varphi x)=&\left\{\begin{array}{l@{\quad:\quad}l}  \lambda(x)& \kappa=0\\0&1\le \kappa\le \xi\end{array}\right..\label{hochter00}\end{align}It suffices to prove formula (\ref{hochter00}) for all $x=t^{\theta}
    \varphi^{\iota-1}z$ with $\theta\ge0$ and $\iota\ge1$,
    with the additional assumption $\theta< q^{\iota-2}{\xi}$ if
    $\iota>1$. This is because $\Delta$ is the $k$-span of such elements, by
    formula (\ref{wassermuehle00}). Now, as $\eta z=t^{\xi}\varphi z$ and
   $t^{\kappa+\xi}\varphi x=t^{\kappa+\xi+q\theta}\varphi^{\iota} z$ in $\Delta$, the claim
     follows from the definition of $\lambda$ (if $\theta=0$ and $\iota=1$ and
     $\kappa=0$ both sides are $=1$, otherwise both sides vanish). 

To prove formula (\ref{altneu1}) we compute
$(\gamma*\lambda)(z)=\lambda(\gamma^{-1}*z)=a_{\gamma,1}^{-w}$.

That a $k((t))$-basis vector $\lambda$ is uniquely determined
up to a $k^{\times}$-scalar by formula (\ref{altneu0}) follows from the fact
that $k$ is the space of invariants in $k[[t]]$ under $\varphi$.\hfill$\Box$\\

\begin{lem}\label{fassageo} Let $a, b\in k^{\times}$ and $\nu , \mu, n_x, n_y\in{\mathbb Z}_{>0}$ such that\begin{gather} \nu \le qn_x-n_y,\quad \mu\le qn_y-n_x.\label{univorso1}\end{gather}Write $V=kx\oplus ky$, define the $k[[t]][\varphi]$-submodule $\nabla$ of $M=k[[t]][\varphi]\otimes_kV=k[[t]][\varphi]V$ to be generated by$$t^{n_x}x,\quad\quad\quad t^{n_y}y,\quad\quad\quad t^{\mu}\varphi y- bx,\quad\quad\quad t^{\nu }\varphi x-ay$$and put $\Delta=M/\nabla$. Identify $V$ with its image in $\Delta$. 

(a) The $k[[t]]$-module $\Delta^*$ is free of rank two.

(b) Let $u\in{\mathbb Z}$. Assume that we have the modulo $p{\mathbb Z}$-congruences\footnote{If $\mu>\xi$ then we read condition (\ref{james0}) to mean $${\mu-\xi\choose m}\equiv 0\quad\mbox{ for all integers }1\le m<{\rm min}\{\frac{n_x}{\xi}, \mu+2-q\},$$and similarly for condition (\ref{james1}).}\begin{gather}{\xi-\mu\choose m}\equiv 0\quad\mbox{ for all integers }1\le m<{\rm min}\{\frac{n_x}{\xi}, q-\mu\},\label{james0}\\{\xi-\nu \choose m}\equiv 0\quad\mbox{ for all integers }1\le m<{\rm min}\{\frac{n_y}{\xi}, q-\nu\}.\label{james1}\end{gather} The $k[[t]][\varphi]$-action on $\Delta$ extends in a unique way to a $k[[t]][\varphi,\Gamma]$-action such that$$\gamma*x=a^{\mu+u+1}_{\gamma,1}\frac{t}{[\gamma]}x,\quad\quad\gamma*y=a^{u+1}_{\gamma,1}\frac{t}{[\gamma]}y$$for all $\gamma\in\Gamma$. In this way, $\Delta^*\otimes_{{k[[t]]}}{k((t))}$ becomes a $(\varphi,\Gamma)$-module of rank two, by the
  construction of Theorem \ref{semferiend}. It has a $k((t))$-basis $\lambda_{x}, \lambda_{y}$ for which $$\varphi(\lambda_{y})=a^{-1}t^{\mu-\xi}\lambda_{x},\quad\quad \varphi(\lambda_{x})=b^{-1}t^{\nu-\xi}\lambda_{y},$$$$\gamma*\lambda_{x}=a^{-\mu-{u}}_{\gamma,1}\lambda_{x},\quad\quad\gamma*\lambda_{y}=a^{-{u}}_{\gamma,1}\lambda_{y}$$for all $\gamma\in\Gamma$ with $[\gamma]=a_{\gamma,1}t$.
\end{lem}

 {\sc Proof:} (a) Inside $\Delta$ consider$$C=\sum_{s>0}\sum_{0\le\theta<q^{s-1}\nu }t^{\theta}\varphi^skx+\sum_{s>0}\sum_{0\le\theta<q^{s-1}\mu}t^{\theta}\varphi^sky.$$The assumptions (\ref{univorso1}) imply $t^{n_x-1}x\notin \nabla$ and $t^{n_y-1}y\notin \nabla$, so that (the classes of) $t^{n_x-1}x$ and $t^{n_y-1}y$ do not vanish in $\Delta$. Repeated application of the relation $\varphi t=t^q\varphi$ in $k[[t]][\varphi]$ and observing the defining relations for $\nabla$ yields the $k$-vector space decomposition\begin{gather}\Delta= C\bigoplus (\bigoplus_{i=0}^{n_x-1}k.t^ix)\bigoplus (\bigoplus_{i=0}^{n_y-1}k.t^iy).\label{wassermuehle0georg1nord}\end{gather}We use it to define the elements $\widetilde{\lambda}_x$, $\widetilde{\lambda}_y$ of $\Delta^*$ by demanding $\widetilde{\lambda}_x(t^{n_x-1}x)=\widetilde{\lambda}_y(t^{n_y-1}y)=1$ and $\widetilde{\lambda}_x(C\oplus(\oplus_{i=0}^{n_x-2}k.t^ix)\oplus(\oplus_{i=0}^{n_y-1}k.t^iy))=0$ and $\widetilde{\lambda}_y(C\oplus(\oplus_{i=0}^{n_x-1}k.t^ix)\oplus(\oplus_{i=0}^{n_y-2}k.t^iy))=0$. We claim that these form a $k[[t]]$-basis of $\Delta^*$. To see this, define for $h>0$ in $\Delta$ the elements$$t^{-h}x:=(ab)^{-h}t^{(1+q^2+q^4+\ldots+q^{2h-2})(\mu+q\nu)-h}\varphi^hx,$$$$t^{-h}y:=(ab)^{-h}t^{(1+q^2+q^4+\ldots+q^{2h-2})(\nu+q\mu)-h}\varphi^hy,$$and on the other hand consider its elements $x, tx,\ldots, t^{n_x-1}x$ and $y, ty,\ldots, t^{n_y-1}y$. Then $$\{t^ix\}_{i\in{\mathbb Z}_{\le n_x-2}}\cup\{t^iy\}_{i\in{\mathbb Z}_{\le n_y-2}}$$is a $k$-basis of $C\oplus(\oplus_{i=0}^{n_x-2}k.t^ix)\oplus(\oplus_{i=0}^{n_y-2}k.t^iy)$, and both $t(t^ix)=t^{i+1}x$ and $t(t^iy)=t^{i+1}y$ for all $i\in{\mathbb Z}_{\le n_y-2}$, as well as $t^{n_x}x=0$ and $t^{n_y}y=0$ (in $\Delta$). This implies the claim. 
 
 (b) The assumption (\ref{james0}) and formula (\ref{amrit}) imply $(\frac{t}{[\gamma]})^{\xi-\mu}t^{\mu}\varphi y\equiv a^{\mu}_{\gamma,1}t^{\mu}\varphi y$ modulo $\nabla$, and thus $$\gamma*(t^{\mu}\varphi y-bx)=a_{\gamma,1}^{u+1}\frac{t}{[\gamma]}((\frac{t}{[\gamma]})^{\xi-\mu}t^{\mu}\varphi y-a^{\mu}_{\gamma,1}bx)\in\nabla.$$ Similarly, the assumption (\ref{james1}) and formula (\ref{amrit}) imply $\gamma*(t^{\nu}\varphi x-ay)\in\nabla$. To see that the stated formulae define a group action of $\Gamma$ on
$\Delta$, and hence a $k[[t]][\varphi,\Gamma]$-action
on $\Delta$, we take $\gamma,
\gamma'\in \Gamma$ and compute$$\gamma'*(\gamma*y)=a_{\gamma',1}^{u+1}a_{\gamma,1}^{u+1}\frac{[\gamma']}{[\gamma\gamma']}\frac{t}{[\gamma']}y=a_{\gamma'\gamma,1}^{u+1}\frac{t}{[\gamma\gamma']}y=(\gamma'\gamma)*y$$and similarly $\gamma'*(\gamma*x)=(\gamma'\gamma)*x$.

Next, we use the decomposition (\ref{wassermuehle0georg1nord}) to define $\lambda_{x}, \lambda_{y}\in\Delta^*$ by demanding$$\lambda_x(C\bigoplus (\bigoplus_{i=1}^{n_x-1}k.t^ix)\bigoplus (\bigoplus_{i=0}^{n_y-1}k.t^iy))=0,$$$$\lambda_y(C\bigoplus (\bigoplus_{i=0}^{n_x-1}k.t^ix)\bigoplus (\bigoplus_{i=1}^{n_y-1}k.t^iy))=0,$$$$\lambda_x(x)=\lambda_y(y)=1.$$To prove $\varphi(\lambda_{y})=a^{-1}t^{\mu-\xi}\lambda_{x}$ we need to see that the image of $\lambda_{y}$ under the map (\ref{pruefeaus1}) is the same as the image
    of $a^{-1}t^{\mu-\xi}\lambda_{x}$ under the map (\ref{pruefeaus0}). Of course, we may
    just as well multiply with $t^{\xi}$ and prove that the image of $t^{\xi}\lambda_{y}$ under the map (\ref{pruefeaus1}) is the same as the image
    of $a^{-1}t^{\mu}\lambda_{x}$ under the map (\ref{pruefeaus0}). Both these image
    elements can be seen as $k$-linear forms on
    $\Delta\otimes_{k[[t]],\varphi}k[[t]]$, and to show that they are the
    same, we show that they take the same value at $z\otimes t^{\kappa}$ for
    all $z\in\Delta$, all $0\le\kappa\le \xi$. This comes down to proving\begin{align}(a^{-1}\lambda_x)(t^{\kappa+\mu}\varphi z)=&\left\{\begin{array}{l@{\quad:\quad}l}  \lambda_y(z)& \kappa=0\\0&1\le\kappa\le \xi\end{array}\right..\label{hochter00marsch}\end{align}It suffices to prove formula (\ref{hochter00marsch}) for all $z=t^{\theta}
    \varphi^{\iota-1}x$ with $\theta\ge0$ and $\iota\ge1$,
    with the additional assumption $\theta< q^{\iota-2}\nu$ if
    $\iota>1$, resp. for all $z=t^{\theta}\varphi^{\iota-1}y$ with $\theta\ge0$ and $\iota\ge1$,
    with the additional assumption $\theta< q^{\iota-2}\mu$ if
    $\iota>1$. This is because $\Delta$ is the $k$-span of such elements, by
    formula (\ref{wassermuehle0georg1nord}). The claim
     follows immediately from the definitions (if $z=y$ and
     $\kappa=0$ both sides are $=1$, otherwise both sides vanish).
     
    The proof of $\varphi(\lambda_{x})=b^{-1}t^{\nu-\xi}\lambda_{y}$ runs the same way. 
 
 To prove $\gamma*\lambda_{x}=a^{-{u}}_{\gamma,1}\lambda_{x}$ for $\gamma\in\Gamma$ with $[\gamma]=a_{\gamma,1}t$ we compute $(\gamma*\lambda_{x})(x)=\lambda_{x}(\gamma^{-1}*x)=a_{\gamma,1}^{-u}$. The proof of $\gamma*\lambda_{y}=a^{-{u}}_{\gamma,1}\lambda_{y}$ is the same.\hfill$\Box$\\

\section{$(\varphi,\Gamma)$-modules for reducible two dimensional Galois representations}

\subsection{The case $0<\ell<q-1$}

\label{viervor}

Let us try to motivate the long list of integers and of congruences satisfied by them we are starting this subsection with (cf. also our discussions of formulae (\ref{revf2}) and (\ref{mariaeverkuendigung})). We will soon be facing the following situation. A $k[[t]][\varphi]$-module $M$ will be given to us, with a semilinear action by $\Gamma$. We want to define a $k[[t]][\varphi]$-submodule $\nabla$ of $M$ by giving explicit generators, so that $\nabla$ is stable under $\Gamma$. Writing the said generators as $\sum_{m\ge0}t^mz_m$ with {\it certain} $z_m\in M$, we need to see $\gamma*\sum_{m\ge0}t^mz_m\in\nabla$ for $\gamma\in\Gamma$. But $\gamma*\sum_{m\ge0}t^mz_m=\sum_{m\ge0}[\gamma]^m\gamma*z_m$, which requires us to understand the series $[\gamma]^m\in k[[t]]$. And indeed, the key for seeing $\sum_{m\ge0}[\gamma]^m\gamma*z_m\in \nabla$ will be to see that for certain "good" exponents $m$ (only these $m$'s should then occur in the said generators), the series $[\gamma]^m$ has many vanishing coefficients. These vanishings are precisely implied by the vanishing (modulo $p{\mathbb Z}$) formulae of Lemma \ref{maryannun}.\\

For $1\le j\le f$ put$$g^{(j)}={\xi}(\sum_{s=0}^{j-1}p^{2f+s+1-j}-p^{f+1-j}).$$For $m,i\in{\mathbb Z}$ with $i\ge0$ put$$H(m,i)=1-(m+2)\sum_{j=1}^ft^{p^ig^{(j)}}\in k[t].$$Fix an integer $0<\ell<\xi$. Write $\ell=\sum_{j=0}^{f-1}m_jp^j$ with
$0\le m_j<p$, write $m_{j+f}=m_j$ for all $j\in{\mathbb Z}$. For $i\in{\mathbb Z}_{\ge0}$ put\begin{gather}\ell_i=\ell_{[i]}=\sum_{j=0}^{f-1}m_{j+i}p^{j+i}.\label{movopfin}\end{gather}Fix a natural number $N$. Put $$r=q^{{N}}-1-\frac{p^{f-1}\xi}{p-1}$$and for $0\le i< f$
put$$o_i=q^{{N}}-1-p^i(\xi+q\frac{\xi}{p-1}).$$Let ${\mathcal D}$ and ${\mathcal E}$ be finite subsets of ${\mathbb
  Z}_{\ge0}$. For $i\in{\mathcal D}$ put$$r_i=q^{{N}}-1+\ell_i-p^i{\xi}.$$For $i\in{\mathcal E}$ and $1\le
j\le f$ put\begin{align}{{n}}_i^{(0)}&=q^N-1+q^2\ell_i-2p^{2f+i}{\xi},\notag\\{{n}}_i^{(j)}&=q^N-1+q^2\ell_i-p^i{\xi}(q^2+p^{f+1-j}-\sum_{s=0}^{j-2}p^{2f+s+1-j})\notag\\{}&={n}_i^{(0)}+p^ig^{(j)}.\notag\end{align}

Assume that the numbers $r$, $o_i$, $r_i$ and ${{n}}_i^{(j)}$ are all non-negative; as ${\mathcal D}$ and ${\mathcal E}$ are finite, this may be achieved by suitably increasing $N$ if necessary.\footnote{The role of $N$ will only be auxiliary.} 

\begin{lem}\label{maryannun} We have
  ${{n}}_i^{(j)}<{{n}}_i^{(j')}$ for all $0\le j<j'\le f$,
  and\begin{gather}{{n}}_i^{(1)}-{{n}}_i^{(0)}=qp^i{\xi}^2,\label{wiederlegaut}\\{{n}}_i^{(j)}+p^{i+1-j+f}{\xi}(p^{f-1}+1)={{n}}_i^{(j+1)}+p^{i-j+f}{\xi}\quad\mbox{
      for }1\le j<f.\label{fuhleresel000}\end{gather}We have the modulo $p{\mathbb
     Z}$-congruences\begin{gather}{r+1\choose m}\equiv0\quad\mbox{ for }m\in[1,p^{f-1}-1].\label{papa000}\end{gather}For $0\le i<f$ we have the modulo $p{\mathbb
     Z}$-congruences\begin{gather}{o_i+1\choose p^i}\equiv 1\quad\mbox{ and }\quad {o_i+1\choose m}\equiv0\quad\mbox{ for }m\in[1,p^iq-1]-\{p^i\}.\label{voroanti}\end{gather}For $i\in{\mathcal D}$ we have the modulo $p{\mathbb
    Z}$-congruences\begin{gather}{r_i+1\choose m}\equiv 0\quad\mbox{ for }m\in[1,p^i-1].\label{papa00}\end{gather}For  $i\in{\mathcal E}$ and $1\le j\le f$ we have the modulo $p{\mathbb
    Z}$-congruences\begin{gather}{{{n}}_i^{(0)}+1\choose p^{i+2f}}\equiv
  m_i+2\quad\mbox{ and }\quad {{{n}}_i^{(0)}+1\choose
    m}\equiv0\quad\mbox{ for
  }m\in[1,2p^{i+2f}-1]-\{p^{i+2f}\},\label{papa201}\\{{{n}}_i^{(j)}+1\choose
    p^{i+1-j+f}}\equiv 1\quad\mbox{ and }\quad {{{n}}_i^{(j)}+1\choose
    m}\equiv0\quad\mbox{ for
  }m\in[1,p^{i+1-j+2f}-1]-\{p^{i+1-j+f}\}.\label{papa200}\end{gather}For $i\in{\mathcal E}$ and $0\le j\le f$ we have the modulo ${\xi}{\mathbb
    Z}$-congruences \begin{gather}\ell\equiv {{n}}_i^{(j)}.\label{soju0}\end{gather}For $i\in{\mathcal D}$ we have the modulo ${\xi}{\mathbb
    Z}$-congruences \begin{gather}\ell\equiv r_i.\label{soju1}\end{gather}If $\ell=\frac{(p-2){\xi}}{p-1}$ we
  have for all $0\le i< f$ the modulo ${\xi}{\mathbb
    Z}$-congruences\begin{gather}\ell\equiv r\equiv o_i.\label{soju2}\end{gather}\end{lem}

 {\sc Proof:} For formula (\ref{papa00}) use $0\equiv r_i+1$ modulo
$p^i{\mathbb Z}$. For formula (\ref{papa000}) use $0\equiv r+1$
modulo $p^{f-1}{\mathbb Z}$. For formula (\ref{papa201}) use that $p^{i+2f}$ divides
${{n}}_i^{(0)}+1$ and $\frac{{{n}}_i^{(0)}+1}{p^{i+2f}}\equiv m_i+2$ modulo
$p{\mathbb Z}$. For formula (\ref{papa200}) use that $p^{i+1-j+f}$ divides
${{n}}_i^{(j)}+1$ and $\frac{{{n}}_i^{(j)}+1}{p^{i+1-j+f}}\equiv 1$ modulo $p{\mathbb Z}$, and that $p^{i+1-j+2f}$ divides ${{n}}_i^{(j)}+1-p^{i+1-j+f}$. Formula (\ref{voroanti}) is proven similarly. Formulae (\ref{soju0}), (\ref{soju1}) and (\ref{soju2}) follow from $\ell\equiv\ell_i$ modulo ${\xi}{\mathbb
    Z}$.\hfill$\Box$\\

Besides $\ell$ fix an integer $0\le u<\xi$ and scalars $\alpha,\beta\in
k^{\times}$. Fix scalars ${\bf b}_i\in k$ for all $i\in{\mathcal D}$, scalars
${\bf c}_i\in k$ for all $i\in{\mathcal E}$, and a scalar ${\bf c}\in
k$. However, if $$\ell\ne\frac{(p-2){\xi}}{p-1}\quad\quad\mbox{ or }\quad\quad
\alpha\ne \beta$$ then we demand ${\bf
  c}=0$ (so that ${\bf c}\ne0$ is allowed only if
$(\ell,\frac{\alpha}{\beta})=(\frac{(p-2){\xi}}{p-1},1)$).\footnote{Notice
  that if $\ell=\frac{(p-2){\xi}}{p-1}$ then $p\ne2$, as we demand $\ell>0$.}

Let $n_y$, $n_x$ be positive integers and assume\begin{gather}\ell_i\ge n_y\quad\mbox{
    for }i\in{\mathcal D},\label{sommerdornwald5}\\qn_x+\ell_i-{\xi}(p^i+1)\ge n_y\quad\mbox{ for }i\in{\mathcal
    D},\label{marufmin00}\\q^2\ell_i-p^{i+f}{\xi}(q-\sum_{s=0}^{f-1}p^{s})\ge n_y\quad\mbox{ for }i\in{\mathcal E},\label{fuhleresel200}\\qn_x+ q^2\ell_i-{\xi}(2p^{2f+i}+1)\ge n_y\quad\mbox{ for }i\in{\mathcal
    E},\label{maruf200}\\\xi\ge n_y\quad\mbox{ if }{\bf c}\ne0,\label{sommerdorne}\\p-2 \ge n_y\quad\mbox{ if }{\bf c}\ne0\mbox{ and }q=p,\label{sommerdornwald1}\\n_x-\frac{p^{f-1}\xi}{p-1}\ge n_y\quad\mbox{ if }{\bf c}\ne0,\label{sommerdornwald6}\\p(n_x-3)+2\ge n_y\quad\mbox{ if }{\bf c}\ne0\mbox{ and }q=p.\label{maruf2000}\end{gather}Let $\gamma\in\Gamma$. recall the coefficients $a_{\gamma,i}$ of $[\gamma]$, cf. Lemma \ref{doju}. In $k[[t]]$ define$$G_{\gamma}=G_{\gamma}(t)=a_{\gamma,q}^{p^{f-1}}a_{\gamma,1}^{-p^{f-1}-r}[\gamma]^r.$$Let $V$ be the two dimensional $k$-vector space on the basis $x,y$. Inside $M=k[[t]][\varphi]\otimes_kV=k[[t]][\varphi]V$ denote by $\nabla_y$ the $k[[t]][\varphi]$-submodule generated by $t^{n_y}y$ and $t^{\xi}\varphi y-\beta y$. Denote by $\nabla$ the
$k[[t]][\varphi]$-submodule of $M$ generated by $\nabla_y$ and $t^{n_x}x$ and\begin{align}R=&t^{\xi}\varphi x-\alpha x+\alpha\beta^{1-N}(\sum_{i\in{\mathcal D}}{\bf b}_it^{r_i}+{\bf c}\sum_{i=0}^{f-1}t^{o_i}+\sum_{i\in{\mathcal E}}{\bf c}_it^{{n}_i^{(0)}}H(m_i,i))\varphi^{{N}}y.\notag\end{align}For later reference we remark that $H(m_i,i)=H(p^{-i}\ell_{[i]},i)$ for $i\in{\mathcal E}$. Put $$\Delta=M/\nabla.$$

{\bf Remark:} The definition of the element $R$ is really the heart of the present section. Except for the scalar $\beta$ (the numbers  $N$, $n_x$ and $n_y$ are only auxiliary), $R$ encapsulates the entire datum of the \'{e}tale $(\varphi,\Gamma)$-module $\Delta^*\otimes_{{k[[t]]}}{k((t))}$ to be defined. The complexity in the definition of the exponents $r_i$, $o_i$, $n_i^{(0)}$ and those appearing in the $H(m_i,i)$ is forced on us by the requirement to have the key Proposition \ref{vacet10} functioning, i.e. it reflects the obstruction of extending the $k[[t]][\varphi]$-action to a $k[[t]][\varphi,\Gamma]$-action.

To cover all {\it isomorphism classes} of reducible \'{e}tale $(\varphi,\Gamma)$-module of rank two it would actually be sufficient to take ${\mathcal D}={\mathcal D}(\ell)$ and ${\mathcal E}={\mathcal E}(\ell)$, with ${\mathcal D}(\ell)$ and ${\mathcal E}(\ell)$ as in formula (\ref{regenyuval}). However, for later purposes it is important to allow more general ${\mathcal D}$, ${\mathcal E}$.\\

{\bf Remark:} It might be helpful (or consoling) to recall that the generic case is where ${\mathcal D}(\ell)={\mathcal F}$, ${\mathcal E}(\ell)=\emptyset$, ${\bf c}=0$. In this case the definition of $R$ greatly simplifies to become$$R=t^{\xi}\varphi x-\alpha x+\alpha\beta^{1-N}\sum_{i\in{\mathcal F}}{\bf b}_it^{r_i}\varphi^{{N}}y$$and the major part of the (rather ramified) subsequent subdiscussions evaporates.

\begin{pro}\label{vacet10}  The $k[[t]][\varphi]$-action on $\Delta$ extends in a unique way to a $k[[t]][\varphi,\Gamma]$-action satisfying the following formulae for all $\gamma\in\Gamma$:\begin{align}\gamma*x&=a_{\gamma,1}^{u}(\frac{a_{\gamma,1}t}{[\gamma]}x-{\bf c}\beta^{1-N}G_{\gamma}\varphi^{{N}} y),\notag\\\gamma*y&=a_{\gamma,1}^{u-\ell+1}\frac{t}{[\gamma]}y.\notag\end{align}
\end{pro}

{\sc Proof:} {\bf (A)} It is easy to see $\gamma*\nabla_y=\nabla_y$ for all $\gamma\in\Gamma$. Next we claim that the stated formulae define a group action of $\Gamma$ on
$M/\nabla_y$, and hence a $k[[t]][\varphi,\Gamma]$-action
on $M/\nabla_y$.

To see this, let $\gamma,\gamma'\in \Gamma$. We compute$$\gamma'*(\gamma*y)=a_{\gamma',1}^{u-\ell+1}a_{\gamma,1}^{u-\ell+1}\frac{[\gamma']}{[\gamma\gamma']}\frac{t}{[\gamma']}y=a_{\gamma'\gamma,1}^{u-\ell+1}\frac{t}{[\gamma\gamma']}y=(\gamma'\gamma)*y$$and\begin{align}\gamma'*(\gamma*x)&=a_{\gamma,1}^{u}\frac{a_{\gamma,1}[\gamma']}{[\gamma\gamma']}a_{\gamma',1}^{u}(\frac{a_{\gamma',1}t}{[\gamma']}x-{\bf
  c}\beta^{1-N}a_{\gamma',q}^{p^{f-1}}a_{\gamma',1}^{-p^{f-1}-r}[\gamma']^r\varphi^{{N}}y)\notag\\{}&\quad\quad - {\bf
  c} a_{\gamma,1}^{u}\beta^{1-N}a_{\gamma',1}^{u-\ell+1}a_{\gamma,q}^{p^{f-1}}a_{\gamma,1}^{-p^{f-1}-r}[\gamma\gamma']^r\varphi^{{N}}\frac{t}{[\gamma']}y\notag\\{}&=a_{\gamma\gamma',1}^{u}(\frac{a_{\gamma\gamma',1}t}{[\gamma\gamma']}x-{\bf c}\beta^{1-N}a_{\gamma',q}^{p^{f-1}}a_{\gamma',1}^{-p^{f-1}-r}a_{\gamma,1}\frac{[\gamma']^{r+1}}{[\gamma\gamma']}\varphi^{{N}}y)\notag\\{}&\quad\quad -a_{\gamma\gamma',1}^{u}a_{\gamma',1}^{1-\ell}a_{\gamma,q}^{p^{f-1}}a_{\gamma,1}^{-p^{f-1}-r}[\gamma\gamma']^{r}\varphi^{{N}}\frac{t}{[\gamma']}y.\notag\end{align}We
claim that this shows $\gamma'*(\gamma*x)\equiv(\gamma'\gamma)*x$ modulo
$\nabla_y$. Indeed, if ${\bf c}=0$ this is clear, so we may now assume
${\bf c}\ne0$. By formula (\ref{sommerdorne}) this implies $a_{\gamma',1}\frac{t}{[\gamma']}y\equiv y$ modulo $\nabla_y$. Thus, we need to see$$(a_{\gamma',q}^{p^{f-1}}a_{\gamma',1}^{-p^{f-1}-r}a_{\gamma,1}\frac{[\gamma']^{r+1}}{[\gamma\gamma']}+a_{\gamma',1}^{-\ell}a_{\gamma,q}^{p^{f-1}}a_{\gamma,1}^{-p^{f-1}-r}[\gamma\gamma']^{r})\varphi^{{N}}y\equiv a_{\gamma\gamma',q}^{p^{f-1}}a_{\gamma\gamma',1}^{-p^{f-1}-r}[\gamma\gamma']^r\varphi^{{N}}y$$modulo $\nabla_y$. In view of the decompositions (\ref{floos}) of $\gamma$ and $\gamma'$, it is enough to consider special cases for $\gamma$, $\gamma'$, as follows. Also, recall that $a_{\gamma\gamma',1}=a_{\gamma',1}a_{\gamma,1}$.

If $[\gamma]=a_{\gamma,1}t$ we use that $a_{\gamma,q}=0$ and $a_{\gamma\gamma',q}=a_{\gamma',q}a_{\gamma,1}$ and $[\gamma\gamma']=[\gamma']a_{\gamma,1}$. If $[\gamma']=a_{\gamma',1}t$ we use that $a_{\gamma',q}=0$ and $a_{\gamma\gamma',q}=a_{\gamma,q}a_{\gamma',1}$, as well as $\ell\equiv r$ modulo ${\xi}{\mathbb Z}$ (formula (\ref{soju2})). Finally, if $a_{\gamma,1}=a_{\gamma',1}=1$ we use formulae (\ref{freivorquinqua1}) and the following: For
all $\gamma\in\Gamma$ with $a_{\gamma,1}=1$ we have\begin{gather}t^r\varphi^{{N}}y\equiv\frac{[\gamma]^{r+1}}{t}\varphi^{{N}}y \quad\mbox{ modulo }\nabla_y.\label{florian}\end{gather}To prove formula (\ref{florian}), we recall $r=q^N-1-\frac{p^{f-1}\xi}{p-1}$ and first observe\begin{align}\frac{[\gamma]^{r+1}}{t}  &= t^{-1}(t+a_{\gamma,q}t^q+t^{(p^{f-1}+1)\xi+1}F)^{r+1}\quad\mbox{ for some }F\in k[[t]]\notag\\{}&\equiv t^r+\sum_{1\le m\le r+1}{r+1\choose m}a_{\gamma,q}^mt^{r+m\xi} \notag\\{}&= t^r+\sum_{1\le m\le r+1\atop m\ge p^{f-1}}{r+1\choose m}a_{\gamma,q}^mt^{q^N-1-\frac{p^{f-1}\xi}{p-1}+m\xi} \notag\end{align}modulo $t^{(p^{f-1}+1)\xi+q^N-1-\frac{p^{f-1}\xi}{p-1}}k[[t]]$, where the first equality follows from formula (\ref{amrit}) and the last equality follows from formula (\ref{papa000}). On the other hand, $t^{\xi}\varphi y-\beta y\in\nabla_y$ implies $t^{q^N-1}\varphi^N y-\beta^N y\in\nabla_y$. Now formula (\ref{florian}) follows from formula (\ref{sommerdorne}) (resp. (\ref{sommerdornwald1}) if $p=q$).

{\bf (B)} The $k[[t]][\varphi,\Gamma]$-action on $M/\nabla_y$ preserves $\nabla/\nabla_y$.

To see this, let $\gamma\in \Gamma$. That $\nabla/\nabla_y$ is preserved under
the action of $\gamma$ follows from formulae (\ref{soju0}), (\ref{soju1}),
(\ref{soju2}) and (\ref{freivorquinqua1}) if $[\gamma]=a_{\gamma,1}t$. Thus, as in {\bf (A)}, we may now assume
$a_{\gamma,1}=1$. Formula (\ref{sommerdornwald6}) implies $[\gamma]^{n_x}G_{\gamma}\varphi^Ny\in\nabla$ and hence $\gamma*t^{n_x}x\in\nabla$. It remains to see
$\gamma*R\in\nabla$. As $R\in
\nabla$ it suffices to see $\gamma*R-\frac{t}{[\gamma]}R\in\nabla$.

Let us write $H(m_i,i)|_{[\gamma]}$ for the result of substituting $[\gamma]$ for $t$ in $H(m_i,i)$, i.e. $H(m_i,i)|_{[\gamma]}=1-(m_i+2)\sum_{j=1}^f[\gamma]^{p^ig^{(j)}}$. Modulo $\nabla$ we
compute\begin{align}\gamma*(t^{\xi}\varphi
x)&\equiv[\gamma]^{\xi}\varphi(\frac{t}{[\gamma]}x-{\bf
  c}\beta^{1-N}G_{\gamma}\varphi^{{N}}y)\notag\\{}&{\equiv}\frac{t}{[\gamma]}t^{\xi}\varphi
x-[\gamma]^{\xi}\varphi {\bf
  c}\beta^{1-N}G_{\gamma}\varphi^{{N}}y,\label{wienschneider00}\\\gamma*(-\alpha
x)&\equiv-\alpha(\frac{t}{[\gamma]}x-{\bf
  c}\beta^{1-N}G_{\gamma}\varphi^{{N}}y),\label{wienschneider10}\\\gamma*(R-t^{\xi}\varphi
x+\alpha x)&\equiv\alpha\beta^{1-N}(\sum_{i\in{\mathcal D}}{\bf b}_i[\gamma]^{r_i}+{\bf
  c}\sum_{i=0}^{f-1}[\gamma]^{o_i}+\sum_{i\in {\mathcal
             E}}[\gamma]^{n_i^{(0)}}H(m_i,i)|_{[\gamma]})\varphi^{{N}}y,\label{wienschneider20}\end{align}\begin{align}-\frac{t}{[\gamma]}R\equiv&-\frac{t}{[\gamma]}t^{\xi}\varphi
  x+\alpha\frac{t}{[\gamma]}x-\alpha\beta^{1-N}\frac{t}{[\gamma]}(\sum_{i\in{\mathcal D}}{\bf
    b}_it^{r_i}+{\bf c}\sum_{i=0}^{f-1}t^{o_i}+\sum_{i\in {\mathcal
      E}}{\bf c}_it^{n_i^{(0)}}H(m_i,i))\varphi^{{N}}y.\label{wienschneider30}\end{align}Adding
formulae (\ref{wienschneider00}), (\ref{wienschneider10}), (\ref{wienschneider20}) and (\ref{wienschneider30}) we obtain that $\gamma*R-\frac{t}{[\gamma]}R$ is congruent modulo $\nabla$ to the sum\begin{gather}\sum_{i\in{\mathcal D}}{\bf b}_i\alpha\beta^{1-N}([\gamma]^{r_i}-\frac{t}{[\gamma]}t^{r_i})\varphi^{{N}}y\label{savorad00000}\\+\sum_{i\in {\mathcal E}}{\bf c}_i\alpha\beta^{1-N}([\gamma]^{{{n}}_i^{(0)}}-\frac{t}{[\gamma]}t^{{{n}}_i^{(0)}}-(m_i+2)\sum_{j=1}^f([\gamma]^{{{n}}_i^{(j)}}-\frac{t}{[\gamma]}t^{{{n}}_i^{(j)}}))\varphi^{{N}}y\label{savorad0000}\\+{\bf c}\beta^{1-N}(\alpha G_{\gamma}-[\gamma]^{\xi}\varphi G_{\gamma}+\alpha\sum_{i=0}^{f-1}([\gamma]^{o_i}-\frac{t}{[\gamma]}t^{o_i}))\varphi^{{N}}y.\label{savorad10}\end{gather}

It is therefore enough to show:

(a) For each $i\in{\mathcal D}$ the corresponding summand in (\ref{savorad00000}) belongs to $\nabla$.

(b) For each $i\in {\mathcal E}$ the corresponding summand in (\ref{savorad0000}) belongs to $\nabla$.

(c) the sum (\ref{savorad10}) belongs to $\nabla$.

To see (a) it is enough to
prove$$\frac{[\gamma]^{r_i+1}}{t}\varphi^{{N}}y\equiv
t^{r_i}\varphi^{{N}}y\quad\mbox{ modulo }\nabla.$$But this follows from formulae (\ref{sommerdornwald5}) and (\ref{amrit}) and (\ref{papa00}), similarly to the proof of formula (\ref{florian}).

We turn to (b). Multiplying
  with the element $t^{-q^{{N}}}[\gamma]$ of $t^{1-q^{{N}}}+t^{2-q^{{N}}}k[[t]]$, we see
  that it will be enough to
  show\begin{gather}t^{-q^N}([\gamma]^{{{n}}_i^{(0)}+1}  -{t}^{{{n}}_i^{(0)}+1})\equiv (m_i+2)t^{-q^N}\sum_{j=1}^f{t}^{{{n}}_i^{(j)}+1}-[\gamma]^{{{n}}_i^{(j)}+1}\label{zweivors00}\end{gather}modulo $t^{n_y}k[[t]]$. It follows from formulae (\ref{amrit}) and (\ref{papa201}) and (\ref{fuhleresel200} that $$t^{-q^N}([\gamma]^{{{n}}_i^{(0)}+1}-{t}^{{{n}}_i^{(0)}+1})\equiv (m_i+2)a_{\gamma,q}^{p^i}t^{q^2\ell_i-p^{i+2f}{\xi}}$$modulo $t^{n_y}k[[t]]$. On the other hand, for $1\le j\le f$ we have\begin{align}t^{-q^{{N}}}([\gamma]^{{{n}}_i^{(j)}+1}-t^{{{n}}_i^{(j)}+1})\equiv&
  a_{\gamma,q}^{p^{i+1-j}}t^{{{n}}_i^{(j)}+1-q^{{N}}+p^{i+1+f-j}{\xi}}-a_{\gamma,q}^{p^{i-j}}t^{{{n}}_i^{(j)}+1-q^{{N}}+p^{i+1+f-j}{\xi}(p^{f-1}+1)}\label{mami1000}\end{align}modulo
$t^{n_y}k[[t]]$, as follows from formulae (\ref{amrit00}), (\ref{fuhleresel200}) and (\ref{papa200}), similarly to the proof of formula (\ref{florian}). Thus, \begin{align}t^{-q^{{N}}}\sum_{j=1}^{f}
   [\gamma]^{{{n}}_i^{(j)}+1}-t^{{{n}}_i^{(j)}+1}
   &\stackrel{(i)}{\equiv}\sum_{j=1}^{f}a_{\gamma,q}^{p^{i+1-j}}t^{{{n}}_i^{(j)}+1-q^{{N}}+p^{i+1+f-j}{\xi}}-a_{\gamma,q}^{p^{i-j}}t^{{{n}}_i^{(j)}+1-q^{{N}}+p^{i+1+f-j}{\xi}(p^{f-1}+1)}\notag\\{}&\stackrel{(ii)}{\equiv}a_{\gamma,q}^{p^{i+1-f}}t^{{{n}}_i^{(f)}+1-q^{{N}}+p^{i+1}{\xi}}\notag\\{}&\quad\quad+\sum_{j=1}^{f-1}a_{\gamma,q}^{p^{i+1-j}}t^{{{n}}_i^{(j)}+1-q^{{N}}+p^{i+1+f-j}{\xi}}-a_{\gamma,q}^{p^{i-j}}t^{{{n}}_i^{(j+1)}+1-q^{{N}}+p^{i+f-j}{\xi}}\notag\\{}&\stackrel{(iii)}{\equiv}a_{\gamma,q}^{p^{i}}t^{{{n}}_i^{(1)}+1-q^{{N}}+p^{f+i}{\xi}}\notag\\{}&{\equiv}a_{\gamma,q}^{p^i}t^{q^2\ell_i-p^{i+2f}{\xi}}\notag \end{align}where
 in $(i)$ we used formula (\ref{mami1000}), in $(ii)$ we used formulae (\ref{fuhleresel200}) and (\ref{fuhleresel000}), and in $(iii)$ we took advantage of cancellations. Formula
 (\ref{zweivors00}) is proven.

  We turn to (c). We may assume ${\bf c}\ne0$ and hence $(\ell,\frac{\alpha}{\beta})=(\frac{(p-2){\xi}}{p-1},1)$. As $\alpha=\beta$ we find$$[\gamma]^{\xi}\varphi
  G_{\gamma} \varphi^Ny\equiv\alpha G_{\gamma}^q(t^{-q^N}[\gamma])^{\xi}\varphi^Ny$$modulo $\nabla$, and it will therefore be enough to show$$\sum_{i=0}^{f-1}([\gamma]^{o_i}-t^{o_i}\frac{t}{[\gamma]}))\varphi^{{N}}y\equiv G_{\gamma}^q(t^{-q^{{N}}}[\gamma])^{\xi}-G_{\gamma}$$modulo $t^{q^{N+n_y-1}}k[[t]]$. Multiplying with the element $t^{-q^{{N}}}[\gamma]$ of $t^{1-q^{{N}}}+t^{2-q^{{N}}}k[[t]]$, we see that this is equivalent with showing$$t^{-q^{{N}}}\sum_{i=0}^{f-1}([\gamma]^{o_i+1}-t^{o_i+1})\varphi^{{N}}y\equiv (t^{-q^{{N}}}[\gamma]G_{\gamma})^{q}-t^{-q^{{N}}}[\gamma]G_{\gamma}$$modulo $t^{n_y}k[[t]]$. With formulae (\ref{amrit00}) and (\ref{voroanti}) and (\ref{sommerdorne}), similarly to the proof of formula (\ref{florian}) we compute$$t^{-q^{{N}}}([\gamma]^{o_i+1}-t^{o_i+1})\equiv a_{\gamma,q}^{p^i}t^{-\frac{p^iq\xi}{p-1}}-a_{\gamma,q}^{p^{i-1}}t^{-\frac{p^{i-1}q\xi}{p-1}}$$modulo $t^{n_y}k[[t]]$ for each $0\le i<f$. Summing up yields$$t^{-q^{{N}}}\sum_{i=0}^{f-1}([\gamma]^{o_i+1}-t^{o_i+1})\equiv a_{\gamma,q}^{p^{f-1}}t^{-\frac{p^{f-1}q\xi}{p-1}}-a_{\gamma,q}^{p^{f-1}}t^{-\frac{p^{f-1}\xi}{p-1}}$$modulo $t^{n_y}k[[t]]$. On the other hand we use formula (\ref{amrit}) and (\ref{sommerdorne}) (resp. (\ref{sommerdornwald1}) if $q=p$) to see $$t^{-q^{{N}}}[\gamma]G_{\gamma}\equiv a_{\gamma,q}^{p^{f-1}}t^{-\frac{p^{f-1}\xi}{p-1}}$$and hence$$t^{-q^{{N}}}[\gamma]G_{\gamma}-(t^{-q^{{N}}}[\gamma]G_{\gamma})^{q}\equiv a_{\gamma,q}^{p^{f-1}}t^{-\frac{p^{f-1}\xi}{p-1}}-a_{\gamma,q}^{p^{f-1}}t^{-\frac{p^{f-1}q\xi}{p-1}}$$modulo $t^{n_y}k[[t]]$. Combining we get what we wanted.\hfill$\Box$\\

\begin{lem}\label{jole} Let $m,s\ge0\in{\mathbb Z}_{\ge0}$ and write $\delta_{m,q^{s}-1}=1$ if $m=q^{s}-1$, but $\delta_{m,q^{s}-1}=0$ otherwise. Modulo $C$ we have in $\Delta$ the congruence\begin{align}t^m\varphi^sx\equiv\delta_{m,q^{s}-1}\alpha^sx&-\sum_{i\in{\mathcal D}}\sum_{0\le v<s\atop m=q^v(q^{N}-r_i-1)+q^s-1}{\bf b}_i\beta^{1+v}\alpha^{s-v}y\notag\\{}&-\sum_{0\le i<f}\sum_{0\le v<s\atop m=q^v(q^{N}-o_i-1)+q^s-1}{\bf c}\beta^{1+v}\alpha^{s-v}y\notag\\{}&-\sum_{i\in{\mathcal E}}\sum_{0\le v<s\atop m=q^v(q^{N}-{{n}}^{(0)}_i-1)+q^s-1}{\bf c}_i\beta^{1+v}\alpha^{s-v}y\notag\\{}&+\sum_{i\in{\mathcal E}}\sum_{j=1}^{f}\sum_{0\le v<s\atop m=q^v(q^{N}-{{n}}^{(j)}_i-1)+q^s-1}{\bf c}_i(m_i+2)\beta^{1+v}\alpha^{s-v}y\notag.\end{align}  
\end{lem}

{\sc Proof:} This follows from looking at the generators of $\nabla$.\hfill$\Box$\\

In Theorem \ref{etaleformeln33} we work out the formulae for the actions of $\varphi$ and $\Gamma$ on the \'{e}tale $(\varphi,\Gamma)$-module $\Delta^*\otimes_{{k[[t]]}}{k((t))}$. They are not really needed for our later purposes, but we will use them as a tool to prove Theorem \ref{classsep0} (and hence Corollary \ref{frfurcht}).

\begin{satz}\label{etaleformeln33} (a) $\Delta^*\otimes_{{k[[t]]}}{k((t))}$ is in a natural way an \'{e}tale $(\varphi,\Gamma)$-module of rank two.

(b) There is a $k((t))$-basis $\lambda_x,\lambda_y$ of $\Delta^*\otimes_{k[[t]]}{k((t))}$ such that for the $\varphi$-linear endomorphism $\varphi$ of $\Delta^*\otimes_{{k[[t]]}}{k((t))}$ we have the formulae $\varphi(\lambda_x)=\alpha^{-1}\lambda_x$ and \begin{align}\varphi(\lambda_y)=&\beta^{-1}\lambda_y-\sum_{0\le i<f}{\bf c}t^{o_i+1-q^N}\lambda_x-\sum_{i\in{\mathcal D}}{\bf b}_it^{r_i+1-q^N}\lambda_x\notag\\{}&\quad\quad-(\sum_{i\in{\mathcal
        E}}{\bf c}_i(t^{{{n}}_i^{(0)}+1-q^N}-(m_i+2)(\sum_{j=1}^{f}t^{{{n}}_i^{(j)}+1-q^N})))\lambda_x.\label{lothoch}\end{align}

  (c) For all $\gamma\in\Gamma$ with $[\gamma]=a_{\gamma,1}t$ we have \begin{gather}\gamma*\lambda_x=a_{\gamma,1}^{-u}\lambda_x\quad\quad\mbox{ and }\quad\quad\gamma*\lambda_y=a_{\gamma,1}^{\ell-u}\lambda_y.\label{theoan}\end{gather}

(d) The ${\mathcal G}_F$-representation $W$ attached to $\Delta^*\otimes_{{k[[t]]}}{k((t))}$ sits in an exact sequence\begin{gather}0\longrightarrow \omega^{-u}\mu_{\alpha^{-1}}\longrightarrow W\longrightarrow \omega^{\ell-u}\mu_{\beta^{-1}}\longrightarrow0.\label{mvpf}\end{gather}
\end{satz}

{\sc Proof:} Statement (d) results from Lemma \ref{allerseelmai} together with statements (b) and (c).

For statement (a) and the fact that $\Delta^*\otimes_{{k[[t]]}}{k((t))}$ has $k((t))$-dimension two (in statement (b)) we use may Lemma \ref{fassa}. Indeed, formula (\ref{univorso})  becomes the set of inequalities $$m+qn_x\ge\xi+\xi(1+q+\ldots+q^{N-1})+n_y,$$$$m\in\{r_i\,|\,i\in{\mathcal D}\}\cup\{n_i^{(0)}\,|\,i\in{\mathcal E}\}\cup\{o_i\,|\,0\le i<f\}$$where however $\{o_i\,|\,0\le i<f\}$ can be omitted if ${\bf c}=0$. But these inequalities are implied by formulae (\ref{marufmin00}), (\ref{maruf200}), (\ref{sommerdornwald6}) (resp. (\ref{maruf2000}) if $q=p$). Notice that formula (\ref{sommerdornwald6}) resp. (\ref{maruf2000}) implies the weaker statement $qn_x-\xi(1+p^{f-1}+\frac{p^{2f-1}}{p-1})\ge n_y$, which is enough for the present purpose.   

Next, identify $V$ with its image in $\Delta$, and inside $\Delta$ define $$\widetilde{V}=tk[[t]]V\quad \mbox{ and }\quad
C=\sum_{s>0}\sum_{0\le\theta<q^{s-1}{\xi}}t^{\theta}\varphi^sV.$$We have the
$k$-vector space decomposition\begin{gather}\Delta= V\oplus\widetilde{V}\oplus
C.\label{wassermuehle}\end{gather}We use it to define
$\lambda_x,\lambda_y\in\Delta^*$ by
requiring $$\lambda_x(\widetilde{V}+ C)=\lambda_y(\widetilde{V}+ C)=0,\quad
\lambda_y(y)=\lambda_x(x)=1,\quad \lambda_y(x)=\lambda_x(y)=0.$$That $\lambda_x,\lambda_y$ is a $k((t))$-basis of $\Delta^*\otimes_{k[[t]]}{k((t))}$ follows e.g. from formula (\ref{cnord}) in the proof of Lemma \ref{fassa} which reduces the statement to a statement on two one-dimensional $k((t))$-vector spaces, proven in Lemma \ref{modvacuumonedim0}. For $\gamma\in\Gamma$ with $[\gamma]=a_{\gamma,1}t$ we compute$$(\gamma*\lambda_x)(x)=\lambda_x(\gamma^{-1}*x)=\lambda_x(a_{\gamma,1}^{-u}x)=a_{\gamma,1}^{-u}=(a_{\gamma,1}^{-u}\lambda_x)(x),$$$$(\gamma*\lambda_y)(y)=\lambda_y(\gamma^{-1}*y)=\lambda_y(a_{\gamma,1}^{\ell-u}y)=a_{\gamma,1}^{\ell-u}=(a_{\gamma,1}^{\ell-u}\lambda_y)(y).$$In the same way we compute $(\gamma*\lambda_x)(y)=0=(\gamma*\lambda_y)(x)$, and moreover $(\gamma*\ell_z)(C+\widetilde{V})=0$ for $z=x$ and $z=y$. This proves formula (\ref{theoan}).

To prove formula (\ref{lothoch}) we need to see that the image of $\lambda_y$ under the map (\ref{pruefeaus1}) is the same as the image
    of the right hand side in formula (\ref{lothoch}) under the map (\ref{pruefeaus0}). Of course, we may
    just as well multiply with $t^{q^N-1}$ and prove that the image of $t^{q^N-1}\lambda_y$ under the map (\ref{pruefeaus1}) is the same as the image
    of $$t^{q^N-1}\,\,\times\,\,(\mbox{the right hand side in formula }(\ref{lothoch}))$$under the map (\ref{pruefeaus0}). Both these image
    elements can be seen as $k$-linear forms on
    $\Delta\otimes_{k[[t]],\varphi}k[[t]]$, and to show that they are the
    same, we show that they take the same value at $z\otimes t^{\kappa}$ for
    all $z\in\Delta$, all $0\le\kappa\le \xi$. This comes down to proving\begin{align}{}&\beta^{-1}\lambda_y(t^{\kappa+q^N-1}\varphi z)-\sum_{i\in{\mathcal D}}{\bf b}_i\lambda_x(t^{\kappa+r_i}\varphi z)-\sum_{0\le i<f}{\bf c}\lambda_x(t^{\kappa+o_i}\varphi z)\notag\\{}&-\sum_{i\in{\mathcal
        E}}{\bf c}_i(\lambda_x(t^{\kappa+{{n}}_i^{(0)}})-(m_i+2)\lambda_x(\sum_{j=1}^{f}t^{\kappa+{{n}}_i^{(j)}}))\varphi z\notag\\=&\left\{\begin{array}{l@{\quad:\quad}l}  \lambda_y(t^{q^{N-1}-1}z)& \kappa=0\\0&1\le\kappa\le \xi\end{array}\right..\label{hochter33}\end{align}It suffices to prove formula (\ref{hochter33}) for all $z=t^{\theta}x$ with $\theta\ge0$ and for all $z=t^{\theta}\varphi^{\iota-1}y$ with $0\le\theta< q^{\iota-2}{\xi}$ and $\iota>1$. This is because $\Delta$ is the $k$-span of such elements, by
    formula (\ref{wassermuehle}).

   For $z=t^{\theta}\varphi^{\iota-1}y$ only the first summand on the left hand side of formula
   (\ref{hochter33}) is possibly non-zero, and as $\beta y=t^{\xi}\varphi y$ and
   $t^{\kappa+q^N-1}\varphi z=t^{\kappa+q^N-1+q\theta}\varphi^{\iota}y$ in $\Delta$ the claim
     follows from the definition of $\lambda_y$ (if $\theta=0$ and $\iota=N$ and
     $\kappa=0$ both sides are $=\beta^{N-1}$, otherwise both sides vanish). 

Now assume $z=t^{\theta}\varphi^{\iota-1}x$. Lemma \ref{jole} with $t^{\kappa+q^N-1}\varphi z=t^{\kappa+q^N+q\theta-1}\varphi^{\iota}x$ shows\begin{align}(\beta^{-1}\lambda_y)(t^{\kappa+q^N-1}\varphi z)=&-\sum_{i\in{\mathcal D}\atop \kappa= q^N-r_i-q\theta-1}{\bf b}_i\alpha-\sum_{0\le i<f\atop \kappa= q^N-o_i-q\theta-1}{\bf c}\alpha\notag\\{}&-\sum_{i\in{\mathcal E}\atop \kappa= q^N-{{n}}_i^{(0)}-q\theta-1}{\bf c}_i\alpha\notag\\{}&+\sum_{i\in{\mathcal E}}\sum_{1\le j\le f\atop \kappa= q^N-{{n}}_i^{(j)}-q\theta-1}{\bf c}_i(m_i+2)\alpha\label{diurnale0}\end{align}if $\iota=N$, but $(\beta^{-1}\lambda_y)(t^{\kappa+q^N-1}\varphi z)=0$ if $\iota\ne N$. On the other hand, Lemma \ref{jole} also shows\begin{gather}\lambda_x(t^{\kappa+r_i}\varphi z)=\alpha\quad\mbox{ if }\kappa=q^N-q\theta-r_i-1,\label{0diurnale1}\\\lambda_x(t^{\kappa+o_i}\varphi z)=\alpha\quad\mbox{ if }\kappa=q^N-q\theta-o_i-1,\label{0diurnale2}\\\lambda_x(t^{\kappa+{{n}}^{(j)}_i}\varphi z)=\alpha\quad\mbox{ if }\kappa=q^N-q\theta-{{n}}^{(j)}_i-1\mbox{ and }0\le j\le f\label{0diurnale4}\end{gather}if $\iota=N$, but that if $\iota\ne N$ or if $\kappa$ is not as stated, then the respective left hand side vanishes. Finally we observe $\lambda_y(t^{q^{N-1}-1}z)=0$ for such $z$, since $t^{\theta}\varphi^{\iota-1}x\in C$ unless $\iota=1$ and $\theta=0$. This, together with formulae (\ref{0diurnale1}), (\ref{0diurnale2}), (\ref{0diurnale4}) proves formula (\ref{hochter00}) for $z=t^{\theta}\varphi^{\iota-1}x$. The above argument concerning $z=t^{\theta}
     \varphi^{\iota-1}y$ also proves $\varphi(\lambda_x)=\alpha^{-1}\lambda_x$.\hfill$\Box$\\

{\bf Definition:} Let $${\mathcal E}(\ell)=\{0\le i<f\,|\,\mbox{ there is some }r<i\mbox{ with }m_{r}=p-1\mbox{ and }m_{s}=p-2\mbox{ for all
}r<s<i\}$$(here $r$ is allowed to be negative) and then ${\mathcal D}(\ell)={{{\mathcal F}}}-{\mathcal E}(\ell)$, so that we have\begin{gather}{{{\mathcal F}}}={\mathcal E}(\ell)\coprod {\mathcal D}(\ell)\label{regenyuval}\end{gather}(disjoint union). For $j\in{\mathbb Z}$
          define \begin{gather}\sigma^{\ell}(j)=\left\{\begin{array}{l@{\quad:\quad}l}1&\Pi(j)\in{\mathcal D}({\ell})\\2& \Pi(j)\in{\mathcal E}({\ell})\end{array}\right..\label{unionbvb}\end{gather}
                                          
\begin{lem}\label{letzterst} (a) For all $i\in {\mathcal E}(\ell)$ we have\begin{gather}{\ell_i}\ge {p^i}(q-\sum_{s=0}^{f-1}p^s).\label{yuvber}\end{gather}
  
  (b) Formula (\ref{fuhleresel200}) is fulfilled if $qp^i(q-\sum_{s=0}^{f-1}p^s)\ge n_y$.  
   
\end{lem}

{\sc Proof:} Formula (\ref{yuvber}) is easily proven. Statement (b) follows from statement (a).\hfill$\Box$\\

Let $i'\in {\mathbb Z}$. Define $i''\le i'$ by $m_{i''-1}\ne p-1$ and
$m_j=p-1$ for all $i''\le j<i'$. If $\Pi(i'')\in {\mathcal D}(\ell)$ put
$i_{\mathcal D}(i')=\Pi(i'')$ and define the integer $m_{\mathcal D}(i')$ by
requiring $i''=i_{\mathcal D}(i')+m_{\mathcal D}(i')f$. If $\Pi(i'')\in
{\mathcal E}(\ell)$ put $i_{\mathcal D}(i')=\Pi(i''-1)$ (which then also
belongs to ${\mathcal E}(\ell)$) and define the integer $m_{\mathcal D}(i')$ by
requiring $i''-1=i_{\mathcal D}(i')+m_{\mathcal D}(i')f$.

Define $i_{\mathcal E}(i')=\Pi(i')$ and then the integer $m_{\mathcal E}(i')$ by requiring $i'=i_{\mathcal E}(i')+m_{\mathcal E}(i')f$.

We have thus defined maps \begin{gather}i_{\mathcal D}^{\ell}:=i_{\mathcal D}:{\mathbb Z}\longrightarrow{{{\mathcal F}}}\label{stetigefreude}\end{gather}and similarly $m_{\mathcal D}, i_{\mathcal E}, m_{\mathcal E}:{\mathbb Z}\to{{{\mathcal F}}}$, all depending on $\ell$.

For $i\in{{{\mathcal F}}}$ let $z_{\mathcal D}(i)$ denote the set of all pairs $(m,i')\in {\mathbb Z}\times{\mathcal D}$ with $i=i_{\mathcal D}(i')$ and $m=m_{\mathcal D}(i')$ and put$$\underline{\bf b}_i=\sum_{(m,i')\in z_{\mathcal D}(i)}(\frac{\alpha}{\beta})^m{\bf b}_{i'}.$$For $i\in {\mathcal E}(\ell)$ let $z_{\mathcal E}(i)$ denote the set of all pairs $(m,i')\in {\mathbb Z}\times{\mathcal E}$ with $i=i_{\mathcal E}(i')$ and $m=m_{\mathcal E}(i')$ and put$$\underline{\bf c}_i=\sum_{(m,i')\in z_{\mathcal E}(i)}(\frac{\alpha}{\beta})^m{\bf c}_{i'}+\sum_{(m,i')\in z_{\mathcal D}(i)}(\frac{\alpha}{\beta})^m{\bf b}_{i'}.$$Put $\underline{\bf c}={\bf c}$. Inside $M$ define the element $$\underline{R}=t^{\xi}\varphi x-\alpha x+\alpha\beta^{1-N}(\sum_{i\in{\mathcal D}(\ell)}\underline{\bf b}_it^{r_i}+\underline{\bf c}\sum_{i=0}^{f-1}t^{o_i}+\sum_{i\in{\mathcal
    E}(\ell)}\underline{\bf c}_it^{n_i^{(0)}}H(m_i,i))\varphi^{{N}}y.$$ Denote by $\underline{\nabla}$ the
$k[[t]][\varphi]$-submodule of $M$ generated by $\nabla_y$, by
$t^{n_x}x$ and by $\underline{R}$. Put $\underline{\Delta}=M/\underline{\nabla}$.

\begin{pro}\label{comme} The same formulae as used for $\Delta$ define an action of $\Gamma$ on $\underline{\Delta}$. The $k[[t]][\varphi,\Gamma]$-modules $\Delta$ and $\underline{\Delta}$ are isomorphic.
\end{pro}

{\sc Proof:} {\it Step 1:}

We make the following claims (a), (b), (c):\\(a) For $i'\in{\mathcal D}$ with $i=i_{\mathcal D}(i')\in {\mathcal D}(\ell)$, putting $m=m_{\mathcal D}(i')$, we have\begin{gather}r_{i'}=q^m(r_i-q^N+1)+q^N-1.\label{yuvdruck1}\end{gather}(b) For $i'\in{\mathcal E}$, putting $i=i_{\mathcal E}(i')$ and $m=m_{\mathcal E}(i')$, we have $i\in {\mathcal E}(\ell)$ and\begin{gather}{{n}}^{(j)}_{i'}=q^m({{n}}^{(j)}_i-q^N+1)+q^N-1\quad\mbox{ for all }0\le j\le f.\label{yuvdruck3}\end{gather}(c) For $i'\in{\mathcal D}$ with $i=i_{\mathcal D}(i')\in {\mathcal E}(\ell)$, putting $m=m_{\mathcal D}(i')$, we have\begin{gather}r_{i'+2f}=q ^m(n_i^{(0)}-q^N+1)+q^N-1\label{neuj}\end{gather}and hence $$t^{n_i^{(0)}}H(m_i,i)=t^{q^{-m}(r_{i'+2f}-q^N+1)+q^N-1}.$$Proofs of these claims:\\(a) To
  see formula (\ref{yuvdruck1})
  observe $$q^m\ell_i-\ell_{i'}={\xi}(p^{mf+i}-p^{i'}).$$(b) That
  $i$ belongs to ${\mathcal E}(\ell)$ follows from formula (\ref{fuhleresel200}). To see formula (\ref{yuvdruck3}) observe $$q^m\ell_i-\ell_{i'}={\xi}(p^{mf+i}-p^{i'}).$$(c) Also formula (\ref{neuj}) is seen in this way.

{\it Step 2:} Let $b,m\ge0$ be such that $q^s(b-q^N+1)+q^N-1\ge0$ for $0\le s<m$ and put $$z=\sum_{s=0}^{m-1}(\frac{\beta}{\alpha})^st^{q^s(b-q^N+1)+q^N-1}\varphi^Ny.$$Claim: Modulo the $k[[t]][\varphi]$-submodule of $M$ generated by $t^{\xi}\varphi y-\beta y$ we have\begin{gather}t^{\xi}\varphi z-\alpha z\equiv \alpha((\frac{\beta}{\alpha})^mt^{q^m(b-q^N+1)+q^N-1}-t^b)\varphi^Ny.\label{yuvalpain}\end{gather}

Proof of this claim: Written out, the left hand side of formula (\ref{yuvalpain}) reads$$\sum_{s=0}^{m-1}(\frac{\beta}{\alpha})^st^{\xi}\varphi t^{q^s(b-q^N+1)+q^N-1}\varphi^Ny-\sum_{s=0}^{m-1}\alpha(\frac{\beta}{\alpha})^st^{q^s(b-q^N+1)+q^N-1}\varphi^Ny.$$But for each $0\le s\le m-1$ we have modulo the $k[[t]][\varphi]$-module generated by $t^{\xi}\varphi y-\beta y$ the congruence$$t^{\xi}\varphi t^{q^s(b-q^N+1)+q^N-1}\varphi^Ny\equiv \beta t^{q^{s+1}(b-q^N+1)+q^N-1}\varphi^Ny.$$Cancelling yields the stated congruence.

{\it Step 3:} Consider the $k[[t]][\varphi]$-module automorphism $\theta$ of $M$ determined by $\theta(y)=y$ and \begin{align}\theta(x)=x&+\sum_{i\in{\mathcal D}(\ell)}\sum_{(m,i')\in z_{\mathcal D}(i)}\alpha^m\beta^{1-N-m}{\bf b}_{i'}\sum_{s=0}^{m-1}(\frac{\beta}{\alpha})^st^{q^s(r_i-q^N+1)+q^N-1}\varphi^Ny\notag\\{}&+\sum_{i\in{\mathcal E}(\ell)}\sum_{(m,i')\in z_{\mathcal E}(i)}\alpha^m\beta^{1-N-m}{\bf c}_{i'}\sum_{j=0}^f\sum_{s=0}^{m-1}(\frac{\beta}{\alpha})^st^{q^s({{n}}_i^{(j)}-q^N+1)+q^N-1}\varphi^Ny\notag\\{}&+\sum_{i\in{\mathcal E}(\ell)}\sum_{(m,i')\in z_{\mathcal D}(i)}\alpha^m\beta^{1-N-m}{\bf b}_{i'}\sum_{j=0}^f\sum_{s=0}^{m-1}(\frac{\beta}{\alpha})^st^{q^s({{n}}_i^{(j)}-q^N+1)+q^N-1}\varphi^Ny.\notag\end{align}This $\theta$ takes $R$ to $\underline{R}$ modulo the $k[[t]][\varphi]$-module generated by $t^{\xi}\varphi y-\beta y$, by what we saw in steps 1 and 2.\hfill$\Box$\\

Proposition \ref{comme} says that assuming ${\mathcal D}={\mathcal D}(\ell)$
and ${\mathcal E}={\mathcal E}(\ell)$ would not restrict the generality of our
constructions.

\begin{lem}\label{voeller} There do exist $n_x, n_y$ satisfying formulae (\ref{sommerdornwald5}), (\ref{marufmin00}), (\ref{fuhleresel200}), (\ref{maruf200}), (\ref{sommerdorne}), (\ref{sommerdornwald1}), (\ref{sommerdornwald6}), (\ref{maruf2000}). In particular, to given scalars $\{{\bf b}_i\}_{i\in{\mathcal D}(\ell)}$, $\{{\bf c}_i\}_{i\in{\mathcal E}(\ell)}$, ${\bf c}$ (now with ${\mathcal D}={\mathcal D}(\ell)$ and ${\mathcal E}={\mathcal E}(\ell)$), our construction yields an \'{e}tale $(\varphi,\Gamma)$-module $\Delta^*\otimes_{{k[[t]]}}{k((t))}$.\end{lem}

  {\sc Proof:} We may take $n_y=1$ (cf. Lemma \ref{letzterst}). Any sufficiently large $n_x$ will then do.\hfill$\Box$\\

We fix $\ell, u,\alpha,\beta$ as before. We give ourselves two sets of scalars $\{{\bf b}_i\}_{i\in{\mathcal D}(\ell)}$, $\{{\bf c}_i\}_{i\in{\mathcal E}(\ell)}$, ${\bf c}$ and $\{{\bf b}'_i\}_{i\in{\mathcal D}(\ell)}$, $\{{\bf c}'_i\}_{i\in{\mathcal E}(\ell)}$, ${\bf c}'$ as before, defining \'{e}tale $(\varphi,\Gamma)$-modules $\Delta^*\otimes_{{k[[t]]}}{k((t))}$ and $(\Delta')^*\otimes_{{k[[t]]}}{k((t))}$ by Lemma \ref{voeller}. 
  
\begin{satz}\label{classsep0} The \'{e}tale $(\varphi,\Gamma)$-modules $\Delta^*\otimes_{{k[[t]]}}{k((t))}$ and $(\Delta')^*\otimes_{{k[[t]]}}{k((t))}$ are isomorphic if and only if there is some $\epsilon\in k^{\times}$ with ${\bf b}_i=\epsilon {\bf b}'_i$ for all $i\in{\mathcal D}(\ell)$, with ${\bf c}_i=\epsilon {\bf c}'_i$ for all $i\in{\mathcal E}(\ell)$ and with ${\bf c}=\epsilon {\bf c}'$.
\end{satz}

{\sc Proof:} If $\epsilon$ as stated does exist, an isomorphism $\Delta\cong\Delta'$ respecting the $\varphi$- and $\Gamma$-actions is obtained by sending $x$ to $\epsilon x$ and $y$ to $y$. We prove the converse. Put $${\mathcal D}=\{0\le i<f\,|\,m_i\ne p-1\},$$$${\mathcal E}=\{i\in {\mathcal E}(\ell)\,|\,m_i\ne p-2\}.$$It is easy to see that the maps $i_{\mathcal D}:{\mathcal D}\to{{{\mathcal F}}}$ and $i_{\mathcal E}:{\mathcal D}\to{{{\mathcal F}}}$ considered above induce a bijection$$ {\mathcal D}\coprod {\mathcal E}\stackrel{\cong}{\longrightarrow}{{{\mathcal F}}}.$$By what we saw in the proof of Proposition \ref{comme} we are therefore allowed to replace ${\mathcal D}(\ell)$, ${\mathcal E}(\ell)$ by ${\mathcal D}$, ${\mathcal E}$.

By Theorem \ref{etaleformeln33}, the basis vectors distinguished by Lemma \ref{modvacuumonedim0}
for the two one
dimensional subquotients of $\Delta^*\otimes_{{k[[t]]}}{k((t))}$ are $\lambda_x$
and (the class of) $\lambda_y$, and similarly we have $\lambda'_x$
and (the class of) $\lambda'_y$ in
$(\Delta')^*\otimes_{{k[[t]]}}{k((t))}$. Let $$\Theta:\Delta^*\otimes_{{k[[t]]}}{k((t))}\stackrel{\cong}{\longrightarrow}(\Delta')^*\otimes_{{k[[t]]}}{k((t))}$$be an isomorphism of \'{e}tale
$(\varphi,\Gamma)$-modules. It sends the unique rank-one (over $k((t))$) sub object of $(\Delta')^*\otimes_{{k[[t]]}}{k((t))}$ spanned by $\lambda_x$ to the unique rank-one sub object of $\Delta^*\otimes_{{k[[t]]}}{k((t))}$ spanned by $\lambda'_x$. The uniqueness-up-to-$k^{\times}$-scalar of the element $\lambda$ in Lemma \ref{modvacuumonedim0} therefore implies that $\Theta$ takes $\lambda_x$ to a $k^{\times}$-scalar multiple of $\lambda'_x$. The same holds analogously for the unique rank-one quotients. Thus, after $k^{\times}$-rescaling the basis vectors if necessary, we have$$\Theta(\lambda_x)=\lambda'_{x},\quad\quad\Theta(\lambda_y)=\lambda'_{y}+F\lambda'_x\quad\mbox{
  for some }F\in k((t)).$$Inserting formula (\ref{lothoch}) into either side
of the equality
$$\Theta(\varphi(\lambda_y))=\varphi(\Theta(\lambda_y))=\varphi(\lambda'_y)+\varphi(F\lambda'_x),$$
we see that with the scalars $b_i={\bf b}'_i-{\bf b}_i$ for $i\in{\mathcal
  D}$, the scalars $c_i={\bf c}'_i-{\bf c}_i$ for $i\in{\mathcal E}$ and the
scalar $c={\bf c}'-{\bf c}$, we have\begin{align}(\varphi-\beta^{-1})F=&\sum_{i\in {\mathcal D}}b_it^{r_i+1-q^N}+c\sum_{i=0}^{f-1}t^{o_i+1-q^N}\notag\\{}&+\sum_{i\in {\mathcal E}}c_i(t^{{{n}}^{(0)}_i+1-q^N}-(m_i+2)\sum_{j=1}^{f}t^{{{n}}^{(j)}_i+1-q^N}).\notag\end{align}We need to see $b_i=0$ for all $i\in{\mathcal D}$, and $c_i=0$ for
all $i\in{\mathcal E}$, and $c=0$. Notice that $r_i+1-q^N<0$ for all $i\in{\mathcal D}$, and ${{n}}^{(0)}_i+1-q^N<0$ for all $i\in{\mathcal D}$, and $o_i+1-q^N<0$ for all $0\le i<f$. It is therefore enough to see that there is
no $h=\sum_{s\in S}h_st^s$ (with $h_s\in k$) in $(\varphi-\beta^{-1})k((t))$, with $h_s\ne0$ for at least one $s\in S$ with $s<0$, where

$$S=\{r_i+1-q^N\}_{i\in{\mathcal D}}\cup\{{{n}}_i^{(j)}+1-q^N\,|\,0\le j\le f\}_{i\in{\mathcal E}}\cup\{o_i^{(j)}+1-q^N\,|\,0\le i<f\}.$$But notice that $c\ne 0$ could possibly happen only if
  $\ell=\frac{(p-2){\xi}}{p-1}$ and that in this case ${\mathcal E}$ is
  empty. Thus, we have$$S=\{r_i+1-q^N\}_{i\in{\mathcal D}}\cup\{{{n}}_i^{(j)}+1-q^N\,|\,0\le j\le f\}_{i\in{\mathcal E}}$$or$$S=\{r_i+1-q^N\}_{i\in{\mathcal D}}\cup\{o_i^{(j)}+1-q^N\,|\,0\le i<f\}.$$For $i\in{\mathcal E}$ consider the integers $m_i^{(j)}$ defined by $q^2(m_i^{(0)}+1-q^N)={{n}}^{(0)}_i+1-q^N$, by $q(m_i^{(j)}+1-q^N)={{n}}^{(j)}_i+1-q^N$ for $1\le j\le i+1$, and by $m_i^{(j)}={{n}}^{(j)}_i$ for $i+1<j\le f$. It follows from Lemma \ref{sepnos} that for $i\in{\mathcal E}$
the numbers $m_i^{(j)}+1-q^N$ for $0\le j\le i+1$ are
pairwise distinct and do not belong to $S$. Adding to the putative $h$ from above suitable elements from $(\varphi-\beta^{-1})k((t))$ we may therefore transform our problem into showing that there is no $\sum_{s\in S'}h'_st^s$ in $(\varphi-\beta^{-1})k((t))$, with $h'_s\ne0$ for at least one $s\in S'$ with $s<0$, where$$S'=\{r_i+1-q^N\}_{i\in{\mathcal D}}\cup\{{m}_i^{(j)}+1-q^N\,|\,0\le j\le f\}_{i\in{\mathcal E}}$$or $S'=\{r_i+1-q^N\}_{i\in{\mathcal D}}\cup\{o_i^{(j)}+1-q^N\,|\,0\le i<f\}$. But this also follows from Lemma \ref{sepnos}, which tells us that none of the elements of $S'$ is divisible by $q$.\hfill$\Box$\\

For $i\in{\mathcal E}$ and
$1\le j\le f$
put$$p_{i,j}=\left\{\begin{array}{l@{\quad:\quad}l}p^{i+1-j} & i+1\ge
j\\p^{i+1-j+f} & i+1< j\end{array}\right.$$
        
$$m_i^{(j)}=q^{{N}}-1+p_{i,j}(p^{j-i-1}q\ell_i-{\xi}(p^{j-1}q+1-\sum_{s=0}^{j-2}p^sq)),$$$$m_i^{(0)}=q^{{N}}-1-2p^i{\xi}+\ell_i.$$

\begin{lem}\label{sepnos} (a) For all $i\in{\mathcal D}$ we have $r_i+1\in
  p^i{\mathbb Z}-p^{i+1}{\mathbb Z}$. For all $i\in{\mathcal E}$ we have $m^{(0)}_i+1\in p^i{\mathbb Z}-p^{i+1}{\mathbb Z}$. For all $i\in{\mathcal E}$ and $1\le
  j\le f$ we have $m^{(j)}_i+1\in p_{i,j}{\mathbb Z}-p_{i,j}p{\mathbb Z}$. For
  all $0\le i<f$ we have $o_i+1\in p^i{\mathbb Z}-p^{i+1}{\mathbb Z}$.

(b) For all $i_1, i_2\in{\mathcal D}$ we have $r_{i_1}\ne r_{i_2}$. For all $i_1, i_2\in{\mathcal E}$ and $0\le j_1, j_2\le f$ with
  $(i_1,j_1)\ne (i_2,j_2)$ we have $m^{(j_1)}_{i_1}\ne m^{(j_2)}_{i_2}$. For
  all $i_1\in{\mathcal D}$ and $i_2\in{\mathcal E}$ and $0\le j_2\le f$ we
  have $r_{i_1}\ne m^{(j_2)}_{i_2}$. If $\ell=\frac{(p-2){\xi}}{p-1}$ then we
  have $o_{i_1}\ne r_{i_2}$ for all $0\le i_1<f$ and $i_2\in{\mathcal D}$. 
  \end{lem}

{\sc Proof:} The statements in (a) are straightforward to check (for
$r_i+1\notin p^{i+1}{\mathbb Z}$ use that $m_i\ne p-1$; for $m^{(0)}_i+1\notin
p^{i+1}{\mathbb Z}$ use that $m_i\ne p-2$). We split up the statements in
(b) into the following statements:

(i) For all $i_1, i_2\in{\mathcal E}$ and $1\le j_1, j_2\le f$ with
$(i_1,j_1)\ne (i_2,j_2)$ we have $m^{(j_1)}_{i_1}\ne m^{(j_2)}_{i_2}$.

(ii) For
  all $i_1\in{\mathcal D}$ and $i_2\in{\mathcal E}$ we
  have $r_{i_1}\ne m^{(0)}_{i_2}$.
 
(iii) For
  all $i_1\in{\mathcal D}$ and $i_2\in{\mathcal E}$ and $1\le j_2\le f$ we
  have $r_{i_1}\ne m^{(j_2)}_{i_2}$. 
 
(iv) For all $i_1, i_2\in{\mathcal E}$ and $1\le j_2\le f$ we have
  $m^{(0)}_{i_1}\ne m^{(j_2)}_{i_2}$.

(v) If $\ell=\frac{(p-2){\xi}}{p-1}$ then we
  have $o_{i_1}\ne r_{i_2}$ for all $0\le i_1<f$ and $i_2\in{\mathcal D}$. 

(vi) For all $i_1, i_2\in{\mathcal D}$ we have $r_{i_1}\ne r_{i_2}$.

To prove these statements it is convenient to put $\ell^{(i)}=p^{-i}\ell_i$
for $0\le i<f$, which is an integer with $0<\ell^{(i)}<\xi$.

To prove (i) we may assume, by (a), that $p_{i_1,j_1}=p_{i_2,j_2}$,
i.e. $i_1-j_1\equiv i_2-j_2$ modulo $f{\mathbb Z}$. We then have $j_1\ne j_2$,
so without loss of generality $j_2<j_1$. We
compute$$m_{i_1}^{(j_1)}-m_{i_2}^{(j_2)}=p_{i_1,j_1}({\xi}q(p^{j_2-1}+p^{j_2}+\ldots+p^{j_1-2})+(\ell^{(i_1)}-\xi)p^{j_1-1}q-(\ell^{(i_2)}-\xi)p^{j_2-1}q).$$Here
${\xi}q(p^{j_2-1}+p^{j_2}+\ldots+p^{j_1-2})+(\ell^{(i_1)}-\xi)p^{j_1-1}q-(\ell^{(i_2)}-\xi)p^{j_2-1}q$
is not divisible by $qp^{j_2}$ since
$\ell^{(i_2)}\equiv m_{i_2}$ modulo $p{\mathbb Z}$, which (as
$i_{2}\in{\mathcal E}$) is not congruent to $-2$ modulo $p{\mathbb
  Z}$.

To prove (ii) we may assume, by (a), that $i_1=i_2$, but then the statement is
immediate.

To prove (iii) we may assume, by (a), that $p^{i_1}=p_{i_2,j_2}$. We then have
$p^{-i_1}(r_{i_1}+1)\equiv \ell^{(i_1)}+1$ modulo $q$, and on the other hand
$p^{-i_1}(m^{(j_2)}_{i_2}+1)\equiv 1$ modulo $q$, and these are not the same (as $0<\ell^{(i_1)}<\xi$).

To prove (iv) we may assume, by (a), that $p^{i_1}=p_{i_2,j_2}$. We then have
$p^{-i_1}(m^{(0)}_{i_1}+1)\equiv \ell^{(i_1)}+2$ modulo $q$, and on the other hand
$p^{-i_1}(m^{(j_2)}_{i_2}+1)\equiv 1$ modulo $q$, and these are not the same (as $0<\ell^{(i_1)}<\xi$).

To prove (v) we may assume, by (a), that $i_1=i_2$, but then the statement is
immediate.

To prove (vi) use (a).\hfill$\Box$\\

\begin{kor}\label{frfurcht} (a) $\Delta^*\otimes_{{k[[t]]}}{k((t))}$ is an
  extension\footnote{for the definition of ${\bf E}(.,.)$ see Lemma \ref{modvacuumonedim0}} of ${\bf E}(\beta,u-\ell)$ by ${\bf E}(\alpha,u)$. 

(b) Assume $(\ell,\frac{\alpha}{\beta})\ne(\frac{(p-2){\xi}}{p-1},1)$. Sending the tuple $(({\bf b}_i)_{i\in{\mathcal D}(\ell)},({\bf c}_i)_{i\in{\mathcal E}(\ell)})$ to the class of
  $\Delta^*\otimes_{{k[[t]]}}{k((t))}$ induces an
  isomorphism between the $k$-vector space of scalar tuples $(({\bf b}_i)_{i\in{\mathcal D}(\ell)},({\bf c}_i)_{i\in{\mathcal E}(\ell)})$ and ${\rm Ext}^1({\bf E}(\beta,u-\ell),{\bf E}(\alpha,u))$.

(c) Assume $(\ell,\frac{\alpha}{\beta})=(\frac{(p-2){\xi}}{p-1},1)$. Sending the tuple $(({\bf b}_i)_{i\in{\mathcal D}(\ell)},({\bf c}_i)_{i\in{\mathcal E}(\ell)},{\bf c})$ to the class of
  $\Delta^*\otimes_{{k[[t]]}}{k((t))}$ induces an
  isomorphism between the $k$-vector space of scalar tuples $(({\bf b}_i)_{i\in{\mathcal D}(\ell)},({\bf c}_i)_{i\in{\mathcal E}(\ell)},{\bf c})$ and ${\rm Ext}^1({\bf E}(\beta,u-\ell),{\bf E}(\alpha,u))$.
\end{kor}

{\sc Proof:} (b) and (c): Theorem \ref{classsep0} says that the stated maps to ${\rm
  Ext}^1({\bf E}(\alpha,u),{\bf E}(\beta,u-\ell))$ are injective. To see
surjectivity it suffices to compare the $k$-vector space dimensions. The decomposition (\ref{regenyuval}) shows $|{\mathcal D}(\ell)|+|{\mathcal E}(\ell)|=f$. Thus, it suffices to observe that ${\rm Ext}^1({\bf E}(\beta,u-\ell),{\bf E}(\alpha,u))$ has dimension $f$ resp. $f+1$. But it is well known that the extension spaces of the corresponding Galois characters have the said dimensions.\hfill$\Box$\\

\subsection{The case $\ell=0$}

\label{secthreeone}

Here what was previously $\ell$ is taken to be $\ell=0$. This leads to \'{e}tale $(\varphi,\Gamma)$-modules corresponding to ${\mathcal G}_F$-representations which are (reducible and) unramified after suitable character twist.

The outline and the arguments in this subsection are only a slight variation of those of the preceeding one.

Fix an integer $0\le u<\xi$ and scalars $\alpha,\beta\in k^{\times}$. Fix a finite subset ${\mathcal D}\subset{\mathbb Z}_{\ge0}$, scalars ${\bf d}_i\in k$ for $i\in{\mathcal D}$, and a scalar ${\bf e}\in k$. Fix $N, n_x\in{\mathbb Z}_{>0}$ so that for all $i\in{\mathcal D}$ we have $q^N\ge p^i\xi+1$ and\begin{gather}q(n_x-1)\ge p^i{\xi}-1.\label{freu}\end{gather}Let $V$ be the two dimensional $k$-vector space on the basis $x,y$. Inside the $k[[t]][\varphi]$-module $M=k[[t]][\varphi]\otimes_kV=k[[t]][\varphi]V$ define the $k[[t]][\varphi]$-submodule $\nabla$ to be generated by the elements $t^{n_x}x$ and $ty$ and $t^{\xi}\varphi y-\beta y$ and $$R=t^{\xi}\varphi x-\alpha x-{\bf e}\alpha\beta y-\alpha\beta^{1-N}(\beta-\alpha)\sum_{i\in{\mathcal D}}{\bf d}_it^{q^N-1-p^i{\xi}}\varphi^Ny.$$

\begin{pro}\label{discord0} (a) The $k[[t]][\varphi]$-action on
  $M$ extends to an action by $k[[t]][\varphi,\Gamma]$
  satisfying the following formulae\footnote{For later comparison notice that as $ty\in\nabla$ we have $a_{\gamma,1}^{u}y\equiv a_{\gamma,1}^{u+1}\frac{t}{[\gamma]}y$ modulo $\nabla$.} for all
  $\gamma\in\Gamma$:$$\gamma*x=a_{\gamma,1}^u(\frac{a_{\gamma,1}t}{[\gamma]}x+\alpha\beta\sum_{i\in{\mathcal D}}{\bf d}_i a_{\gamma,q}^{p^i}y),$$$$\gamma*y=a_{\gamma,1}^{u}y.$$This
  action preserves $\nabla$, hence provides $\Delta=M/\nabla$ with the
  structure of a $k[[t]][\varphi,\Gamma]$-module.

  (b) If $\alpha\ne \beta$ then the isomorphism class of $\Delta$ is independent of ${\bf e}$.\footnote{In particular, if $\alpha\ne \beta$ then one may assume without loss of generality that ${\bf e}=0$.}
\end{pro}

{\sc Proof:} For (a) the arguments are similar to those invoked in the proof of Proposition \ref{vacet10}, and run as follows. It is straightforward to check that the stated formulae define an action of $\Gamma$ on $M$; invoke formula (\ref{freivorquinqua1}) for this. The critical claim is the stability of $\nabla$ under the action of $\Gamma$. First, if $\gamma\in\Gamma$ is such that $[\gamma]=a_{\gamma,1}t$ then $\gamma*\nabla=\nabla$ is obvious. Thus, we may now assume that $\gamma\in\Gamma$ is such that $a_{\gamma,1}=1$. We then compute\begin{align}\gamma*R-\frac{t}{[\gamma]}R\equiv&-\alpha\beta^{1-N}(\beta-\alpha)\sum_{i\in{\mathcal D}}{\bf d}_i ([\gamma]^{q^N-1-p^i{\xi}}-\frac{t^{q^N-p^i{\xi}}}{[\gamma]})\varphi^Ny\notag\\{}&+\alpha\beta\sum_{i\in{\mathcal D}} {\bf d}_i a_{\gamma,q}^{p^i}t^{\xi}\varphi y-\alpha^2\beta\sum_{i\in{\mathcal D}}{\bf d}_i a_{\gamma,q}^{p^i}y\notag\end{align}modulo the $k[[t]][\varphi]$-submodule generated by $ty$, where as usual we take advantage of formula (\ref{amrit}). But formula (\ref{amrit}) also yields$$([\gamma]^{q^N-1-p^i{\xi}}-\frac{t^{q^N-p^i{\xi}}}{[\gamma]})\varphi^Ny\equiv a_{\gamma,q}^{p^i}t^{q^N-1}\varphi^Ny$$modulo the $k[[t]][\varphi]$-submodule generated by $ty$ for each ${i\in{\mathcal D}}$ since, similarly as before, we have the modulo $p{\mathbb Z}$ congruences$${q^N-p^i{\xi}\choose p^i}\equiv 1\quad\mbox{ and }\quad {q^N-p^i{\xi}\choose m}\equiv0$$for $m\in[1,2p^i-1]-\{p^i\}$. But$$-\frac{\beta-\alpha}{\beta^N}t^{q^N-1}\varphi^Ny+ t^{\xi}\varphi y-\alpha y$$belongs to the $k[[t]][\varphi]$-submodule generated by $t^{\xi}\varphi y-\beta y$, so that altogether $\gamma*R\in\nabla$, as desired. To see statement (b) put $$\widetilde{x}=x+\frac{{\bf e}\beta\alpha}{\beta-\alpha}y.$$We then find $$t^{\xi}\varphi \widetilde{x}-\alpha \widetilde{x}-\alpha\beta^{1-N}(\beta-\alpha)\sum_{i\in{\mathcal D}}{\bf d}_it^{q^N-1-p^i{\xi}}\varphi^Ny\in\nabla$$which shows that the isomorphism $V\to V$ with $y\mapsto y$ and $x\mapsto\widetilde{x}$ provides the searched-for isomorphism between $\Delta$ for general ${\bf e}$ and $\Delta$ for ${\bf e}=0$.\hfill$\Box$\\

We identify $V$ with its image in $\Delta$, and inside
$\Delta$ we define $$\widetilde{V}=tk[[t]]V\quad \mbox{ and }\quad
C=\sum_{s>0}\sum_{0\le\theta<q^{s-1}{\xi}}t^{\theta}\varphi^sV.$$We then have the $k$-vector space decomposition\begin{gather}\Delta= V\oplus\widetilde{V}\oplus
C.\label{wassermuehle0}\end{gather}We use it to define
$\lambda_x,\lambda_y\in\Delta^*$ by
requiring $$\lambda_x(\widetilde{V}+ C)=\lambda_y(\widetilde{V}+ C)=0,\quad
\lambda_y(y)=\lambda_x(x)=1,\quad \lambda_y(x)=\lambda_x(y)=0.$$As before, $\lambda_x,\lambda_y$ is a basis of the $k((t))$-vector space
$\Delta^*\otimes_{{k[[t]]}}{k((t))}$. We consider the \'{e}tale
$(\varphi,\Gamma)$-module $\Delta^*\otimes_{{k[[t]]}}{k((t))}$, cf. Theorem
  \ref{semferiend}.

\begin{satz}\label{etaleformel} For the $\varphi$-linear endomorphism
  $\varphi$ of $\Delta^*\otimes_{{k[[t]]}}{k((t))}$ we have the
  formulae $$\varphi(\lambda_x)=\alpha^{-1}\lambda_x\quad\mbox{ and }\quad\varphi(\lambda_y)=\beta^{-1}\lambda_y+\sum_{i\in{\mathcal D}}{\bf d}_i(\beta-\alpha)t^{-p^i{\xi}}\lambda_x+{\bf e}\lambda_x.$$For all $\gamma\in\Gamma$ with $[\gamma]=a_{\gamma,1}t$ we have \begin{gather}\gamma*\lambda_x=a_{\gamma,1}^{-u}\lambda_x\quad\quad\mbox{ and }\quad\quad\gamma*\lambda_y=a_{\gamma,1}^{-u}\lambda_y.\label{theoan0}\end{gather}For all $\gamma\in\Gamma$ with $a_{\gamma,1}=1$ we have the congruences \begin{gather}\gamma*\lambda_x\equiv \lambda_x\quad\quad\mbox{ and }\quad\quad\gamma*\lambda_y\equiv \lambda_y+\sum_{i\in{\mathcal D}}{\bf d}_i\alpha\beta a_{\gamma,q}^{p^i}\lambda_x\label{samstagoffen}\end{gather}modulo (the image of) $tk[[t]]\lambda_x+tk[[t]]\lambda_y$. The ${\mathcal G}_F$-representation $W$ attached to $\Delta^*\otimes_{{k[[t]]}}{k((t))}$ sits in an exact sequence\begin{gather}0\longrightarrow \omega^{-u}\mu_{\alpha^{-1}}\longrightarrow W\longrightarrow \omega^{-u}\mu_{\beta^{-1}}\longrightarrow0.\label{mvpf1}\end{gather}
\end{satz}

{\sc Proof:} The formulae for the $\varphi$-action are seen exactly as in the proof of
Theorem \ref{etaleformeln33}. To see formulae (\ref{theoan0}) and (\ref{samstagoffen}) argue as in the proof of
(Theorem \ref{etaleformeln33} and) Lemma \ref{modvacuumonedim0}. The exact sequence (\ref{mvpf1}) results from Lemma \ref{allerseelmai}.\hfill$\Box$\\

Put $${\mathcal D}(0)={{{\mathcal F}}},\quad\quad\quad {\mathcal E}(0)=\emptyset.$$For $i'\in {\mathbb Z}$ define the integer $m(i')=\frac{i'-\Pi(i')}{f}$. For $i\in{\mathcal D}(0)$ let $z(i)$ denote the set of all pairs $(m,i')\in {\mathbb Z}\times{\mathcal D}$ with $i=\Pi(i')$ and $m=m(i')$ and put$$\underline{\bf d}_i=\sum_{(m,i')\in z(i)}(\frac{\alpha}{\beta})^m{\bf d}_{i'}.$$Inside $M$ define the element $$\underline{R}=t^{\xi}\varphi x-\alpha x-{\bf e}\alpha\beta y-\alpha\beta^{1-N}(\beta-\alpha)\sum_{i\in{\mathcal D}(0)}\underline{\bf d}_it^{q^N-1-p^i{\xi}}\varphi^Ny.$$ Denote by $\underline{\nabla}$ the
$k[[t]][\varphi]$-submodule of $M$ generated by the elements $t^{n_x}x$ and $ty$ and $t^{\xi}\varphi y-\beta y$ and $\underline{R}$. Put $\underline{\Delta}=M/\underline{\nabla}$.

\begin{pro}\label{comme000} The same formulae as used for $\Delta$ define an action of $\Gamma$ on $\underline{\Delta}$. The $k[[t]][\varphi,\Gamma]$-modules $\Delta$ and $\underline{\Delta}$ are isomorphic.
\end{pro}

{\sc Proof:} This is only a much easier version of the proof of Proposition \ref{comme}.\hfill$\Box$\\

We give ourselves two sets of scalars $\{{\bf d}_i\}_{i\in {\mathcal D}(0)}$, ${\bf e}$ and $\{{\bf d}'_i\}_{i\in {\mathcal D}(0)}$, ${\bf e}'$, defining \'{e}tale $(\varphi,\Gamma)$-modules $\Delta^*\otimes_{{k[[t]]}}{k((t))}$ and $(\Delta')^*\otimes_{{k[[t]]}}{k((t))}$.

\begin{satz}\label{classsep0ell0} (a) Assume $\alpha\ne \beta$. The \'{e}tale $(\varphi,\Gamma)$-modules $\Delta^*\otimes_{{k[[t]]}}{k((t))}$ and $(\Delta')^*\otimes_{{k[[t]]}}{k((t))}$ are isomorphic if and only if there is some $\epsilon\in k^{\times}$ with ${\bf d}_i=\epsilon {\bf d}'_i$ for all ${i\in {\mathcal D}(0)}$.

(b) Assume $\alpha=\beta$. The \'{e}tale $(\varphi,\Gamma)$-modules $\Delta^*\otimes_{{k[[t]]}}{k((t))}$ and $(\Delta')^*\otimes_{{k[[t]]}}{k((t))}$ are isomorphic if and only if there is some $\epsilon\in k^{\times}$ with ${\bf d}_i=\epsilon {\bf d}'_i$ for all ${i\in {\mathcal D}(0)}$ and with ${\bf e}= \epsilon{\bf e}'$. 
\end{satz}

{\sc Proof:} If $\epsilon$ as stated does exist, an isomorphism $\Delta\cong\Delta'$ respecting the $\varphi$- and $\Gamma$-actions is obtained by sending $x$ to $\epsilon x$ and $y$ to $y$ (in case (a) use that the value of ${\bf e}$ does not influence the isomorphism class of $\Delta$, as seen in Proposition \ref{discord0}).

We prove the converse. By Theorem \ref{etaleformeln33}, the basis vectors distinguished by Lemma \ref{modvacuumonedim0}
for the two one
dimensional subquotients of $\Delta^*\otimes_{{k[[t]]}}{k((t))}$ are $\lambda_x$
and (the class of) $\lambda_y$, and similarly we have $\lambda'_x$
and (the class of) $\lambda'_y$ in
$(\Delta')^*\otimes_{{k[[t]]}}{k((t))}$. Therefore an isomorphism of \'{e}tale
$(\varphi,\Gamma)$-modules $$\Theta:\Delta^*\otimes_{{k[[t]]}}{k((t))}\stackrel{\cong}{\longrightarrow}(\Delta')^*\otimes_{{k[[t]]}}{k((t))}$$
is given (after $k^{\times}$-rescaling the basis vectors if necessary) by the
assignments$$\Theta(\lambda_x)=\ell'_{x},\quad\quad\Theta(\lambda_y)=\ell'_{y}+F\lambda'_x\quad\mbox{
  for some }F\in k((t)).$$Inserting formula (\ref{lothoch}) into either side
of the equality
$$\Theta(\varphi(\lambda_y))=\varphi(\Theta(\lambda_y))=\varphi(\lambda'_y)+\varphi(F\lambda'_x),$$
we see that with the scalars $d_i={\bf d}_i-{\bf d}'_i$ for ${i\in {\mathcal D}(0)}$ and the
scalar $e={\bf e}-{\bf e}'$, we have\begin{gather}(\varphi-\beta^{-1})F=\sum_{{i\in {\mathcal D}(0)}}d_i(\beta-\alpha)t^{-p^i{\xi}}+e.\label{simsam}\end{gather}We need to see $d_i=0$ for all ${i\in {\mathcal D}(0)}$, and if $\alpha=\beta$ then also $e=0$. If $\alpha\ne\beta$ then $d_i=0$ for all ${i\in {\mathcal D}(0)}$ already follows from formula (\ref{simsam}), by exactly the same proof as for Theorem \ref{classsep0} (it is only much easier here). This same proof also shows $e=0$ if $\alpha=\beta$. It remains to see $d_i=0$ for all ${i\in {\mathcal D}(0)}$ if $\alpha=\beta$. Formula (\ref{samstagoffen}) for $\gamma\in\Gamma$ with $a_{\gamma,1}=1$ yields the congruence\begin{align}\lambda'_y+F\lambda'_x+\sum_{{i\in {\mathcal D}(0)}}{\bf d}_i\alpha\beta a_{\gamma,q}^{p^i}\lambda'_x&\equiv\Theta(\gamma*\lambda_y)=\gamma*(\Theta(\lambda_y))\notag\\{}&\equiv \lambda'_y+\sum_{{i\in {\mathcal D}(0)}}{\bf d}'_i\alpha\beta a_{\gamma,q}^{p^i}\lambda'_x+(\gamma*F)\lambda'_x\notag\end{align}modulo $tk[[t]]\lambda_x+tk[[t]]\lambda_y$. Since $F$ and $\gamma*F$ have the same constant term it follows that $\sum_{{i\in {\mathcal D}(0)}}{\bf d}_ia_{\gamma,q}^{p^i}=\sum_{{i\in {\mathcal D}(0)}}{\bf d}'_i a_{\gamma,q}^{p^i}$. As this is true for all $\gamma\in\Gamma$ with $a_{\gamma,1}=1$ we get ${\bf d}_i={\bf d}'_i$ for all ${i\in {\mathcal D}(0)}$ (use formula (\ref{ruofra}) to see that each element of $k$ is of the form $a_{\gamma,q}$ for a suitable $\gamma\in \Gamma$ (with $\gamma\equiv 1$ in $k$)).\hfill$\Box$\\

\begin{kor}\label{frfur} (a) $\Delta^*\otimes_{{k[[t]]}}{k((t))}$ is an
  extension of ${\bf E}(\beta,u)$ by ${\bf E}(\alpha,u)$. 

(b) Let $\alpha\ne \beta$. Sending the tuple $({\bf d}_i)_{{i\in {\mathcal D}(0)}}$ to the class of
  $\Delta^*\otimes_{{k[[t]]}}{k((t))}$ induces an
  isomorphism between the $k$-vector space of scalar tuples $({\bf d}_i)_{{i\in {\mathcal D}(0)}}$ and ${\rm Ext}^1({\bf E}(\beta,u),{\bf E}(\alpha,u))$.

(c) Let $\alpha= \beta$. Sending the tuple $(({\bf d}_i)_{{i\in {\mathcal D}(0)}},{\bf e})$ to the class of $\Delta^*\otimes_{{k[[t]]}}{k((t))}$ induces an isomorphism between the $k$-vector space of scalar tuples $(({\bf d}_i)_{{i\in {\mathcal D}(0)}},{\bf e})$ and ${\rm Ext}^1({\bf E}(\beta,u),{\bf E}(\alpha,u))$.

\end{kor}

{\sc Proof:} (b) and (c): Theorem \ref{classsep0ell0} says that the stated maps to ${\rm
  Ext}^1({\bf E}(\alpha,u),{\bf E}(\beta,u))$ are injective. To see
surjectivity it suffices to observe that ${\rm Ext}^1({\bf E}(\beta,u),{\bf E}(\alpha,u))$ has dimension $f$ resp. $f+1$. But it is well known that the extension spaces of the corresponding Galois characters have the said dimensions.\hfill$\Box$\\
           
\subsection{The stack's strata set $\breve{\mathbb E}$}

\label{stackstrataset}

{\bf Definition:} (a) Let ${\mathbb E}$ denote the set of triples $(\ell,u,{\mathcal
  I})$ with ${\mathcal
  I}\subset{{{\mathcal F}}}$ and $\ell,u\in[0,q-2]$.

(b) For $(\ell,u,{\mathcal
   I})\in{\mathbb E}$ and $\alpha,\beta\in k^{\times}$ we define in ${\rm Ext}^1({\bf E}(\beta,u-\ell),{\bf
  E}(\alpha,u))$ the following $k$-subspace $E^{\alpha,\beta}_{{\mathcal
  I}}$:

$\bullet$ If $\ell\ne0$ and $(\ell,\frac{\alpha}{\beta})\ne(\frac{(p-2){\xi}}{p-1},1)$ it corresponds, under the bijection in Corollary \ref{frfurcht}, to tuples $(({\bf
  b}_i)_{i\in{\mathcal D}(\ell)},({\bf c}_i)_{i\in{\mathcal
  E}(\ell)})$ with ${\bf b}_i=0$ resp. ${\bf
  c}_i=0$ for $i\notin{\mathcal
  I}$.

$\bullet$  If $(\ell,\frac{\alpha}{\beta})=(\frac{(p-2){\xi}}{p-1},1)$ it corresponds, under the bijection in Corollary \ref{frfurcht}\footnote{notice that $\ell=\frac{(p-2){\xi}}{p-1}$ implies ${\mathcal E}(\ell)=\emptyset$ and ${\mathcal D}(\ell)={{{\mathcal F}}}$}, to tuples $(({\bf
  b}_i)_{i\in{\mathcal D}(\ell)},{\bf
  c})$ with ${\bf b}_i=0$ for $i\notin{\mathcal I}$, and with ${\bf c}=0$.

$\bullet$  If $\ell=0$ and ${\alpha}\ne\beta$ it corresponds, under the bijection in Corollary
  \ref{frfur}, to tuples $(({\bf d}_i)_{0\le i<f})$ with ${\bf
  d}_i=0$ for $i\notin {\mathcal I}$.

$\bullet$  If $\ell=0$ and ${\alpha}=\beta$ it corresponds, under the bijection in
    Corollary \ref{frfur}, to tuples $(({\bf
  d}_i)_{0\le i<f},{\bf e})=0$ with ${\bf
  d}_i=0$ for $i\notin {\mathcal I}$, and with ${\bf e}=0$.

(c) For $(\ell,u,{\mathcal
  I})\in{\mathbb E}$ let ${\mathcal V}(\ell,u,{\mathcal
  I})$ denote the set of $(\varphi,\Gamma)$-modules whose class belong to $E^{\alpha,\beta}_{{\mathcal
  I}}$ for some $\alpha,\beta\in k^{\times}$. Elements
    of ${\mathcal V}(\ell,u,{\mathcal
  I})$ correspond to certain\footnote{those whose class in ${\rm Ext}^1$ belong to the span of the basis elements indexed by ${\mathcal I}$} reducible ${\mathcal G}_F$-representations
    $W$ which, when restricted to ${\mathcal I}_F$, sit in an exact sequence$$0\longrightarrow
\omega^{-{u}}\longrightarrow W\longrightarrow
\omega^{{\ell}-{u}}\longrightarrow0.$$For $(\ell,u,{\mathcal
  I})\in{\mathbb E}$ put \begin{gather}{\mathcal V}^0(\ell,u,{\mathcal
  I})={\mathcal V}(\ell,u,{\mathcal
  I})-\bigcup_{{\mathcal I}'\subsetneq {\mathcal I}}{\mathcal V}(\ell,u,{\mathcal
  I}').\label{fernandovertex}\end{gather}For $h\in {\mathbb Z}$ let ${\mathcal V}(h)$ denote the set of $(\varphi,\Gamma)$-modules whose associated ${\mathcal G}_F$-representation, when restricted to ${\mathcal I}_F$, is isomorphic with $\omega_{2f}^h\oplus\omega_{2f}^{qh}$ (thus ${\mathcal V}(h)={\mathcal V}(qh)$).

(d) Put $\breve{\mathbb E}={\mathbb E}\coprod{\mathbb Z}/(q^2-1){\mathbb Z}$.\\

{\bf Remark:} By construction, for $(\ell,u,{\mathcal
  I})\in{\mathbb E}$ the set ${\mathcal V}(\ell,u,{\mathcal
  I})$ is the set of $k$-points of an algebraic family ${\bf D}(\ell,u,{\mathcal
  I})$ of
$(\varphi,\Gamma)$-modules, defined over $({\mathbb G}_m)^2\times {\mathbb A}(\ell,u,{\mathcal
  I})$ for some affine $k$-space ${\mathbb A}(\ell,u,{\mathcal
  I})$ of dimension $|{\mathcal I}|$ (resp. of dimension $|{\mathcal I}|+1$ in
the exceptional cases). Better, we may {\it canonically}
identify ${\mathcal I}$ (resp. ${\mathcal I}$ together with one more element in
the exceptional cases) with a set of {\it distinguished} lines in ${\mathbb A}(\ell,u,{\mathcal
  I})$, which together span ${\mathbb A}(\ell,u,{\mathcal
  I})$. One may ask for the degenerations of the family ${\bf D}(\ell,u,{\mathcal
  I})$ into $(\varphi,\Gamma)$-modules which do {\it not} belong to ${\mathcal V}(\ell,u,{\mathcal
  I})$, i.e. for families of $(\varphi,\Gamma)$-modules ${\bf D}$,
defined over (say) an open subspace ${\mathbb B}$ of a suitable affine
$k$-space, which on an open subspace ${\mathbb B}^0$ of ${\mathbb B}$ restrict
to subfamilies of ${\bf D}(\ell,u,{\mathcal
  I})$, but whose $k$-points in ${\mathbb B}-{\mathbb B}^0$ belong to either
${\mathcal V}(h)$ for some $h$, or to ${\mathcal V}(\ell',u',{\mathcal
  I}')$ for some triple $(\ell',u',{\mathcal
  I}')$ with $(\ell',u')\ne (\ell,u)$.

The purpose of sections \ref{okt3} and \ref{zoomgergelyfengwei} is to describe
such degenerations, and in fact all of them. Namely, in section \ref{okt3},
for any ${\mathcal I}$ with $|{\mathcal I}|=1$ we describe such a degeneration
of ${\bf D}(\ell,u,{\mathcal
  I})$ into elements of ${\mathcal V}(h)$, for a specific $h$ depending on
$\ell, u,{\mathcal I}$. Similarly, in section \ref{zoomgergelyfengwei}, for
any ${\mathcal I}$ with $|{\mathcal I}|\ge 2$ and any $\Pi(i_1)\in {\mathcal I}$ we describe a degeneration
of ${\bf D}(\ell,u,{\mathcal
  I})$ into a
family ${\bf D}(\ell',u',{\mathcal
  I}')$ for some $(\ell',u',{\mathcal
  I}')$ with $(\ell',u')\ne (\ell,u)$ and $|{\mathcal
  I}'|=|{\mathcal
  I}|-1$; namely, we let the $\Pi(i_1)$-coordinate (inside ${\mathbb A}(\ell,u,{\mathcal
  I})$) 'approach and finally attain infinity'.\footnote{In section \ref{zoomgergelyfengwei} we write $(\overline{\ell},\overline{u},\overline{\mathcal
  I})$ and $(\widetilde{\ell},\widetilde{u},\widetilde{\mathcal
  I})$ rather than $(\ell',u',{\mathcal
  I}')$ and $(\ell,u,{\mathcal
  I})$.}\\

{\bf Remark:} We may stratify the
stack ${\mathcal X}_{2,{\rm red}}$ parametrizing $(\varphi,\Gamma)$-modules of rank two, as described in \cite{eg}, resp. \cite{pham}, as follows:

${\bullet}$ Up to the equalities ${\mathcal V}(h)={\mathcal V}(qh)$ resp. ${\mathcal V}(\ell,u,\emptyset)={\mathcal
    V}^0(\ell,u,\emptyset)={\mathcal V}^0(\ell',u',\emptyset)={\mathcal V}(\ell',u',\emptyset)$ if $\ell\equiv
  u-u'\equiv-\ell'$ modulo $\xi{\mathbb Z}$, and

  ${\bullet}$ up to neglecting $(\varphi,\Gamma)$-modules for which ${\bf e}\ne 0$, resp. ${\bf c}\ne0$,\footnote{The degeneration behaviour of the $(\varphi,\Gamma)$-modules on these loci is trivial --- it is captured by moving into the loci with ${\bf e}=0$ resp. ${\bf c}=0$ ---  so we do not loose anything for our purposes.}\\we make an element of $\breve{\mathbb E}$ correspond to a stratum of ${\mathcal X}_{2,{\rm red}}$ according to the following rule: ${\mathcal V}^0(\ell,u,{\mathcal
  I})$ (resp. ${\mathcal V}(h)$) is the set of $k$-points of the stratum
corresponding to $(\ell,u,{\mathcal
  I})\in{\mathbb E}$ (resp. to $[h]\in {\mathbb
  Z}/(q^2-1){\mathbb Z}$). The degeneration properties just
discussed then precisely reflect the closure relation between the said strata
of ${\mathcal X}_{2,{\rm red}}$.

This stratification of ${\mathcal X}_{2,{\rm red}}$ refines the one given by its decomposition into irreducible components. At least ignoring exceptional components, the strata ${\mathcal V}(\ell,u,{\mathcal
  F})$ (or just ${\mathcal V}^0(\ell,u,{\mathcal
  F})$, which comes down to the same) correspond to the irreducible components of the stack; in fact, the 'interior' of the latter are obtained simply by quotienting out from our ${\mathcal V}(\ell,u,{\mathcal
  F})$-family a torus automorphism group.

\section{Degeneration into irreducibles}

\label{okt3}

The purpose of this section is to prove Corollary \ref{cusanus00}.

Let us try to briefly explain the strategy. We would like to non-trivially deform a family of reducible $(\varphi,\Gamma)$-modules as constructed before. Recall that such $(\varphi,\Gamma)$-modules arose in particular (i.e. here we do not need their most general shape) from a collection of defining relations$$t\widetilde{x},\quad\quad  t^{n_{\widetilde{y}}}\widetilde{y},\quad\quad  a\widetilde{y}+ t^{{\xi}}\varphi
\widetilde{y}-t^{q^N-1-w}G(t)\varphi^N\widetilde{x},\quad\quad \widetilde{x}-t^{{\xi}}\varphi
\widetilde{x}$$in $k[[t]][\varphi](k.\widetilde{x}\oplus k.\widetilde{y})$, with $a\in k^{\times}$, certain $G(t)\in 1+tk[[t]]$, certain $n_{\widetilde{y}}\in{\mathbb N}$ and $w\in{\mathbb Z}$. We make this into a family parametrized by $\tau\in k^{\times}$ by instead looking at$$t^{n_{\widetilde{x}}}\widetilde{x},\quad\quad t^{n_{\widetilde{y}}}\widetilde{y},\quad\quad a\widetilde{y}+\tau t^{{\xi}}\varphi\widetilde{y}-\tau^{1-N}t^{q^N-1-w}G(t)\varphi^N\widetilde{x},\quad\quad \tau \widetilde{x}-t^{{\xi}}\varphi
\widetilde{x}.$$Apparently, specifying the relations $a\widetilde{y}+\tau t^{{\xi}}\varphi\widetilde{y}-\tau^{1-N}t^{q^N-1-w}G(t)\varphi^N\widetilde{x}$ and $\tau \widetilde{x}-t^{{\xi}}\varphi \widetilde{x}$ at $\tau=0$ makes the construction collapse and does not yield a rank two $(\varphi,\Gamma)$-module. However, for $\tau\in k^{\times}$ it turns out that instead of stipulating the relations $a\widetilde{y}+\tau t^{{\xi}}\varphi\widetilde{y}-\tau^{1-N}t^{q^N-1-w}G(t)\varphi^N\widetilde{x}$ and $\tau \widetilde{x}-t^{{\xi}}\varphi \widetilde{x}$ one can equivalently (so that the {\it same} $\nabla$ is generated as before!) stipulate the relations $$a{\widetilde{y}}+\tau t^{\xi}\varphi
{\widetilde{y}}-t^{qn_{{\bf x}}-w-1}G(t)\varphi{{{\bf x}}},\quad\quad{{{\bf x}}}-t^{q+w-n_{{\bf x}}}\varphi{{\widetilde{y}}}$$(with $\widetilde{x}$ substituted by $t^{n_{\bf x}-1}{\bf x}$ for a suitable $n_{\bf x}$, see below). But these relations indeed work well even for $\tau=0$, yielding an irreducible $(\varphi,\Gamma)$-module.

\subsection{$k[\tau][[t]][\varphi]$-modules}

\label{ostertouch}

Let $a\in k^{\times}$ and $\tau\in k$. Let $N, n_{{\bf x}}, n_{{\bf y}}\in{\mathbb Z}_{>0}$ and $w\in{\mathbb Z}$ such that\begin{gather}n_{{\bf y}}+n_{{\bf x}}\le w+1\le{\rm min}\{ q(n_{{\bf y}}-1)+1, q^N, qn_{{\bf x}}\}.\label{floodo}\end{gather}Let $$G(t)\in 1+t^{{\rm max}\{n_{{\bf x}},n_{{\bf y}}\}}k[[t]].$$If $\tau\ne0$ then in
$\widetilde{M}=k[[t]][\varphi](k\widetilde{x}\oplus k\widetilde{y})$
define\begin{align}\widetilde{A}&=a\widetilde{y}+\tau t^{{\xi}}\varphi
\widetilde{y}-\tau^{1-N}t^{q^N-1-w}G(t)\varphi^N\widetilde{x},\notag\\\widetilde{B}&=\tau \widetilde{x}-t^{{\xi}}\varphi
\widetilde{x}.\notag\end{align}Apparantly, for $\tau\ne0$ the $k[[t]][\varphi]$-submodule $\langle \widetilde{A},\widetilde{B},t{\widetilde{x}},t^{n_{{\bf y}}}{\widetilde{y}}\rangle$ of $\widetilde{M}$ generated by $\widetilde{A}$, $\widetilde{B}$, $t{\widetilde{x}}$, $t^{n_{{\bf y}}}{\widetilde{y}}$ has the shape of the submodules studied in subsection \ref{viervor}, and in particular yields an \'{e}tale $(\varphi,\Gamma)$-module if a suitable $\Gamma$-action can be added. The (surprising) point of this section is that this construction remains meaningful even if $\tau$ approaches and attains the value $0$. The trick for this is to pass to suitable alternative defining relations and rather push these to their limit for $\tau=0$.

In ${M}=k[[t]][\varphi](k{{{\bf x}}}\oplus k{{{\bf y}}})$
  define\begin{align}{A}&=a{{\bf y}}+\tau t^{\xi}\varphi
{{\bf y}}-t^{qn_{{\bf x}}-w-1}G(t)\varphi{{{\bf x}}},\notag\\{C}&={{{\bf x}}}-t^{q+w-n_{{\bf x}}}\varphi{{{\bf y}}},\notag\\\overline{A}&=
a{{\bf y}}-t^{qn_{{\bf x}}-w-1}\varphi{{{\bf x}}}.\notag\end{align}Let ${\nabla}=\langle
A,C,t^{n_{{\bf x}}}{{\bf x}},t^{n_{{\bf y}}}{{\bf y}}\rangle$ be the $k[[t]][\varphi]$-submodule of $M$ generated by the elements $A$, $C$, $t^{n_{{\bf x}}}{{\bf x}}$ and $t^{n_{{\bf y}}}{{\bf y}}$. Similarly define $\langle
  \overline{A},{C},t^{n_{{\bf x}}}{{\bf x}},t^{n_{{\bf y}}}{{\bf y}}\rangle$ inside $M$. Put $\Delta=M/{\nabla}$. 

\begin{lem}\label{bendent} (a) If $\tau\ne0$ then the $k[[t]][\varphi]$-linear map\begin{gather}\widetilde{M}\longrightarrow M,\quad\quad \widetilde{x}\mapsto t^{{n_{{\bf x}}}-1}{{\bf x}},\quad \widetilde{y}\mapsto {{\bf y}}\label{justus1}\end{gather}induces an isomorphism $ \langle \widetilde{A},\widetilde{B},t{\widetilde{x}},t^{n_{{\bf y}}}{\widetilde{y}}\rangle\cong {\nabla}$.

  (b) If $\tau=0$ then ${\nabla}=\langle
  \overline{A},{C},t^{n_{{\bf x}}}{{\bf x}},t^{n_{{\bf y}}}{{\bf y}}\rangle$.

(c) View $\tau$ as a free variable. The $k[\tau][[t]]$-module $\Delta^*={\rm Hom}_{k[\tau]}(\Delta,k[\tau])$ is free of rank two.  

\end{lem}

{\sc Proof:} For statement (b) we need to see $A\in \langle
  \overline{A},{C},t^{n_{{\bf x}}}{{\bf x}},t^{n_{{\bf y}}}{{\bf y}}\rangle$ and $\overline{A}\in \nabla$. Using that $G(t)\in 1+t^{n_{{\bf y}}}k[[t]]$ (which implies $(1-G(t))a{{\bf y}}\equiv a{{\bf y}}$ modulo $\langle t^{n_{{\bf y}}}{{\bf y}}\rangle$) and $\tau=0$ we find $A\equiv G(t)\overline{A}$ modulo $\langle t^{n_{{\bf y}}}{{\bf y}}\rangle$, and hence also $G(t)^{-1}A\equiv \overline{A}$ modulo $\langle t^{n_{{\bf y}}}{{\bf y}}\rangle$ (with $G(t)^{-1}\in k[[t]]$). This shows us what we want.

  For statement (a) define$${B}=\tau{{{\bf x}}}-t^{\xi n_{{\bf x}}}\varphi{{{\bf x}}}.$$If $\tau\ne0$ then the map (\ref{justus1}) induces an isomorphism $\langle\widetilde{A},\widetilde{B},t{\widetilde{x}},t^{n_{{\bf y}}}{\widetilde{y}}\rangle\cong\langle {A},{B},t^{n_{{\bf x}}}{{\bf x}},t^{n_{{\bf y}}}{{\bf y}}\rangle$. Thus, we need to see $B\in \nabla$ and $C\in\langle {A},{B},t^{n_{{\bf x}}}{{\bf x}},t^{n_{{\bf y}}}{{\bf y}}\rangle$. But we have$$\tau
    C+t^{w+1-n_{{\bf x}}}A\equiv B\quad\mbox{ modulo }\langle t^{n_{{\bf y}}}{{\bf y}},
    t^{n_{{\bf x}}}{{\bf x}}\rangle$$due to $n_{{\bf y}}+n_{{\bf x}}\le w+1$ (which implies $t^{w+1-n_{{\bf x}}}a{{\bf y}}\equiv 0$ modulo $\langle t^{n_{{\bf y}}}{{\bf y}}\rangle$) and $G(t)\in 1+t^{n_{{\bf x}}}k[[t]]$ (which implies $t^{\xi n_{{\bf x}}}(G(t)-1)\varphi {{\bf x}}\equiv 0$ modulo $\langle t^{n_{{\bf x}}}{{\bf x}}\rangle$). This shows us what we want.

Now let us look at statement (c). Viewing $\tau$ as a free variable means that now we understand ${M}=k[\tau][[t]][\varphi]\otimes_{k[\tau]}V$ with $V=k[\tau]{{\bf x}}\oplus k[\tau]{{\bf y}}$, furthermore $\Delta=M/{\nabla}$ where ${\nabla}=\langle C,A,t^{n_{{\bf x}}}{{\bf x}},t^{n_{{\bf y}}}{{\bf y}}\rangle$ is the $k[\tau][[t]][\varphi]$-submodule of $M$ generated by the elements $C$, $A$, $t^{n_{{\bf x}}}{{\bf x}}$, $t^{n_{{\bf y}}}{{\bf y}}$. Inside
$\Delta$
consider$$W_-=\sum_{s>0}\sum_{0\le\theta<q^{s-1}(qn_{\bf x}-w-1)}t^{\theta}\varphi^sk[\tau]{\bf
  x}+\sum_{s>0}\sum_{0\le\theta<q^{s-1}(q+w-n_{\bf
    x})}t^{\theta}\varphi^sk[\tau]{\bf y},$$$$W_+=(\sum_{i=0}^{n_{\bf
    x}-2}k[\tau]t^i{\bf x})+ (\sum_{i=0}^{n_{\bf y}-2}k[\tau]t^i{\bf
  y}),$$$$W=W_-+W_+.$$We then claim to have the $k[\tau]$-module
decomposition\begin{gather}\Delta= W\bigoplus k[\tau]t^{n_{\bf x}-1}{\bf
  x}\bigoplus k[\tau]t^{n_{\bf y}-1}{\bf
  y}.\label{wassermuehle0georg1steuer00}\end{gather}Indeed, we first look at the
locus where $\tau=0$. Then we may assume $G(t)=1$ (cf. statement (b)), so we
may apply Lemma \ref{fassageo}, as formulae
(\ref{univorso1}) are then equivalent with $n_{{\bf y}}\le w+1$ and $w\le
q(n_{{\bf y}}-1)$. Above the locus where $\tau\ne 0$ we first claim\begin{gather}W[\tau^{-1}]=W'_-+
W_+[\tau^{-1}]\label{josefmai00}\end{gather}with$$W'_-=\sum_{s>0}\sum_{0\le\theta<q^{s-1}\xi}t^{\theta}\varphi^sk[\tau^{\pm}]{\bf
  y}+\sum_{s>0}\sum_{0\le\theta<q^{s-1}\xi n_{\bf
    x}}t^{\theta}\varphi^sk[\tau^{\pm}]{\bf x}.$$Indeed, that $\sum_{s>0}\sum_{0\le\theta<q^{s-1}\xi n_{\bf
    x}}t^{\theta}\varphi^sk[\tau^{\pm}]{\bf x}$ is contained in $W[\tau^{-1}]$
follows by looking at the generators $A$ and $C$ of $\nabla$; thus $W'_-+
W_+[\tau^{-1}]\subset W[\tau^{-1}]$. Conversely, that $\sum_{s>0}\sum_{0\le\theta<q^{s-1}(q+w-n_{\bf
    x})}t^{\theta}\varphi^sk[\tau^{\pm}]{\bf y}$ is contained in $W'_-+
W_+[\tau^{-1}]$ follows by looking at the generators $A$ and $B$ of $\nabla$ (with $B$ as defined in the proof of
statement (a)); thus $ W[\tau^{-1}]\subset W'_-+
W_+[\tau^{-1}]$ and the claim is proven. It implies that the decomposition
(\ref{wassermuehle0georg1steuer00}) over $k[\tau^{\pm}]$ follows from statement (a) and Lemma
\ref{fassa}, as formula (\ref{univorso}) is then equivalent with the
assumption $w\le q(n_{{\bf y}}-1)$.

In view of the decomposition (\ref{wassermuehle0georg1steuer00}) we may now
define $\lambda_{\bf x},\lambda_{\bf y}\in \Delta^*$ by
requiring $$\lambda_{\bf x}(t^{n_{\bf x}-1}{\bf
  x})=\lambda_{\bf y}(t^{n_{\bf y}-1}{\bf
  y})=1,\quad \quad \lambda_{\bf x}(t^{n_{\bf y}-1}{\bf
  y})=\lambda_{\bf y}(t^{n_{\bf x}-1}{\bf
  x})=\lambda_{\bf x}(W)=\lambda_{\bf y}(W)=0.$$We claim that $\lambda_{\bf
  x},\lambda_{\bf y}$ is a $k[\tau][[t]]$-basis of $\Delta^*$. Indeed, that $\lambda_{\bf
  x},\lambda_{\bf y}$ are linearly independent may be checked after passing to
$k[\tau^{\pm}][[t]]$, as everything in sight is $k[\tau]$-torsion free. It
then follows as in Lemma \ref{fassa}, which, in its version over
$k[\tau^{\pm}]$ rather than $k$, we may invoke due to statement (a) and
formula (\ref{josefmai}). To see that $\lambda_{\bf
  x},\lambda_{\bf y}$ generate $\Delta^*$ notice first that a general element in $\Delta^*$ is a (possibly infinite)
linear combination of $k[\tau]$-linear forms $\lambda$ with support killed by a
power of $t$. It is
enough to see that for any such $\lambda$ we find $\alpha_{\bf
  x},\alpha_{\bf
  y}\in k[\tau][[t]]$ with $\lambda= \alpha_{\bf
  x}\lambda_{\bf
  x}+\alpha_{\bf
  y}\lambda_{\bf
  y}$. We use statement (a) and
formula (\ref{josefmai00}). They allow us to invoke Lemma \ref{fassa}, in its
version over $k[\tau^{\pm}]$ rather than $k$, to provide $\alpha_{\bf
  x},\alpha_{\bf
  y}\in k[\tau^{\pm}][[t]]$ with $\lambda= \alpha_{\bf
  x}\lambda_{\bf
  x}+\alpha_{\bf
  y}\lambda_{\bf
  y}$. But in fact we have $\alpha_{\bf
  x},\alpha_{\bf
  y}\in k[\tau][[t]]$: this follows by
looking at the fibre $\tau=0$, employing Lemma \ref{fassageo}.\hfill$\Box$\\

Let $\widetilde{u}\in[0,q-2]$. Suppose that we have an action of $\Gamma$ on
$\Delta$ such that$$\gamma*{{\bf x}}=a^{\widetilde{u}}_{\gamma,1}{{\bf x}},\quad\quad\gamma*{{\bf y}}=a^{\widetilde{u}}_{\gamma,1}{{\bf y}}$$for all $\gamma\in\Gamma$ with $[\gamma]=a_{\gamma,1}t$. Put $$\overline{\bf D}=\Delta^*\otimes_{k[[t]]}k((t))\quad\quad\mbox{ for }\tau=0.$$

If $\tau=0$ then it follows from Lemma \ref{fassageo} that $\overline{\bf D}$ has a $k((t))$-basis\footnote{Warning: These $\lambda_{{\bf x}}$, $\lambda_{{\bf y}}$ need not be those of the proof of Lemma \ref{bendent}.} $\lambda_{{\bf x}}, \lambda_{{\bf y}}$ for
which $$\varphi(\lambda_{{\bf x}})=t^{q(n_{{\bf x}}-1)-w}\lambda_{{\bf
    y}},\quad\quad \varphi(\lambda_{{\bf y}})=a^{-1}t^{1+w-n_{{\bf
      x}}}\lambda_{{\bf x}},\quad\quad\gamma*\lambda_{{\bf
    y}}=a^{-\widetilde{u}}_{\gamma,1}\lambda_{{\bf y}}$$for all
$\gamma\in\Gamma$ with $[\gamma]=a_{\gamma,1}t$. In particular we see
$\varphi^2(\lambda_{{\bf y}})=a^{-1}t^{{\xi} w}\lambda_{{\bf y}}$. It follows
that, putting $$h=-w-(q+1)\widetilde{u}$$and
$\lambda'=t^{\widetilde{u}}\lambda_{{\bf y}}$ we
have$$\varphi^2(\lambda')=a^{-1}t^{-\xi
  h}\lambda',$$$$\gamma*\lambda'=\lambda'$$for all $\gamma\in\Gamma$ with
$[\gamma]=a_{\gamma,1}t$. With Lemma \ref{endosterkreis} we see that if $q+1$
does not divide $w$ then $\overline{\bf D}\in{\mathcal V}(h)$.

\subsection{Example}

\label{uebgutD}

Let $\widetilde{\ell}\in[1,q-2]$ and $i\in{\mathbb Z}_{\ge f}$. Put\footnote{For the definition of $\widetilde{\ell}_{[i]}$ see formula \ref{movopfin}.}
\begin{gather}w=p^i\xi-\widetilde{\ell}_{[i]},\label{saanruf1}\end{gather}$$n_{{\bf x}}=\frac{w+q}{q},\quad\quad n_{{\bf y}}=w+1-n_{{\bf x}}=\frac{w \xi}{q}$$and $G(t)=1$. It is straightforward to verify formula (\ref{floodo}). Let $\Gamma$ act on $M$ by means
of$$\gamma*{{\bf y}}=a^{\widetilde{u}+1}_{\gamma,1}\frac{t}{[\gamma]}{{\bf y}},\quad\quad\gamma*{{\bf x}}=a^{\widetilde{u}-\widetilde{\ell}+1}_{\gamma,1}(\frac{t}{[\gamma]})^{n_{{\bf x}}}{{\bf x}}$$for
$\gamma\in\Gamma$. To see that the $\Gamma$-action on $M$ preserves $\nabla$ we separately treat the cases $\tau\ne0$ and $\tau=0$.

If $\tau\ne0$ we use Proposition \ref{vacet10}. We translate from the present
setting to the setting in Proposition \ref{vacet10} (with $x$, $n_x$, $y$, $n_y$, ${\mathcal D}$, ${\mathcal E}$,
$\ell$, $u$ as in Proposition \ref{vacet10}) by means of the
assignments$$x=\widetilde{y},\quad\quad n_x=n_{{\bf y}},\quad\quad
y=\widetilde{x},\quad\quad n_y=1,$$$${\mathcal E}=\emptyset,\quad\quad{\mathcal D}=\{i\},\quad\quad
\ell=\widetilde{\ell},\quad\quad u=\widetilde{u}.$$Condition (\ref{sommerdornwald5}) then becomes $\widetilde{\ell}_{[i]}\ge1$ (which clearly holds true), and condition (\ref{marufmin00}) then becomes $q(n_{{\bf y}}-1)\ge w$ (which follows from the definitions). 

If $\tau=0$ we use Lemma \ref{fassageo}. We translate from the present setting
to the setting in Lemma \ref{fassageo} (with $x$, $n_x$, $y$, $n_y$, $\nu$,
$\mu$, $u$ as in Lemma \ref{fassageo}) by means of the assignments$$y={\bf
  y},\quad\quad n_y=n_{{\bf y}},\quad\quad x={\bf x},\quad\quad n_x=n_{{\bf
    x}},$$$$\nu=q n_x-w-1=\xi,\quad\quad\quad\quad
\mu=q+w-n_x=\frac{\xi}{q}(q+w)$$and $u=\widetilde{u}$. We observe that
$p^{i-f}$ divides $\mu-\xi=n_{{\bf y}}$, and as $p^{i-f}\xi\ge\frac{w}{q}+1$
this implies that  condition (\ref{james0}) is valid. For condition
(\ref{james1}) this is trivially the case. The conditions $\nu\le qn_x-n_y$
and $\mu\le qn_y-n_x$ are easily verified.

As $q+1$ does not divide $w$ we see that the corresponding ${\mathcal G}_F$-representation is irreducible.

\subsection{Example}

\label{uebgutE}

Let $\widetilde{\ell}\in[1,q-2]$ and $i\in{\mathcal E}(\widetilde{\ell})$. Write $\widetilde{\ell}=\sum_{i=0}^{f-1}m_ip^i$ with $0\le m_i<p$. Put\begin{gather}w=2p^{2f+i}\xi-q^2\widetilde{\ell}_{[i]},\label{saanruf}\end{gather}
$$n_{\bf y}=p^{i+f}\xi^2,$$$$n_{{\bf
    x}}=1+p^{i+f}(2q+m_i)-q{\widetilde{\ell}_{[i]}}.$$It is straightforward to verify formula (\ref{floodo}). Put $$G(t)=H(m_i,i)$$with $H(m_i,i)$ as in subsection \ref{viervor}. Notice that $$n_{{\bf x}}\le p^{i+f}\xi^2,\quad\quad n_{{\bf y}}\le p^{i+f}\xi^2,$$and hence $G(t)\in 1+t^{{\rm max}\{n_{{\bf x}},n_{{\bf y}}\}}k[[t]]$ since $p^ig^{(1)}={{n}}_i^{(1)}-{{n}}_i^{(0)}=p^{i+f}\xi^2$ (formula (\ref{wiederlegaut})).

Let $\Gamma$ act on $M$ by means
of$$\gamma*{{\bf y}}=a^{\widetilde{u}+1}_{\gamma,1}\frac{t}{[\gamma]}{{\bf y}},\quad\quad\gamma*{{\bf x}}=a^{\widetilde{u}-\widetilde{\ell}+1}_{\gamma,1}(\frac{t}{[\gamma]})^{n_{{\bf x}}}{{\bf x}}$$for
$\gamma\in\Gamma$. To see that the $\Gamma$-action on $M$ preserves $\nabla$ we separately treat the cases $\tau\ne0$ and $\tau=0$.

If $\tau\ne0$ we use Proposition \ref{vacet10}. We translate from the present setting to the setting in Proposition \ref{vacet10} (with $x$, $n_x$, $y$, $n_y$, ${\mathcal D}$, ${\mathcal E}$,
$\ell$, $u$ as in Proposition \ref{vacet10}) by means of the
assignments$$x=\widetilde{y},\quad\quad n_x=n_{{\bf y}},\quad\quad
y=\widetilde{x},\quad\quad n_y=1,$$$${\mathcal
  D}=\emptyset,\quad\quad{\mathcal E}=\{i\},\quad\quad
\ell=\widetilde{\ell},\quad\quad u=\widetilde{u}.$$Condition
(\ref{fuhleresel200}) holds true by Lemma \ref{letzterst}. Condition
(\ref{maruf200}) is equivalent, after inserting the present definitions,
with $(q-3)\xi p^{i+2f}+q^2\widetilde{\ell}_{[i]}\ge q$, which indeed is true.

If $\tau=0$ we use Lemma \ref{fassageo}. We translate from the present setting
to the setting in Lemma \ref{fassageo} (with $x$, $n_x$, $y$, $n_y$, $\nu$,
$\mu$, $u$ as in Lemma \ref{fassageo}) by means of the assignments$$y={\bf
  y},\quad\quad n_y=n_{{\bf y}},\quad\quad x={\bf x},\quad\quad n_x=n_{{\bf
    x}},$$$$\nu=q n_x-w-1,\quad\quad\quad\quad
\mu=q+w-n_x,\quad\quad\quad\quad u=\widetilde{u}.$$Notice that indeed
$\mu,\nu>0$. We observe that $p^{i+f+1}$
divides $\mu-\xi$, and as $p^{i+f+1}\xi\ge n_x$ this implies that condition (\ref{james0}) is valid. Similarly, $p^{i+2f}$
divides $\nu-\xi$, and as $p^{i+2f}\xi\ge n_y$ this implies that condition (\ref{james1}) is valid. The conditions $n_{{\bf y}}+n_{{\bf x}}-1\le w$ and $\nu\le qn_x-n_y$ and $\mu\le qn_y-n_x$ are easily verified. 

We distinguish two cases (cf. Lemma \ref{endosterkreis}):

If $\widetilde{\ell}_{[i]}\ne p^i(q-3)$ then $\overline{\bf D}\in{\mathcal
  V}(h)$. If $\widetilde{\ell}_{[i]}=p^i(q-3)$ then $\overline{\bf D}$ is the
direct sum of two isomorphic $(\varphi,\Gamma)$-modules of rank one,
and the restriction to ${\mathcal I}_F$ of the corresponding ${\mathcal
  G}_F$-representation is isomorphic with $\omega^{-p^i-\widetilde{u}}\bigoplus\omega^{-p^i-\widetilde{u}}$.

\subsection{Reformulation}

\label{lichtherzkarsa}

Recall the set ${\mathbb E}$ defined in subsection \ref{stackstrataset}.

\begin{kor}\label{cusanus00} Let $(\widetilde{\ell}, \widetilde{u},\widetilde{\mathcal I})\in{\mathbb E}$ with $\widetilde{\ell}\ne0$, let $\overline{\bf D}$ be a $(\varphi,\Gamma)$-module and assume that we are in one of the following situations (a) or (b):

  (a) There is some $h\in{\mathbb Z}$ with $\overline{\bf D}\in{\mathcal V}(h)$ and some $i\ge0$ with $\widetilde{\mathcal
  I}=\{\Pi(i)\}$ such that either (a1) or (a2) holds:

(a1) $\Pi(i)\in{\mathcal D}(\widetilde{\ell})$ and $$h\equiv\widetilde{\ell}_{[i]}-p^i\xi-(q+1)\widetilde{u}\quad\mbox{
       modulo }(q^2-1){\mathbb Z}.$$    

     (a2) $\Pi(i)\in{\mathcal E}(\widetilde{\ell})$ and $\widetilde{\ell}_{[\Pi(i)]}\ne p^{\Pi(i)}(q-3)$ and$$h\equiv\widetilde{\ell}_{[i]}-2p^i\xi-(q+1)\widetilde{u}\quad\mbox{ modulo }(q^2-1){\mathbb Z}.$$

(b) There is some $(\overline{\ell}, \overline{u},\overline{\mathcal I})\in{\mathbb E}$ with $\overline{\bf D}\in{\mathcal V}(\overline{\ell},\overline{u},\overline{\mathcal
  I})$ and some $i\ge 2f$ with $\widetilde{\ell}_{[i]}=p^i(q-3)$ and $\overline{\ell}=0$ and $\overline{u}\equiv p^i+\widetilde{u}$ modulo $\xi{\mathbb Z}$ and $\overline{\mathcal I}=\emptyset$ and $\widetilde{\mathcal I}=\{\Pi(i)\}$.

Then there exists a $(\varphi,\Gamma)$-module ${\bf D}$ over $k[\tau]$ with
  ${\bf D}\otimes_{k[\tau]}k\cong \overline{\bf D}$ via $\tau\mapsto 0\in k$, and
  with ${\bf D}\otimes_{k[\tau]}k\in{\mathcal V}(\widetilde{\ell},\widetilde{u},\widetilde{\mathcal
  I})$ via $\tau\mapsto a$ for each $a\in k^{\times}$.

  \end{kor}

  {\sc Proof:} As far as case (a1) is concerned, this is the upshot of subsection \ref{uebgutD}. Namely, choose $\overline{A}$, $C$ as at the beginning of subsection \ref{ostertouch}, and $n_{\bf x}$, $n_{\bf y}$ as in subsection \ref{uebgutD}, such that we have an isomorphism of $(\varphi,\Gamma)$-modules $$\overline{\bf D}\cong (k[[t]][\varphi](k{\bf x}\oplus k{\bf y})/\langle
  \overline{A},{C},t^{n_{{\bf x}}}{{\bf x}},t^{n_{{\bf y}}}{{\bf y}}\rangle)^*\otimes_{k[[t]]}k((t)).$$Use the coefficients of $\overline{A}$ to define $A$ as at the beginning of subsection \ref{ostertouch}. Then, by Lemma \ref{bendent} and the discussion of subsection \ref{uebgutD},$${\bf D}=({\rm Hom}_{k[\tau]}(\Delta,k[\tau]))\otimes_{k[\tau][[t]]}k[\tau]((t))$$with$$\Delta=k[\tau][[t]][\varphi](k{\bf x}\oplus k{\bf y})/\langle{A},{C},t^{n_{{\bf x}}}{{\bf x}},t^{n_{{\bf y}}}{{\bf y}}\rangle$$carries the structure of a $(\varphi,\Gamma)$-module over $k[\tau]$ as desired. 

In the cases (a2) and (b) use the same argument, now employing subsection \ref{uebgutE} rather than subsection \ref{uebgutD}.\hfill$\Box$\\

  \section{Degeneration into reducibles with changing $\ell$}

\label{zoomgergelyfengwei}

From this point on until the end of this paper we assume $p>2$.

The purpose of this section is to prove Corollary \ref{cusanus}. The strategy is very similar to the one explained at the beginning of subsection \ref{ostertouch}.

\subsection{$k[\tau][[t]][\varphi]$-modules}

\label{osterfreitag}

Let $N, n_{{\bf x}}, n_{{\bf y}},w\in{\mathbb
  Z}_{>0}$ with $q^N-1\ge{\max}\{q(2w+1-n_{{\bf x}}),q(n_{{\bf x}}-1) \}$ and assume \begin{gather}{\xi}(2w+1)\ge qn_{{\bf x}},\label{voklau1}\\{\xi}(n_{{\bf x}}-1-w)\ge n_{{\bf y}},\label{voklau2}\\n_{{\bf y}}\ge n_{{\bf x}}.\label{voklau3}\end{gather}Suppose we are also given $a, b\in k^{\times}$ and $\tau\in k$, as well as \begin{align}F(t)&\in 1+t^{qn_{{\bf x}}-{\xi}(w+1)}k[[t]],\notag\\G(t)&\in 1+t^{{\rm max}\{n_{{\bf x}},n_{{\bf y}}-w\}}k[[t]].\notag\end{align}If $\tau\ne0$ then in
$\widetilde{M}=k[[t]][\varphi](k\widetilde{x}\oplus k\widetilde{y})$
define\begin{align}\widetilde{A}&=a\widetilde{y}-\tau t^{\xi}\varphi\widetilde{y}-(b\tau)^{1-N}t^{q^N-1+\xi w-q(n_{{\bf x}}-1)}G(t)
\varphi^N\widetilde{x}-{\tau}(b\tau)^{1-N}t^{q^N-1-q(n_{{\bf x}}-1)}F(t)\varphi^N\widetilde{x},\label{vierwoend}\\\widetilde{B}&=b\tau \widetilde{x}-t^{\xi}\varphi
\widetilde{x}.\notag\end{align}Just as in section \ref{okt3}: For $\tau\ne0$ the $k[[t]][\varphi]$-submodule $\langle \widetilde{A},\widetilde{B},t{\widetilde{x}},t^{n_{{\bf y}}}{\widetilde{y}}\rangle$ of $\widetilde{M}$ generated by $\widetilde{A}$, $\widetilde{B}$, $t{\widetilde{x}}$, $t^{n_{{\bf y}}}{\widetilde{y}}$ apparantly has the shape of the submodules studied in subsection \ref{viervor}, and in particular yields an \'{e}tale $(\varphi,\Gamma)$-module if a suitable $\Gamma$-action can be added. The (surprising) point of this section is that this construction remains meaningful even if $\tau$ approaches and attains the value $0$. The trick for this is to pass to suitable alternative defining relations and rather push these to their limit for $\tau=0$. The core of the present section is really just Lemma \ref{pitpaul} (which still ignores $\Gamma$-actions). It is applied in subsections \ref{keydefo}, \ref{keydefo1} where the principal technical difficultis arise in ensuring the list of inequalities (\ref{sommerdornwald5}), (\ref{marufmin00}), (\ref{fuhleresel200}), (\ref{maruf200}), (\ref{sommerdorne}), (\ref{sommerdornwald1}), (\ref{sommerdornwald6}), (\ref{maruf2000}) both generically and for $\tau=0$, and in handling the limiting behaviour of the exponents $r, o_i, r_i, n_i^{(j)}$ of subsection \ref{viervor}.

In ${M}=k[[t]][\varphi](k{{{\bf x}}}\oplus k{{{\bf y}}})$
define\begin{align}C&=a{{\bf y}}-\tau t^{{\xi}}\varphi {{\bf y}}-t^{{\xi}(w+1)}G(t)\varphi
{{\bf y}}+{\tau}t^{{\xi}}F(t)\varphi {{\bf x}},\notag\\D&=b{{\bf x}}-t^{{\xi}(n_{{\bf x}}-w)}\varphi {{\bf x}}+t^{{\xi}(n_{{\bf x}}-w)}F^{-1}(t)\varphi
{{\bf y}}.\notag\end{align}In
$\overline{M}=k[[t]][\varphi](k\overline{x}\oplus k\overline{y})$
define\begin{align}\overline{C}&=a\overline{y}-t^{\xi}\varphi
\overline{y},\notag\\\overline{D}&=b\overline{x}-t^{\xi}\varphi
                                   \overline{x}-a^{1-N}t^{q^N-1+q(n_{{\bf x}}-2w-1)}F^{-1}(t)\varphi^N\overline{y}.\label{frankvirsik}\end{align}Let ${\nabla}=\langle C,D,t^{n_{{\bf x}}}{{\bf x}},t^{n_{{\bf y}}}{{\bf y}}\rangle$ be the $k[[t]][\varphi]$-submodule of $M$ generated by the elements $C$, $D$, $t^{n_{{\bf x}}}{{\bf x}}$ and $t^{n_{{\bf y}}}{{\bf y}}$. Similarly define $\langle\overline{C}, \overline{D}, t^{w+1}\overline{x},t^{n_{\bf y}-w}\overline{y}\rangle$ inside $\overline{M}$. Put $\Delta=M/{\nabla}$.

\begin{lem}\label{pitpaul} (a) If $\tau\ne0$ then the $k[[t]][\varphi]$-linear map\begin{gather}\widetilde{M}\longrightarrow M,\quad\quad \widetilde{x}\mapsto t^{{n_{{\bf x}}}-1}{{\bf x}},\quad \widetilde{y}\mapsto {{\bf y}}\label{justus}\end{gather}induces an isomorphism $\langle\widetilde{A},\widetilde{B},t\widetilde{x},t^{n_{\bf y}}\widetilde{y}\rangle\cong\nabla$.

  (b) If $\tau=0$ then the $k[[t]][\varphi]$-linear map$$\overline{M}\longrightarrow M,\quad\quad \overline{x}\mapsto t^{{n_{{\bf x}}}-w-1}{{\bf x}},\quad \overline{y}\mapsto t^w{{\bf y}}$$ induces an isomorphism $\langle\overline{C}, \overline{D}, t^{w+1}\overline{x},t^{n_{\bf y}-w}\overline{y}\rangle\cong\nabla$.

(c) View $\tau$ as a free variable. The $k[\tau][[t]]$-module $\Delta^*={\rm Hom}_{k[\tau]}(\Delta,k[\tau])$ is free of rank two.
\end{lem}

{\sc Proof:} As $G(t)\in1+t^{n_{\bf y}-w}k[[t]]$ we have $G(t)\overline{C}\equiv a\overline{y}-t^{\xi}G(t)\varphi\overline{y}$ modulo $\langle t^{n_{\bf y}-w}\overline{y}\rangle$, thus statement (b) is clear. Define\begin{align}{A}&=a{{\bf y}}-\tau t^{\xi}\varphi {{\bf y}}-t^{{\xi}(w+1)}G(t)\varphi
{{\bf x}}+{\tau}t^{{\xi}}F(t)\varphi {{\bf x}},\notag\\B&=b\tau {{\bf x}}-t^{{\xi}n_{{\bf x}}}\varphi
{{\bf x}}.\notag\end{align}If $\tau\ne0$ then the map (\ref{justus}) induces an
isomorphism $\langle\widetilde{A},\widetilde{B},t\widetilde{x},t^{n_{\bf y}}\rangle\cong\langle
{A},{B},t^{n_{{\bf x}}}{{\bf x}},t^{n_{{\bf y}}}{{\bf y}}\rangle$. Thus, we need to see $A, B\in\nabla$ and $C,D\in\langle
{A},{B},t^{n_{{\bf x}}}{{\bf x}},t^{n_{{\bf y}}}{{\bf y}}\rangle$.

Using formula (\ref{voklau1}) and that $F(t)\in 1+t^{qn_{{\bf
      x}}-{\xi}(w+1)}k[[t]]$ (notice that
$n_{{\bf x}}\le qn_{{\bf x}}-{\xi}(w+1)$) we see \begin{gather}A\equiv
  C+t^{{\xi}(2w+1-n_{{\bf x}})}G(t)F(t)D\quad\mbox{ modulo }\langle t^{n_{{\bf y}}}{{\bf y}},
  t^{n_{{\bf x}}}{{\bf x}}\rangle.\label{maisonne1}\end{gather}Using formula
(\ref{voklau2}) and that $F(t)\in 1+t^{n_{{\bf x}}}k[[t]]$ (notice that
$n_{{\bf x}}\le qn_{{\bf x}}-{\xi}(w+1)$) and $G(t)\in1+t^{n_{\bf x}}k[[t]]$ we see \begin{gather}\tau D+ t^{{\xi}(n_{{\bf x}}-1-w)}F^{-1}(t)A\equiv B\quad\mbox{ modulo }\langle t^{n_{{\bf y}}}{{\bf y}},
  t^{n_{{\bf x}}}{{\bf x}}\rangle.\label{maisonne2}\end{gather}Combining formulae (\ref{maisonne1}) and (\ref{maisonne2}) we get what we wanted.

Now let us look at statement (c). Viewing $\tau$ as a free variable means that now
we understand ${M}=k[\tau][[t]][\varphi]\otimes_{k[\tau]}V$ with
$V=k[\tau]{{\bf x}}\oplus k[\tau]{{\bf y}}$, furthermore $\Delta=M/{\nabla}$
where ${\nabla}=\langle C,D,t^{n_{{\bf x}}}{{\bf x}},t^{n_{{\bf y}}}{{\bf
    y}}\rangle$ is the $k[\tau][[t]][\varphi]$-submodule of $M$ generated by
the elements $C$, $D$, $t^{n_{{\bf x}}}{{\bf x}}$, $t^{n_{{\bf y}}}{{\bf
    y}}$. Inside $\Delta$
consider$$W_-=\sum_{s>0}\sum_{0\le\theta<q^{s-1}\xi(w+1)}t^{\theta}\varphi^sk[\tau]{\bf
  y}+\sum_{s>0}\sum_{0\le\theta<q^{s-1}\xi(n_{\bf
    x}-w)}t^{\theta}\varphi^sk[\tau]{\bf x},$$$$W_+=(\sum_{i=0}^{n_{\bf
    x}-2}k[\tau]t^i{\bf x})+ (\sum_{i=0}^{n_{\bf y}-2}k[\tau]t^i{\bf
  y}),$$$$W=W_-+ W_+.$$We then claim to have the $k[\tau]$-module
decomposition\begin{gather}\Delta= W\bigoplus k[\tau]t^{n_{\bf x}-1}{\bf
  x}\bigoplus k[\tau]t^{n_{\bf y}-1}{\bf
  y}.\label{wassermuehle0georg1steuer}\end{gather}Indeed, we first look at the
locus where $\tau=0$. Then we may assume $G(t)=1$ (cf. statement (b)), so we may apply Lemma \ref{fassa}, as formula
(\ref{univorso}) is then equivalent with $(w+1)\xi+n_{\bf y}\le q n_{\bf x}$, which is
implied by hypothesis (\ref{voklau2}). Above the locus where $\tau\ne 0$ we
first claim\begin{gather}W[\tau^{-1}]=W'_-+
W_+[\tau^{-1}]\label{josefmai}\end{gather}with$$W'_-=\sum_{s>0}\sum_{0\le\theta<q^{s-1}\xi}t^{\theta}\varphi^sk[\tau^{\pm}]{\bf
  y}+\sum_{s>0}\sum_{0\le\theta<q^{s-1}\xi n_{\bf
    x}}t^{\theta}\varphi^sk[\tau^{\pm}]{\bf x}.$$Indeed, that $\sum_{s>0}\sum_{0\le\theta<q^{s-1}\xi n_{\bf
    x}}t^{\theta}\varphi^sk[\tau^{\pm}]{\bf x}$ is contained in $W[\tau^{-1}]$
follows by looking at the generators $C$ and $D$ of $\nabla$; thus $W'_-+
W_+[\tau^{-1}]\subset W[\tau^{-1}]$. Conversely, that
$\sum_{s>0}\sum_{0\le\theta<q^{s-1}\xi(w+1)}t^{\theta}\varphi^sk[\tau^{\pm}]{\bf
  y}$ is contained in $W'_-+
W_+[\tau^{-1}]$ follows by looking at $A$ and $B$ as defined in the proof of
statement (a); thus $ W[\tau^{-1}]\subset W'_-+
W_+[\tau^{-1}]$ and the claim is proven. It implies that the decomposition
(\ref{wassermuehle0georg1steuer}) over $k[\tau^{\pm}]$ follows from statement (a) and Lemma
\ref{fassa}, as formula (\ref{univorso}) is then equivalent with $w\le
q(n_{{\bf y}}-1)$, which is clearly satisfied.

In view of the decomposition (\ref{wassermuehle0georg1steuer}) we may now
define $\lambda_{\bf x},\lambda_{\bf y}\in \Delta^*$ by
requiring $$\lambda_{\bf x}(t^{n_{\bf x}-1}{\bf
  x})=\lambda_{\bf y}(t^{n_{\bf y}-1}{\bf
  y})=1,\quad \quad \lambda_{\bf x}(t^{n_{\bf y}-1}{\bf
  y})=\lambda_{\bf y}(t^{n_{\bf x}-1}{\bf
  x})=\lambda_{\bf x}(W)=\lambda_{\bf y}(W)=0.$$We claim that $\lambda_{\bf
  x},\lambda_{\bf y}$ is a $k[\tau][[t]]$-basis of $\Delta^*$. Indeed, that $\lambda_{\bf
  x},\lambda_{\bf y}$ are linearly independent may be checked after passing to
$k[\tau^{\pm}][[t]]$, as everything in sight is $k[\tau]$-torsion free. It
then follows as in Lemma \ref{fassa}, which, in its version over
$k[\tau^{\pm}]$ rather than $k$, we may invoke due to statement (a) and
formula (\ref{josefmai}). To see that $\lambda_{\bf
  x},\lambda_{\bf y}$ generate $\Delta^*$ notice first that a general element in $\Delta^*$ is a (possibly infinite)
linear combination of $k[\tau]$-linear forms $\lambda$ with support killed by a
power of $t$. It is
enough to see that for any such $\lambda$ we find $\alpha_{\bf
  x},\alpha_{\bf
  y}\in k[\tau][[t]]$ with $\lambda= \alpha_{\bf
  x}\lambda_{\bf
  x}+\alpha_{\bf
  y}\lambda_{\bf
  y}$. We use statement (a) and
formula (\ref{josefmai}). They allow us to invoke Lemma \ref{fassa}, in its
version over $k[\tau^{\pm}]$ rather than $k$, to provide $\alpha_{\bf
  x},\alpha_{\bf
  y}\in k[\tau^{\pm}][[t]]$ with $\lambda= \alpha_{\bf
  x}\lambda_{\bf
  x}+\alpha_{\bf
  y}\lambda_{\bf
  y}$. But in fact we have $\alpha_{\bf
  x},\alpha_{\bf
  y}\in k[\tau][[t]]$: this follows by
looking at the fibre $\tau=0$, again employing Lemma \ref{fassa}.\hfill$\Box$\\

\subsection{Example}

\label{keydefo}

Keeping $a,b\in k^{\times}$ and $\tau\in k$ as before, we now consider specific choices of $n_{\bf x}$, $n_{\bf y}$, $w$, $F(t)$ and $G(t)$. Recall formula (\ref{unionbvb}) (and the definition of ${\ell}_{[i]}$ for $\ell\in[0,q-2]$, formula (\ref{movopfin})).

\begin{lem}\label{hanny} Let $\ell\in[0,q-2]$. For integers $0\le j_1<j_2$ we have$$\sigma^{{\ell}}(j_2)\xi
  p^{j_2}-{\ell}_{[j_2]}>2(\sigma^{{\ell}}(j_1)\xi
  p^{j_1}-{\ell}_{[j_1]}).$$
\end{lem}

{\sc Proof:} Written out, with $\ell=\sum_{i=0}^{f-1}m_ip^i$ as usual, the statement becomes$$\sum_{j=j_1+f}^{j_2+f-1}(\sigma^{{\ell}}(j_2)(p-1)-m_j)p^j>\sum_{j=j_2}^{j_1+f-1}((2\sigma^{{\ell}}(j_1)-\sigma^{{\ell}}(j_2))(p-1)-m_j)p^j+\sum_{j=j_1}^{j_2-1}2(\sigma^{{\ell}}(j_1)(p-1)-m_j)p^j.$$If $\sigma^{{\ell}}(j_2)=2$ and $\sigma^{{\ell}}(j_1)=1$ this is clear.

If $\sigma^{{\ell}}(j_2)=2=\sigma^{{\ell}}(j_1)$ the $j=j_1+f$-summand on the left hand side is $\ge(p-1)p^{j_1+f}$, and this exceeds the right hand side: Indeed, $\sigma^{{\ell}}(j_1)=2$ implies $m_{j_1-1}\ge p-2$, so the $j=j_1+f-1$-summand on the right hand side is $\le p^{j_1+f}$; the remaining summands on the right hand side sum up to $<2p^{j_1+f-1}$, and as $p>2$ the claim follows. 

If $\sigma^{{\ell}}(j_2)=1=\sigma^{{\ell}}(j_1)$ then $m_{j_2-1}\le p-2$ and hence the $j=j_2+f-1$-summand on the left hand side is $\ge p^{j_2+f-1}$, and this exceeds the right hand side.

If $\sigma^{{\ell}}(j_2)=1$ and $\sigma^{{\ell}}(j_1)=2$ then there is some $j_1\le i<j_2$ with $m_{i}\le p-3$ and hence the $j=i+f$-summand on the left hand side is $\ge 2 p^{i+f}$, and this exceeds the right hand side: Indeed, using $\sigma^{{\ell}}(j_1)=2$ and hence $m_{j_1-1}\ge p-2$ again, the $j=j_1+f-1$-summand on the right hand side is $\le p^{j_1+f}$; the remaining summands on the right hand side sum up to $<3p^{j_1+f-1}$, and as $p>2$ the claim follows.\hfill$\Box$\\

Fix integers $2f\le i_1<
i_2< i_1+f$ and a set $\widetilde{\mathcal I}$ with $\Pi(\{i_1,i_2\})\subset
\widetilde{\mathcal I}\subset \Pi([i_1,i_2])$. For $\ell\in{\mathbb Z}$
define$${\mathcal J}^{\ell}_{\mathcal
  D}=[i_1+1,i_2-1]\cap \Pi^{-1}(\widetilde{\mathcal I}\cap{\mathcal
  D}({\ell})),$$$${\mathcal J}^{\ell}_{\mathcal E}=[i_1+1,i_2-1]\cap
\Pi^{-1}(\widetilde{\mathcal I}\cap{\mathcal E}({\ell})),$$$${\mathcal J}={\mathcal J}^{\ell}_{\mathcal
  D}\cup{\mathcal J}^{\ell}_{\mathcal
  E}=[i_1+1,i_2-1]\cap
\Pi^{-1}(\widetilde{\mathcal I}).$$Fix $\widetilde{\ell}\in[1,q-2]$ and write $\widetilde{\ell}=\sum_{i=0}^{f-1}m_ip^i$ with $0\le m_i\le p-1$, write $m_i=m_{i+f}$ for all $i\in{\mathbb Z}$. For $i_1\le
i\le i_2$ and $\sigma\in\{1,2\}$ put$$e_{\sigma}(i)=\sigma^{\widetilde{\ell}}(i_2){\xi}p^{i_2}-\widetilde{\ell}_{[i_2]}-\sigma\xi
p^{i}+\widetilde{\ell}_{[i]}.$$Let $\{\alpha_{i}\}_{i\in{\mathcal J}}$ be scalars in $k$ and define\begin{align}\check{F}(t)&=\sum_{i\in{\mathcal J}^{\widetilde{\ell}}_{\mathcal
    D}}\alpha_{i}t^{e_1(i)}+\sum_{i\in{\mathcal J}^{\widetilde{\ell}}_{\mathcal
    E}}\alpha_{i}t^{e_2(i)}H(m_i,i-2f),\notag\\F(t)&=\left\{\begin{array}{l@{\quad:\quad}l}
\check{F}(t)+1&\sigma^{\widetilde{\ell}}(i_2)=1\\\check{F}(t)+H(m_{i_2},i_2-2f)&
\sigma^{\widetilde{\ell}}(i_2)=2\end{array}\right.,\notag\\G(t)&=\left\{\begin{array}{l@{\quad:\quad}l}
1&\sigma^{\widetilde{\ell}}(i_1)=1\\H(m_{i_1},i_1-2f)&\sigma^{\widetilde{\ell}}(i_1)=2\end{array}\right..\notag\end{align}Define the integers$$n_{\bf y}=n_{\bf
  x}=1+\frac{\sigma^{\widetilde{\ell}}(i_2)\xi p^{i_2}-\widetilde{\ell}_{[i_2]}}{q},$$$$w=\sigma^{\widetilde{\ell}}(i_2)
p^{i_2}-\sigma^{\widetilde{\ell}}(i_1)
p^{i_1}-\frac{\widetilde{\ell}_{[i_2]}-\widetilde{\ell}_{[i_1]}}{\xi}.$$Define
$\overline{\ell}$ to be the unique element of $[0,q-2]$ satisfying the modulo
$\xi{\mathbb Z}$ congruence\begin{gather}\overline{\ell}\equiv
2\sigma^{\widetilde{\ell}}(i_1)p^{i_1}-2\sigma^{\widetilde{\ell}}(i_2)p^{i_2}-\sum_{i=i_2}^{i_1+f-1}m_ip^i+\sum_{i=i_1}^{i_2-1}m_ip^i.\label{kohlebrennt00}\end{gather}We
assume $\overline{\ell}\ne0$.\footnote{The case $\overline{\ell}=0$ would bring us back to Corollary \ref{cusanus00}.} For $i\in{\mathbb Z}$ put$$\sigma^{\overline{\ell},\widetilde{\ell}}(i)=\left\{\begin{array}{l@{\quad:\quad}l}
                                                                                                               1&\sigma^{\widetilde{\ell}}(i)+\sigma^{\overline{\ell}}(i)\ne4\\2&\sigma^{\widetilde{\ell}}(i)+\sigma^{\overline{\ell}}(i)=4\end{array}\right.$$$${\mathcal J}^{\overline{\ell},\widetilde{\ell}}_{\mathcal E}={\mathcal J}^{\overline{\ell}}_{\mathcal E}\cap{\mathcal J}^{\widetilde{\ell}}_{\mathcal E},\quad\quad\quad\quad{\mathcal J}^{\overline{\ell},\widetilde{\ell}}_{\mathcal D}={\mathcal J}-{\mathcal J}^{\overline{\ell},\widetilde{\ell}}_{\mathcal E}.$$
      
\begin{lem}\label{unionreuss} We have $\sigma^{\overline{\ell}}(i_2)-\sigma^{\overline{\ell},\widetilde{\ell}}(i_2)=1$ if and only if $\sigma^{\overline{\ell}}(i_2)=2$ and $\sigma^{\widetilde{\ell}}(i_2)=1$, and this happens if and only if one of the following three equivalent conditions hold:\begin{align}2\sigma^{\widetilde{\ell}}(i_1)-2\sigma^{\widetilde{\ell}}(i_2)-\sum_{i=i_2}^{i_1+f-1}m_ip^i+\sum_{i=i_1+f}^{i_2+f-1}m_ip^i&>\xi p^{i_2},\notag\\2\sigma^{\widetilde{\ell}}(i_1)-2\sigma^{\widetilde{\ell}}(i_2)-\sum_{i=i_2}^{i_1+f-1}m_ip^i+\sum_{i=i_1+f}^{i_2+f-1}m_ip^i&=\xi p^{i_2}+\overline{\ell}_{[i_2]},\notag\\2\sigma^{\widetilde{\ell}}(i_1)-2\sigma^{\widetilde{\ell}}(i_2)-\sum_{i=i_2}^{i_1+f-1}m_ip^i+\sum_{i=i_1+f}^{i_2+f-1}m_ip^i&\ne\overline{\ell}_{[i_2]}.\notag\end{align}
\end{lem}                                                                                                       
{\sc Proof:} This is easy.\hfill$\Box$\\

\begin{lem}\label{khw80brusch} (a)  We have\begin{align}\xi w-q(n_{{\bf x}}-1)&=\widetilde{\ell}_{[i_1]}-\sigma^{\widetilde{\ell}}(i_1){\xi}p^{i_1},\label{julgard000}\\e_{\sigma}(i)-q(n_{{\bf x}}-1)&=\widetilde{\ell}_{[i]}-\sigma{\xi}p^i\quad
\mbox{ for }i_1\le
i\le i_2\mbox{ and }\sigma\in\{1,2\},\label{qqjulsa1}\\q(n_{{\bf
    x}}-2w-1)&=\overline{\ell}_{[i_2]}-\sigma^{\overline{\ell},\widetilde{\ell}}(i_2){\xi}p^{i_2},\label{allerseel}\\q(n_{{\bf
    x}}-2w-1)+e_{\sigma^{\widetilde{\ell}}(i)}(i)&=\overline{\ell}_{[i]}-\sigma^{\overline{\ell},\widetilde{\ell}}(i){\xi}p^i\quad\mbox{
  for }i_1< i<i_2\label{qqjulsa3}.\end{align}

(b) The inequalities (\ref{voklau1}), (\ref{voklau2}), (\ref{voklau3}) are
  satisfied, as well as\begin{gather}n_{\bf x}\le\xi^2p^{i_1-f}.\label{verz}\end{gather}
  
(c) $F(t)$ belongs to $1+t^{qn_{{\bf x}}-{\xi}(w+1)}k[[t]]$.

(d) If $\tau=0$ then modulo $\nabla$ we have\begin{gather}\check{F}(t)\equiv\sum_{i\in{\mathcal J}^{\overline{\ell},\widetilde{\ell}}_{\mathcal D}}\alpha_{i}t^{e_1(i)}+\sum_{   i\in{\mathcal J}^{\overline{\ell},\widetilde{\ell}}_{\mathcal E}}\alpha_{i}t^{e_2(i)}H(p^{-i}\overline{\ell}_{[i]},i-2f),\label{attendedomine}\end{gather}and in the definition of $\overline{D}$ (formula (\ref{frankvirsik})) we may replace $F^{-1}(t)$ by\begin{gather}\left\{\begin{array}{l@{\quad:\quad}l}-\check{F}(t)+1&\sigma^{\overline{\ell},\widetilde{\ell}}(i_2)=1\\-\check{F}(t)+H(p^{-i_2}\overline{\ell}_{[i_2]},i_2-2f)&\sigma^{\overline{\ell},\widetilde{\ell}}(i_2)=2\end{array}\right.\label{allsouls}\end{gather}
without changing $\nabla$.

(e) $G(t)$ belongs to $1+t^{{\rm max}\{n_{{\bf x}},n_{{\bf y}}-w\}}k[[t]]$.

\end{lem}

{\sc Proof:} (a) Formulae (\ref{julgard000}) and (\ref{qqjulsa1}) are trivial. The
left hand side in formula (\ref{qqjulsa3}) is congruent modulo $\xi{\mathbb Z}$ with the right hand side. Since the right hand belongs to $-p^i[1,q-2]$ if $\sigma^{\overline{\ell},\widetilde{\ell}}(i)=1$, resp. to $-p^i[1,q-2]-p^i\xi$ if $\sigma^{\overline{\ell},\widetilde{\ell}}(i)=2$, it will be enough to show the same for the
left hand side. Thus, it will be enough to show\begin{gather}q(n_{{\bf x}}-2w-1)+e_{1}(i)\in\left\{\begin{array}{l@{\quad:\quad}l}-p^i[1,q-2]&\sigma^{\widetilde{\ell}}(i)=\sigma^{\overline{\ell},\widetilde{\ell}}(i)\\\xi p^i-p^i[1,q-2]&\sigma^{\widetilde{\ell}}(i)=\sigma^{\overline{\ell},\widetilde{\ell}}(i)+1\end{array}\right.\label{notag}\end{gather}We
compute$$q(n_{{\bf
    x}}-2w-1)+e_{1}(i)=2\sigma^{\widetilde{\ell}}(i_1)p^{i_1+f}-2\sigma^{\widetilde{\ell}}(i_2)p^{i_2}+\widetilde{\ell}_{[i]}-\xi
p^i+\frac{2}{\xi}\widetilde{\ell}_{[i_2]}-\frac{2q}{\xi}\widetilde{\ell}_{[i_1]}$$$$=2\sigma^{\widetilde{\ell}}(i_1)p^{i_1+f}-2\sigma^{\widetilde{\ell}}(i_2)p^{i_2}+\sum_{j=i}^
{i+f-1}(m_j+1-p)p^j+\sum_{j=i_2}^{i_1+f-1}(-2m_j)p^j.$$This is the same
as$$\sum_{j=i}^{i_2-1}(m_j+1-p)p^j+(1-2\sigma^{\widetilde{\ell}}(i_2)-m_{i_2})p^{i_2}+\sum_{j=i_2+1}^{i_1+f-1}(-m_j)p^j+(2\sigma^{\widetilde{\ell}}(i_1)-p+m_{i_1})p^{i_1+f}+\sum_{j=i_1+f+1}^{i+f-1}(m_j+1-p)p^j.$$The
claim follows.

In the same way formula (\ref{allerseel}) is proven.

(b) Formula (\ref{voklau1}) is equivalent
with $$\sigma^{\widetilde{\ell}}(i_2){\xi}p^{i_2}-\widetilde{\ell}_{[i_2]}\ge
2(\sigma^{\widetilde{\ell}}(i_1){\xi}p^{i_1}-\widetilde{\ell}_{[i_1]})+1,$$and formula
(\ref{voklau2}) is equivalent with$$\sigma^{\widetilde{\ell}}(i_1)\xi
p^{i_1+f}-q\widetilde{\ell}_{[i_1]}\ge 2(\sigma^{\widetilde{\ell}}(i_2)\xi
p^{i_2}-\widetilde{\ell}_{[i_2]})+q.$$Both these follow from Lemma \ref{hanny}. Formula (\ref{voklau3}) is trivial. Formula (\ref{verz}) is equivalent with$$q+\sigma^{\widetilde{\ell}}(i_2)\xi
p^{i_2}-\widetilde{\ell}_{[i_2]} \le\xi^2p^{i_1}$$which is clear.

(c) It is enough to see\begin{gather}p^{i_2-f}\xi^2\ge qn_{\bf x}-\xi(w+1)\quad\mbox{ if }\sigma^{\widetilde{\ell}}(i_2)=2,\label{sana2}\\e_{\sigma^{\widetilde{\ell}}(i)}(i)\ge qn_{\bf x}-\xi(w+1)\quad\mbox{ for all }i_1<i<i_2.\label{sana0}\end{gather}Formula (\ref{sana0}) is equivalent with$$(\sigma^{\widetilde{\ell}}(i_2)\xi p^{i_2}-\widetilde{\ell}_{[i_2]})-(\sigma^{\widetilde{\ell}}(i)\xi p^{i}-\widetilde{\ell}_{[i]})-(\sigma^{\widetilde{\ell}}(i_1)\xi p^{i_1}-\widetilde{\ell}_{[i_1]})\ge1.$$That this is satisfied
follows from Lemma \ref{hanny}. Formula (\ref{sana2}) is equivalent with$$p^{i_2+f}-2p^{i_2}+p^{i_2-f}\ge 1+\sigma^{\widetilde{\ell}}(i_1)\xi p^{i_1}-\widetilde{\ell}_{[i_1]}$$which is clearly okay.

(d) Formula (\ref{attendedomine}) follows from \begin{gather}\overline{\ell}_{[i]}-\widetilde{\ell}_{[i]}\in p^{i+1}{\mathbb
    Z}\quad\mbox{ for }i\in{\mathcal J}^{\overline{\ell},\widetilde{\ell}}_{\mathcal
  E},\label{schreck0}\\H(p^{-i}\widetilde{\ell}_{[i]},i-2f)\equiv 1\mbox{ modulo
}\nabla\quad\mbox{ for }i\in{\mathcal J}_{\mathcal E}^{\widetilde{\ell}}\cap
{\mathcal J}_{\mathcal D}^{\overline{\ell},\widetilde{\ell}}={\mathcal J}_{\mathcal E}^{\widetilde{\ell}}\cap
{\mathcal J}_{\mathcal D}^{\overline{\ell}}.\label{schreck1}\end{gather}The
proof of formula (\ref{schreck0}) is straightforward, the one of formula
(\ref{schreck1}) follows in a moment. Next, we observe that if $\sigma^{\overline{\ell},\widetilde{\ell}}(i_2)=2$ (and
hence $\sigma^{\widetilde{\ell}}(i_2)=2$), then we have \begin{gather}H(m_{i_2},i_2-2f)-1=1-H(p^{-i_2}\overline{\ell}_{[i_2]},i_2-2f)\label{schreckfast}\end{gather}as follows from\begin{gather}\overline{\ell}_{[i_2]}+\widetilde{\ell}_{[i_2]}+\sigma^{\widetilde{\ell}}(i_2)+\sigma^{\overline{\ell},\widetilde{\ell}}(i_2)\in
p{\mathbb Z},\label{schreckfa}\notag\end{gather}which indeed is a consequence of
$\sigma^{\overline{\ell},\widetilde{\ell}}(i_2)=2$.

Next, we claim\begin{gather}H(m_{i_2},i_2-2f)\equiv 1\mbox{ modulo
}\nabla\quad\mbox{ if
}\sigma^{\widetilde{\ell}}(i_2)-\sigma^{\overline{\ell},\widetilde{\ell}}(i_2)=1.\label{fastlo}\end{gather}The proof of
formulae (\ref{fastlo}) and (\ref{schreck1}) is the same, so we concentrate on
that for formulae (\ref{fastlo}). For this we need to see\begin{gather}p^{i_2-f}\xi^2+q(n_{\bf
  x}-2w-1)\ge n_{\bf x}-w\quad\mbox{ if
}\sigma^{\widetilde{\ell}}(i_2)-\sigma^{\overline{\ell},\widetilde{\ell}}(i_2)=1,\label{lent5}\end{gather}to
be proven in a moment.

To see formula (\ref{allsouls}) we need to see
$$F(t)^{-1}=(1+(F(t)-1))^{-1}\equiv 1-(F(t)-1)$$modulo
$\nabla$ (still assuming $\tau=0$). Given the claims (\ref{schreck0}),
(\ref{schreck1}), (\ref{schreckfast}), (\ref{fastlo}), this is a
matter of showing that the squares (and hence all mutual products with at least two factors) of all the
monomials (in $t$) with positive exponent appearing in $F(t)-1$ vanish modulo
$\nabla$. For the monomials appearing in $H(p^{-i_2}\overline{\ell}_{[i_2]},i_2-2f)-1$ and
$H(m_{i_2},i_2-2f)-1$ this follows from\begin{gather}2
p^{i_2-f}\xi^2+q(n_{\bf x}-2w-1)\ge n_{\bf x}-w\quad\mbox{ if }\sigma^{\widetilde{\ell},\overline{\ell}}(i_2)=2\label{sana5}.\end{gather}For
the other monomials this follows from\begin{gather}2 e_{\sigma^{\overline{\ell},\widetilde{\ell}}(i)}(i)+q(n_{\bf x}-2w-1)\ge n_{\bf
    x}-w\quad\mbox{ for all }i_1<i<i_2\label{sana4}.\end{gather}To
see formulae (\ref{sana5}) and (\ref{lent5}) we need to
see$$Ap^{i_2-f}\xi^2+q(\frac{2\xi
  p^{i_2}-\widetilde{\ell}_{[i_2]}}{q}-2(2p^{i_2}-\sigma^{\widetilde{\ell}}(i_1)p^{i_1}-\frac{\widetilde{\ell}_{[i_2]}-\widetilde{\ell}_{[i_1]}}{\xi}))$$$$\ge
1+\frac{2\xi
  p^{i_2}-\widetilde{\ell}_{[i_2]}}{q}-2p^{i_2}+\sigma^{\widetilde{\ell}}(i_1)p^{i_1}+\frac{\widetilde{\ell}_{[i_2]}-\widetilde{\ell}_{[i_1]}}{\xi}$$with
$A=1$ if $\sigma^{\widetilde{\ell}}(i_2)-\sigma^{\overline{\ell},\widetilde{\ell}}(i_2)=1$ and
$A=2$ if $\sigma^{\overline{\ell},\widetilde{\ell}}(i_2)=2$. This is equivalent
with $$(1+\frac{2q-1}{\xi
  q})\widetilde{\ell}_{[i_2]}-p^{i_2}\frac{6\xi}{q}+\frac{1-2q}{\xi}\widetilde{\ell}_{[i_1]}+p^{i_1}\sigma^{\widetilde{\ell}}(i_1)(2q-1)\ge1$$if
$\sigma^{\overline{\ell},\widetilde{\ell}}(i_2)=2$, and with$$(1+\frac{2q-1}{\xi
  q})\widetilde{\ell}_{[i_2]}-p^{i_2}(q+\frac{4\xi}{q})+\frac{1-2q}{\xi}\widetilde{\ell}_{[i_1]}+p^{i_1}\sigma^{\widetilde{\ell}}(i_1)(2q-1)\ge1$$if
$\sigma^{\widetilde{\ell}}(i_2)-\sigma^{\overline{\ell},\widetilde{\ell}}(i_2)=1$. The first of
these statements is clear, for the second one use Lemma \ref{unionreuss}.

Formula (\ref{sana4}) is equivalent with $$(\sigma^{\widetilde{\ell}}(i_2)\xi p^{i_2}-\widetilde{\ell}_{[i_2]})-2(\sigma^{\overline{\ell},\widetilde{\ell}}(i)\xi p^{i}-\widetilde{\ell}_{[i]})+\frac{2q-1}{\xi}(\sigma^{\widetilde{\ell}}(i_2)\xi p^{i_2}-\widetilde{\ell}_{[i_2]}-\sigma^{\widetilde{\ell}}(i_1)\xi p^{i_1}+\widetilde{\ell}_{[i_1]})\ge1.$$That this is satisfied
follows from Lemma \ref{hanny}.

(e) If $\sigma^{\widetilde{\ell}}(i_1)=1$ this is trivial. If $\sigma^{\widetilde{\ell}}(i_1)=2$ then it
follows from formula (\ref{verz}).\hfill$\Box$\\

Let $\widetilde{u}, \overline{u}\in[0,q-2]$ be such
  that \begin{gather}2(\widetilde{u}-\overline{u})\equiv\widetilde{\ell}-\overline{\ell}\quad\mbox{
      modulo }\xi{\mathbb
      Z},\label{1rook2}\\\widetilde{u}-\overline{u}-\widetilde{\ell}\in\xi{\mathbb
      Z}+p^{i_1}[0,\sum_{j=0}^{i_2-i_1}p^{j}].\label{2rook2}\end{gather}

{\bf Remark:} For each $\widetilde{u}\in[0,q-2]$ (resp. each
  $\overline{u}\in[0,q-2]$) there is a unique $\overline{u}\in[0,q-2]$ (resp. a unique
  $\widetilde{u}\in[0,q-2]$) such
  that formulae (\ref{1rook2}) and (\ref{2rook2}) both hold true. Indeed, formula (\ref{1rook2}) (and the definition of
  $\overline{\ell}$) means that either $\widetilde{u}-\overline{u}-\widetilde{\ell}
\equiv
\sigma^{\widetilde{\ell}}(i_2)p^{i_2}-\sigma^{\widetilde{\ell}}(i_1)p^{i_1}-\sum_{i=i_1}^{i_2-1}m_ip^i$
or $\widetilde{u}-\overline{u}-\widetilde{\ell}
\equiv \frac{1}{2}\xi+
\sigma^{\widetilde{\ell}}(i_2)p^{i_2}-\sigma^{\widetilde{\ell}}(i_1)p^{i_1}-\sum_{i=i_1}^{i_2-1}m_ip^i$ modulo $\xi{\mathbb
  Z}$. The first of these two alternatives is always solvable by $\widetilde{u}, \overline{u}\in[0,q-2]$ satisfying formula
    (\ref{2rook2}), but the second one is never solvable by $\widetilde{u}, \overline{u}\in[0,q-2]$ satisfying formula
    (\ref{2rook2}).

    It also follows that we could have replaced the joint conditions (\ref{1rook2}) and (\ref{2rook2}) by the single condition $\widetilde{u}-\overline{u}-\widetilde{\ell}
\equiv\sigma^{\widetilde{\ell}}(i_2)p^{i_2}-\sigma^{\widetilde{\ell}}(i_1)p^{i_1}-\sum_{i=i_1}^{i_2-1}m_ip^i$ modulo $\xi{\mathbb
  Z}$. The point of sticking with the joint conditions (\ref{1rook2}) and (\ref{2rook2}) instead is that they are not formulated in terms of the digits of $\widetilde{\ell}$.   

\begin{lem}\label{focault} We have the modulo $\xi{\mathbb Z}$-congruences\begin{gather}\overline{\ell}\equiv-\widetilde{\ell}-2w,\label{rook0}\\n_{{\bf x}}-1\equiv-\widetilde{\ell},\label{rook1}\\\widetilde{u}- \overline{u}\equiv \widetilde{\ell}+w.\label{rook2}\end{gather}
  \end{lem} 

  {\sc Proof:} The congruence (\ref{rook1}) is clear, the congruence (\ref{rook0}) follows from formula (\ref{allerseel}). To see the congruence (\ref{rook2}) begin with the congruences\begin{gather}2(\widetilde{u}- \overline{u})\equiv\widetilde{\ell}-\overline{\ell}\equiv 2(\sigma^{\widetilde{\ell}}(i_2)p^{i_2}-\sigma^{\widetilde{\ell}}(i_1)p^{i_1}+\sum_{i=i_2}^{i_1+f-1}m_ip^i),\notag\\w=\sigma^{\widetilde{\ell}}(i_2)p^{i_2}-\sigma^{\widetilde{\ell}}(i_1)p^{i_1}-\frac{\widetilde{\ell}_{[i_2]}-\widetilde{\ell}_{[i_1]}}{\xi}= \sigma^{\widetilde{\ell}}(i_2)p^{i_2}-\sigma^{\widetilde{\ell}}(i_1)p^{i_1}-\sum_{i=i_1}^{i_2-1}m_ip^i\label{sawachdi}\end{gather}which show $2(\widetilde{u}- \overline{u}-\widetilde{\ell})\equiv 2w$ modulo $\xi{\mathbb Z}$. But formula (\ref{sawachdi}) also shows $w\in p^{i_1}[0,\sum_{j=0}^{i_2-i_1}p^{j}]$, which together with formula (\ref{2rook2}) then gives the claim.\hfill$\Box$\\

Let $\Gamma$ act on $M$ by means
of$$\gamma*{{\bf y}}=a^{\widetilde{u}+1}_{\gamma,1}\frac{t}{[\gamma]}{{\bf y}},\quad\quad\gamma*{{\bf x}}=a^{\widetilde{u}-\widetilde{\ell}+1}_{\gamma,1}(\frac{t}{[\gamma]})^{n_{{\bf x}}}{{\bf x}}$$for $\gamma\in\Gamma$.

\begin{lem} (a) Writing
  $\overline{x}=t^{n_{{\bf x}}-w-1}{{\bf x}}$ and $\overline{y}=t^{w}{{\bf y}}$ and $\widetilde{x}=t^{n_{{\bf x}}-1}{{\bf x}}$ and $\widetilde{y}={{\bf y}}$, the induced $\Gamma$-action on
$M/\langle t^{n_{{\bf y}}}{{\bf y}}, t^{n_{{\bf x}}}{{\bf x}}\rangle$ satisfies\begin{align}\gamma*\overline{y}&=a_{\gamma,1}^{\overline{u}-\overline{\ell}+1}\frac{t}{[\gamma]}\overline{y},\label{oktsa1}\\\gamma*\overline{x}&=a_{\gamma,1}^{\overline{u}+1}\frac{t}{[\gamma]}\overline{x},\label{oktsa2}\\\gamma*\widetilde{x}&=a^{\widetilde{u}-\widetilde{\ell}+1}_{\gamma,1}\frac{t}{[\gamma]}\widetilde{x},\label{oktsa3}\\\gamma*\widetilde{y}&=a^{\widetilde{u}+1}_{\gamma,1}\frac{t}{[\gamma]}\widetilde{y}\label{oktsa4}.\end{align} (b) The $\Gamma$-action on $M$ preserves $\nabla$.

\end{lem}

{\sc Proof:} (a) This is due to formulae (\ref{rook0}), (\ref{rook1}), (\ref{rook2}) as well as the following: From ${\rm ord}_p(w)\ge i_1$ and formula
(\ref{verz}) we get $w+\xi p^{{\rm ord}_p(w)}\ge n_{{\bf y}}$. As
$[\gamma]\in a_{\gamma,1}t+t^qk[[t]]$ this implies
$(\frac{[\gamma]}{a_{\gamma,1}t})^w\overline{y}\equiv \overline{y}$ modulo
$\langle t^{n_{{\bf y}}}{{\bf y}}\rangle$, and hence formula
(\ref{oktsa1}). Similarly, from ${\rm ord}_p(w)\ge i_1$ and formula
(\ref{verz}) we get $\xi p^{{\rm ord}_p(w)}\ge w+1$. As $[\gamma]\in a_{\gamma,1}t+t^qk[[t]]$ this implies
$(\frac{[\gamma]}{a_{\gamma,1}t})^w\overline{x}\equiv \overline{x}$ modulo
$\langle t^{n_{{\bf x}}}{{\bf x}}\rangle$, and hence formula (\ref{oktsa2}). Formula (\ref{oktsa3}) is clear, formula (\ref{oktsa4}) is trivial.

(b) We use Proposition \ref{vacet10}, separately for the cases $\tau\ne0$ and $\tau=0$. 

If $\tau=0$ we translate from the present setting to the setting in Proposition \ref{vacet10} (with $x$, $n_x$, $y$, $n_y$, ${\mathcal D}$, ${\mathcal E}$,
$\ell$, $u$ as in Proposition \ref{vacet10}) by means of the
assignments$$y=\overline{y},\quad\quad n_y=n_{{\bf y}}-w,\quad\quad
x=\overline{x},\quad\quad n_x=w+1,\quad\quad \ell=\overline{\ell},\quad\quad
u=\overline{u},$$\begin{gather}{\mathcal E}=-2f+\left\{\begin{array}{l@{\quad:\quad}l}{\mathcal J}^{\overline{\ell},\widetilde{\ell}}_{\mathcal E}&\sigma^{\overline{\ell},\widetilde{\ell}}(i_2)=1\\{\mathcal J}^{\overline{\ell},\widetilde{\ell}}_{\mathcal E}\cup\{i_2\}& \sigma^{\overline{\ell},\widetilde{\ell}}(i_2)=2\end{array}\right.,\notag\\{\mathcal D}=\left\{\begin{array}{l@{\quad:\quad}l}{\mathcal J}^{\overline{\ell},\widetilde{\ell}}_{\mathcal D}\cup\{i_2\}&\sigma^{\overline{\ell},\widetilde{\ell}}(i_2)=1\\{\mathcal J}^{\overline{\ell},\widetilde{\ell}}_{\mathcal D}&\sigma^{\overline{\ell},\widetilde{\ell}}(i_2)=2\end{array}\right..\label{nokalender}\end{gather}By
formula (\ref{qqjulsa3}), condition (\ref{sommerdornwald5}) is equivalent
with$$q(n_{{\bf x}}-2w-1)+e_1(i)+p^i\xi\ge n_{{\bf y}}-w\quad\mbox{ for all
}i\in {\mathcal D}.$$This is equivalent with$$\frac{2q-1}{\xi}(\sigma^{\widetilde{\ell}}(i_1)\xi
p^{i_1}-\widetilde{\ell}_{[i_1]})-\frac{2q-1}{\xi q}(\sigma^{\widetilde{\ell}}(i_2)\xi
p^{i_2}-\widetilde{\ell}_{[i_2]})+\widetilde{\ell}_{[i]}\ge1.$$That this is satisfied
follows from Lemma \ref{hanny}.

To prove condition (\ref{fuhleresel200}) for $i-2f\in{\mathcal E}$ it is
enough, by Lemma \ref{letzterst}, to prove$$p^{i-f}(q-\sum_{s=0}^{f-1}p^s)\ge
n_{\bf y}-w.$$This is equivalent
with$$p^{i}(q-\sum_{s=0}^{f-1}p^s)\ge\frac{q}{\xi}(\sigma^{\widetilde{\ell}}(i_1)\xi
p^{i_1}-\widetilde{\ell}_{[i_1]})-\frac{1}{\xi}(\sigma^{\widetilde{\ell}}(i_2)\xi
p^{i_2}-\widetilde{\ell}_{[i_2]})+q.$$As $i>i_1$ this is clear. Condition (\ref{maruf200}) reads $$q(w+1)+
\overline{\ell}_{[i]}-\xi(2p^{i}+1)\ge n_{\bf y}-w\quad\mbox{ for
}i-2f\in{\mathcal E}$$which using formula (\ref{qqjulsa3}) is seen to be
equivalent with$$\frac{\xi}{q}(\sigma^{\widetilde{\ell}}(i_2)\xi
p^{i_2}-\widetilde{\ell}_{[i_2]})-(2\xi
p^{i}-\widetilde{\ell}_{[i]})+(\sigma^{\widetilde{\ell}}(i_1)\xi
p^{i_1}-\widetilde{\ell}_{[i_1]})\ge 0.$$That this is satisfied
follows from Lemma \ref{hanny}. 

Inserting formulae (\ref{allerseel}), (\ref{qqjulsa3}), (\ref{allsouls}) into the defining formula (\ref{frankvirsik}) shows that the setting of Proposition \ref{vacet10} indeed applies here.

If $\tau\ne 0$ we translate from the present setting to the setting in Proposition \ref{vacet10} (with $x$, $n_x$, $y$, $n_y$, ${\mathcal D}$, ${\mathcal E}$,
$\ell$, $u$ as in Proposition \ref{vacet10}) by means of the
assignments$$x=\widetilde{y},\quad\quad n_x=n_{{\bf y}},\quad\quad
y=\widetilde{x},\quad\quad n_y=1,\quad\quad\ell=
\widetilde{\ell},\quad\quad\ u= \widetilde{u},$$$${\mathcal E}=-2f+\left\{\begin{array}{l@{\quad:\quad}l}
{\mathcal J}^{\widetilde{\ell}}_{\mathcal E} &(\sigma^{\widetilde{\ell}}(i_1),\sigma^{\widetilde{\ell}}(i_2))=(1,1)\\{\mathcal J}^{\widetilde{\ell}}_{\mathcal E}\cup\{i_2\}&(\sigma^{\widetilde{\ell}}(i_1),\sigma^{\widetilde{\ell}}(i_2))=(1,2)\\{\mathcal J}^{\widetilde{\ell}}_{\mathcal E}\cup\{i_1\}&(\sigma^{\widetilde{\ell}}(i_1),\sigma^{\widetilde{\ell}}(i_2))=(2,1)\\{\mathcal J}^{\widetilde{\ell}}_{\mathcal E}\cup\{i_1,i_2\}&(\sigma^{\widetilde{\ell}}(i_1),\sigma^{\widetilde{\ell}}(i_2))=(2,2)\end{array}\right.,$$$${\mathcal D}=\left\{\begin{array}{l@{\quad:\quad}l}
{\mathcal J}^{\widetilde{\ell}}_{\mathcal D}\cup\{i_1,i_2\} &(\sigma^{\widetilde{\ell}}(i_1),\sigma^{\widetilde{\ell}}(i_2))=(1,1)\\{\mathcal
  J}^{\widetilde{\ell}}_{\mathcal D}\cup\{i_1\}&(\sigma^{\widetilde{\ell}}(i_1),\sigma^{\widetilde{\ell}}(i_2))=(1,2)\\{\mathcal
  J}^{\widetilde{\ell}}_{\mathcal D}\cup\{i_2\}&(\sigma^{\widetilde{\ell}}(i_1),\sigma^{\widetilde{\ell}}(i_2))=(2,1)\\{\mathcal
  J}^{\widetilde{\ell}}_{\mathcal D}&(\sigma^{\widetilde{\ell}}(i_1),\sigma^{\widetilde{\ell}}(i_2))=(2,2)\end{array}\right..$$By formula (\ref{qqjulsa1})
resp. (\ref{qqjulsa3}), condition (\ref{sommerdornwald5}) is equivalent with$$q-qn_{{\bf
    x}}+e_1(i)+p^i\xi\ge 1\quad \mbox{ for all }i\in {\mathcal D}.$$This is equivalent with
$\xi+\widetilde{\ell}_{[i]}\ge 0$, which is okay. By formula (\ref{qqjulsa1})
resp. (\ref{qqjulsa3}) again, condition (\ref{marufmin00}) is equivalent with $$q(n_{{\bf y}}-1)\ge
qn_{{\bf x}}-e_1(i)-q\quad\mbox{ for all }i\in {\mathcal D}.$$This is a
consequence of formula (\ref{voklau3}). Formula (\ref{fuhleresel200}) holds true for all $i-2f\in{\mathcal E}$ by Lemma \ref{letzterst} as
$n_y=1$ and ${\mathcal E}\subset \Pi^{-1}({\mathcal
  E}(\widetilde{\ell}))$. Formula (\ref{maruf200}) for $i-2f\in{\mathcal E}$ reads
$qn_{\bf y}+\widetilde{\ell}_{[i]}-\xi(2p^{i}+1)\ge1$. This is equivalent
with$$\sigma^{\widetilde{\ell}}(i_2)\xi p^{i_2}-\widetilde{\ell}_{[i_2]}\ge 2\xi
p^{i}-\widetilde{\ell}_{[i]}.$$If $i=i_2$ then $\sigma^{\widetilde{\ell}}(i_2)=2$ so this is fine. If
$i<i_2$ then it follows from Lemma \ref{hanny}.

Inserting formulae (\ref{julgard000}), (\ref{qqjulsa1}) into the defining formula (\ref{vierwoend}) shows that the setting of Proposition \ref{vacet10} indeed applies here.
\hfill$\Box$\\

\subsection{Example}

\label{keydefo1}

This is the $\widetilde{\ell}=0$ version of subsection \ref{keydefo}. Fix integers $f\le i_1<
i_2< i_1+f$ and a set $\widetilde{\mathcal I}$ with $\Pi(\{i_1,i_2\})\subset
\widetilde{\mathcal I}\subset \Pi([i_1,i_2])$. Define$${\mathcal J}=[i_1+1,i_2-1]\cap \Pi^{-1}(\widetilde{\mathcal I}).$$For $i_1\le
i\le i_2$ put $e(i)={\xi}p^{i_2}-\xi
p^{i}$. Let $\{\alpha_{i}\}_{i\in{\mathcal J}}$ be scalars in $k$ and define$${F}(t)=1+\sum_{i\in{\mathcal J}}\alpha_{i}t^{e(i)},\quad\quad\quad\quad G(t)=1.$$Define the integers$$n_{\bf y}=n_{\bf
  x}=1+\frac{\xi p^{i_2}}{q},\quad \quad w=p^{i_2}-p^{i_1}.$$Define
$\overline{\ell}$ to be the unique element of $[0,q-2]$ satisfying the modulo
$\xi{\mathbb Z}$-congruence\begin{gather}\overline{\ell}\equiv
2p^{i_1}-2p^{i_2}.\label{kohlebrennt000}\end{gather}

\begin{lem}\label{khw80brusch000} (a)  We have\begin{align}\xi w-q(n_{{\bf x}}-1)&=-{\xi}p^{i_1},\notag\\e(i)-q(n_{{\bf x}}-1)&=-{\xi}p^i\quad
\mbox{ for }i_1\le
i\le i_2,\notag\\q(n_{{\bf
    x}}-2w-1)+e(i)&=-{\xi}p^i\quad\mbox{
  for }i_1< i\le i_2\notag.\end{align}

(b) The inequalities (\ref{voklau1}), (\ref{voklau2}), (\ref{voklau3}) are
  satisfied, as well as\begin{gather}n_{\bf x}\le\xi^2p^{i_1-f}.\notag\end{gather}
  
(c) $F(t)$ belongs to $1+t^{qn_{{\bf x}}-{\xi}(w+1)}k[[t]]$.

(d) If $\tau=0$ then we may replace $F^{-1}(t)$ by $1-\sum_{i\in{\mathcal J}}\alpha_{i}t^{e(i)}$ without changing $\nabla$.

\end{lem}

{\sc Proof:} The same as for Lemma \ref{khw80brusch}, only much easier.\hfill$\Box$\\

Again we demand that $\widetilde{u}, \overline{u}\in[0,q-2]$ satisfy formulae (\ref{1rook2}) and (\ref{2rook2}).

\begin{lem} We have the modulo $\xi{\mathbb Z}$ congruences$$\overline{\ell}\equiv-2w,\quad \quad n_{{\bf x}}\equiv 1,\quad \quad\widetilde{u}-\overline{u}\equiv w.$$
  \end{lem} 

{\sc Proof:} The same as for Lemma \ref{focault}, only much easier.\hfill$\Box$\\

Assume $b\tau^2\ne a$ and put$${\bf d}_{i_1}=\frac{-\tau a^{-1}}{b\tau^2-a},\quad\quad {\bf d}_{i_2}=\frac{-\tau^2 a^{-1}}{b\tau^2-a},\quad\quad {\bf d}_{i}=\frac{-\tau^2 a^{-1}\alpha_i}{b\tau^2-a}\quad\mbox{ for }i\in{\mathcal J}.$$Let $\Gamma$ act on $M$ by means
of$$\gamma*{{\bf y}}=a^{\widetilde{u}}_{\gamma,1}(\frac{a^{\widetilde{u}}_{\gamma,1} t}{[\gamma]}{{\bf y}}+ab\sum_{i\in{\mathcal J}\cup\{i_1, i_2\}}{\bf d}_ia_{\gamma,1}^{p^i}t^{n_{\bf x}-1}{\bf x}),$$$$\gamma*{{\bf x}}=a^{\widetilde{u}-\widetilde{\ell}+1}_{\gamma,1}(\frac{t}{[\gamma]})^{n_{{\bf x}}}{{\bf x}}$$for
$\gamma\in\Gamma$.

\begin{lem}\label{focault0} (a) Writing
  $\overline{x}=t^{n_{{\bf x}}-w-1}{{\bf x}}$ and $\overline{y}=t^{w}{{\bf y}}$ and $\widetilde{x}=t^{n_{{\bf x}}-1}{{\bf x}}$ and $\widetilde{y}={{\bf y}}$, the induced $\Gamma$-action on
$M/\langle t^{n_{{\bf y}}}{{\bf y}}, t^{n_{{\bf x}}}{{\bf x}}\rangle$ satisfies\begin{align}\gamma*\overline{y}&=a_{\gamma,1}^{\overline{u}-\overline{\ell}+1}\frac{t}{[\gamma]}\overline{y}\quad\quad\mbox{ if }\tau=0,\label{oktsa1000}\\\gamma*\overline{x}&=a_{\gamma,1}^{\overline{u}+1}\frac{t}{[\gamma]}\overline{x},\notag\\\gamma*\widetilde{x}&=a^{\widetilde{u}+1}_{\gamma,1}\frac{t}{[\gamma]}\widetilde{x},\notag\\\gamma*\widetilde{y}&=a^{\widetilde{u}}_{\gamma,1}(\frac{a^{\widetilde{u}}_{\gamma,1} t}{[\gamma]}{\widetilde{y}}+ab\sum_{i\in{\mathcal J}\cup\{i_1, i_2\}}{\bf d}_ia_{\gamma,1}^{p^i}\widetilde{x})\notag.\end{align} (b) The $\Gamma$-action on $M$ preserves $\nabla$.

\end{lem}

{\sc Proof:} (a) The same as in Lemma \ref{focault0}, except that for formula (\ref{oktsa1000}) we in addition observe ${\bf d}_i=0$ if $\tau=0$, for all $i\in{\mathcal J}\cup\{i_1, i_2\}$.

(b) We use Proposition \ref{vacet10}, separately for the cases $\tau\ne0$ and $\tau=0$. 

If $\tau=0$ the argument is exactly as in Lemma \ref{focault0}. 

If $\tau\ne 0$ we translate from the present setting to the setting in Proposition \ref{vacet10} (with $x$, $n_x$, $y$, ${\mathcal D}$, $u$ as in Proposition \ref{discord0}) by means of the
assignments$$x=\widetilde{y},\quad\quad n_x=n_{{\bf y}},\quad\quad
y=\widetilde{x},\quad\quad u= \widetilde{u},$$$${\mathcal D}={\mathcal J}\cup\{i_1,i_2\}.$$Similarly as before we observe that formula (\ref{freu}) is satisfied, and a small computation shows that the scalars ${\bf d}_i$ are precisely so designed as to make Proposition \ref{discord0} applicable.\hfill$\Box$\\

\subsection{Reformulation}

\label{lichtherz}

\begin{lem}\label{hanny11} Let $\ell\in[0,q-2]$. Let $0\le m_i< p$ and $1\le a_i\le p$ for $0\le i<f$ be such that modulo $\xi{\mathbb Z}$ we have $$\ell\equiv\sum_{i=0}^{f-1}m_ip^i\equiv-\sum_{i=0}^{f-1}a_ip^i$$and assume $a_i\ne p$ for at least one $i$. Write $a_i=a_{i+f}$ and $m_i=m_{i+f}$  for all $i\in{\mathbb Z}$.

(a) For each $0\le i<f$ we have\begin{gather}a_{i}=\sigma^{{\ell}}(i+1)p-\sigma^{{\ell}}(i)-m_{i}.\label{waltergub}\end{gather}

(b) We have\begin{gather}{\mathcal E}({\ell})=\{0\le i<f\,|\,a_{j}=p\mbox{ and }a_{s}=p-1\mbox{ for all
  }j<s<i,\mbox{ for some }j<i\}.\label{lagathanks2}\end{gather}

(c) For $i\in{\mathbb Z}$ we have\footnote{The map $i_{\mathcal
    D}^{{\ell}}$ was defined in formula (\ref{stetigefreude}); for $\ell=0$ we take $i_{\mathcal
    D}^{{\ell}}=\Pi$, in particular $i_{\mathcal
    D}^{0}(x)=x$ for $x\in{{{\mathcal F}}}$.}\begin{gather}i_{\mathcal D}^{{\ell}}(i)=\Pi({\rm min}(\{j<i\,|\,a_j=p\mbox{ and }a_s=p-1\mbox{ for all }j<s<i \}\cup\{i\})).\label{lagathanks1}\end{gather}

(d) For $i_1<i_2<i_1+f$ we have the modulo $\xi{\mathbb Z}$-congruence \begin{gather}2\sigma^{{\ell}}(i_1)p^{i_1}-2\sigma^{{\ell}}(i_2)p^{i_2}-\sum_{i=i_2}^{i_1+f-1}m_ip^i+\sum_{i=i_1}^{i_2-1}m_ip^i\equiv\sum_{i=i_2}^{i_1+f-1}a_ip^i-\sum_{i=i_1}^{i_2-1}a_ip^i.\label{lichtmess1}\end{gather}

(e) For each $i\ge0$ we have\begin{gather}\ell_{[i]}-\sigma^{{\ell}}(i)p^i\xi=-\sum_{j=i}^{i+f-1}a_{j}p^j.\label{ashwed}\end{gather} 

\end{lem}

{\sc Proof:} (a) Formula (\ref{waltergub}) is seen to be equivalent with\begin{gather}a_i=\left\{\begin{array}{l@{\quad:\quad}l}p-m_i-1 &i\in{\mathcal D}({\ell})\mbox{ and }m_i<p-1\\p-m_i-2 &i\in{\mathcal E}({\ell})\mbox{ and }m_i<p-2\\2p-m_i-1 &i\in{\mathcal D}({\ell})\mbox{ and }m_i=p-1\\2p-m_i-2 &i\in{\mathcal E}({\ell})\mbox{ and }m_i\ge p-2\end{array}\right..\label{allein}\end{gather}For $0\le i<f$ denote by $\gamma(i)\in\{1,2,3,4\}$ the case number (in the given order from top to bottom) which applies in the stated formula for $a_i$. We find $$(\gamma(\Pi(i)),\gamma(\Pi(i+1)))\in\{(1,1), (1,3), (3,4), (4,4), (4,2), (2,1)\}$$for each $0\le i<f$. This implies statement (a).

Statements (b) and (c) can be derived from formula (\ref{allein}).
  
To prove statement (d) we fix $i_2$ and induct on $i_1$. We start with $i_1=i_2-f$ in which case $\sum_{i=i_2}^{i_1+f-1}$ is the empty sum and the statement is trivial. For the induction step we write the left hand side of formula (\ref{lichtmess1}) modulo $\xi{\mathbb Z}$ as$$2\sigma^{{\ell}}(i_1)p^{i_1}-2\sigma^{{\ell}}(i_1-1)p^{i_1-1}-2m_{i_1-1}p^{i_1-1}+[2\sigma^{{\ell}}(i_1-1)p^{i_1-1}-2\sigma^{{\ell}}(i_2)p^{i_2}-\sum_{i=i_2}^{i_1+f-2}m_ip^i+\sum_{i=i_1-1}^{i_2-1}m_ip^i]$$and the right hand side as$$2a_{i_1-1}p^{i_1-1}+[\sum_{i=i_2}^{i_1+f-2}a_ip^i-\sum_{i=i_1-1}^{i_2-1}a_ip^i].$$By induction hypothesis, the respective summands in squared brackets are congruent modulo $\xi{\mathbb Z}$. Therefore we may conclude with formula (\ref{waltergub}), i.e. with$$a_{i_1-1}= \sigma^{{\ell}}(i_1)p-\sigma^{{\ell}}(i_1-1)-m_{i_1-1}.$$Statement (e) is easy to see.\hfill$\Box$\\ 

{\bf Remark:} For each $0\le i<f$ we also have$$m_i=\left\{\begin{array}{l@{\quad:\quad}l}p-a_i-1 &i\in{\mathcal D}({\ell})\mbox{ and }a_i<p\notag\\p-a_i-2 &i\in{\mathcal E}({\ell})\mbox{ and }a_i< p-1\notag\\2p-a_i-1 &i\in{\mathcal D}({\ell})\mbox{ and }a_i=p\notag\\2p-a_i-2 &i\in{\mathcal E}({\ell})\mbox{ and }a_i\ge p-1\notag\end{array}\right.$$This can be deduced from Lemma \ref{hanny11} (or can be proven similarly as Lemma \ref{hanny11}).\\

Let $(\widetilde{\ell}, \widetilde{u},\widetilde{\mathcal I})$ and $(\overline{\ell}, \overline{u},\overline{\mathcal I})$ be elements in ${\mathbb E}$, subject to all the following requirements.

We require $\overline{\ell}\ne0$ and that $\widetilde{u}, \overline{u}$ satisfy formulae (\ref{1rook2}) and (\ref{2rook2}). Define $1\le a_i\le p$ by $\widetilde{\ell}\equiv-\sum_{i=0}^{f-1}a_ip^i$ modulo $\xi{\mathbb Z}$, assuming $a_i\ne p$ for at least one $i$. As usual write $a_{i+f}=a_i$ for all $i$. We require that there are $2f\le i_1<i_2<i_1+f$ satisfying\begin{gather}\overline{\ell}\equiv\sum_{i=i_2}^{i_1+f-1}a_ip^i-\sum_{i=i_1}^{i_2-1}a_ip^i\quad\mbox{ modulo }\xi{\mathbb Z}\label{ftklaus}\end{gather}and\begin{gather}\{\Pi(i_1),\Pi(i_2)\}\subset\widetilde{\mathcal
    I}\subset\Pi([i_1,i_2]).\label{armand}\end{gather}We define
the map\begin{gather}\nu_{\widetilde{\ell}}^{\overline{\ell}}:{{{\mathcal F}}}\longrightarrow{{{\mathcal F}}},\label{lagathanks0}\\\nu_{\widetilde{\ell}}^{\overline{\ell}}(x)=\left\{\begin{array}{l@{\quad:\quad}l}i_{\mathcal
  D}^{\overline{\ell}}(x)&x\in{\mathcal
                           D}(\widetilde{\ell})\notag\\x&x\in{\mathcal  E}(\widetilde{\ell})\notag\end{array}\right.\notag\end{gather}and require\begin{gather}\overline{\mathcal I}=\nu_{\widetilde{\ell}}^{\overline{\ell}}(\widetilde{\mathcal I}-\{\Pi(i_1)\}).\label{armand1}\end{gather}
                     
\begin{kor}\label{cusanus} Let $(\widetilde{\ell}, \widetilde{u},\widetilde{\mathcal I})$ and $(\overline{\ell}, \overline{u},\overline{\mathcal I})$ be as above. Let $\overline{\bf D}\in{\mathcal
    V}(\overline{\ell},\overline{u},\overline{\mathcal I})$. There exists a $(\varphi,\Gamma)$-module ${\bf D}$ over $B$, where $B=k[\tau]$ if $\widetilde{\ell}\ne0$ and $B=k[\tau][\frac{1}{\tau^2-\epsilon}]$ for some $\epsilon\in k^{\times}$ if $\widetilde{\ell}=0$, with the following property: We have ${\bf D}\otimes_{B}k\cong \overline{\bf D}$ via $\tau\mapsto 0\in k$, and ${\bf D}\otimes_{B}k\in {\mathcal V}(\widetilde{\ell},\widetilde{u},\widetilde{\mathcal
  I})$ via $\tau\mapsto a$ for each $a\in k^{\times}$ (resp. each $a\in k^{\times}$ with $a^2\ne\epsilon$ if $\widetilde{\ell}=0$).
  \end{kor}

{\sc Proof:} Define $0\le m_i<p$ by $\widetilde{\ell}=\sum_{i=0}^{f-1}m_ip^i$; write $m_{i+f}=m_i$ for all $i$. From formula (\ref{lichtmess1}) we obtain the modulo
$\xi{\mathbb Z}$-congruence of formula (\ref{kohlebrennt00}). Corollary \ref{cusanus} is thus the upshot of subsections
\ref{keydefo} and \ref{keydefo1}. Namely, choose $\overline{C}$, $\overline{D}$ as at the beginning of subsection \ref{osterfreitag}, and $n_{\bf x}$, $n_{\bf y}$ as in subsection \ref{keydefo} resp. \ref{keydefo1}, such that we have an isomorphism of $(\varphi,\Gamma)$-modules$$\overline{\bf D}\cong(k[[t]][\varphi](k\overline{x}\oplus k\overline{y})/\langle
\overline{C},\overline{D},t^{w+1}\overline{x},t^{n_{\bf y}-w}\overline{y}\rangle)^*\otimes_{k[[t]]}k((t)).$$Use the coefficients of $\overline{C}$, $\overline{D}$ to define $C$, $D$ as at the beginning of subsection \ref{osterfreitag}. Then, by Lemma \ref{pitpaul} and the discussion of subsection \ref{keydefo} resp. \ref{keydefo1},$${\bf D}=({\rm Hom}_{k[\tau]}(\Delta,k[\tau]))\otimes_{k[\tau][[t]]}k[\tau]((t))$$with$$\Delta=k[\tau][[t]][\varphi](k{\bf x}\oplus k{\bf y})/\langle{C},{D},t^{n_{{\bf x}}}{{\bf x}},t^{n_{{\bf y}}}{{\bf y}}\rangle$$carries the structure of a $(\varphi,\Gamma)$-module over $k[\tau]$ as desired.

Note that the shifts performed by $i_{\mathcal D}^{\overline{\ell}}$, applied to the set given by formula (\ref{nokalender}), account for the shifts in (the proof of) Proposition \ref{comme}.\hfill$\Box$\\

\section{Number theoretic combinatorics}

\label{combi}

For $x\in{\mathbb Z}$ we often denote $\Pi(x)$ simply again by $x$. For $x,y\in{\mathbb Z}$ with $x-f\le y\le x$ we write
$[x,y]=\Pi([x-f,y])$ and $[y,x]=\Pi([y,x])$. For $x,y\in{\mathbb Z}$ with $x-f>y$ we
understand $[x,y]=\emptyset$. For $J\subset {{{\mathcal F}}}$ put$$J^c={{{\mathcal F}}}-J,\quad\quad\quad\quad J^-=\{j\in J\,|\,j-1\notin J\},\quad\quad\quad\quad
J^+=(J^c)^-,$$$$J^{c,1}=(J^c\cup J^-)-J^+=J^c+1=\{j+1\,|\,j\in J^c\}$$(where, as a first instance of our above convention, in the definition of $J^-$ (resp. of $J^{c,1}$) by $j-1$ (resp. $j+1$) we really mean $\Pi(j-1)$ (resp. $\Pi(j+1)$)). We then have a unique bijective map$$J^-\longrightarrow J^+,\quad j\mapsto j_J^+$$such that$$J=\coprod_{j\in J^-}[j,j_J^+-1].$$

{\bf Definition:} (a) For $\widetilde{\ell},\overline{\ell}$ in $[0,q-2]$ and a subset $J$ of ${{{\mathcal F}}}$ we declare the relation $\widetilde{\ell}\succ_J\overline{\ell}$ to hold if and only if we have modulo $\xi{\mathbb Z}$-congruences$$\widetilde{\ell}\equiv-\sum_{j\in{{{\mathcal F}}}}a_jp^j\quad\quad\mbox{ and }\quad\quad\overline{\ell}\equiv\sum_{j\in J}a_jp^j-\sum_{j\in J^c}a_jp^j$$with $1\le a_j\le p$ for all $j$, not all of them $=p$.

(b) On the set$${\mathbb D}=\{(\ell,{\mathcal I})\,\,|\,\,\ell\in[0,q-2],\, {\mathcal I}\subset{{{\mathcal F}}}\}$$ define relations $\succ_J$ for $J\subset{{{\mathcal F}}}$ by declaring $(\widetilde{\ell},\widetilde{\mathcal I})\succ_J(\overline{\ell},\overline{\mathcal I})$ if and only if\footnote{In this section, the map $\nu^{\overline{\ell}}_{\widetilde{\ell}}$ (formula (\ref{lagathanks0})) should be read in light of formulae (\ref{lagathanks2}) and (\ref{lagathanks1}).}

$$\widetilde{\ell}\succ_J\overline{\ell},\quad\quad\quad\quad J^-\cup J^+\subset \widetilde{\mathcal I},\quad\quad\quad\quad\overline{\mathcal I}\subset\nu^{\overline{\ell}}_{\widetilde{\ell}}(\widetilde{\mathcal I}\cap J^{c,1}).$$

{\bf Remarks:} (a) The definitions imply $|\widetilde{\mathcal I}|\ge |\overline{\mathcal I}|+|J^-|$.

(b) In Corollary \ref{cusanus}, we have $(\widetilde{\ell},\widetilde{\mathcal I})\succ_J(\overline{\ell},\overline{\mathcal I})$ with $J=[i_2, i_1-1]=\Pi([i_2, i_1+f-1])$.

\begin{lem}\label{mario1} Let $(\widetilde{\ell},\widetilde{\mathcal
    I})\succ_J(\overline{\ell},\overline{\mathcal I})$.

  (a) We have
  $J^{c,1}\cap{\mathcal D}(\widetilde{\ell})\subset {\mathcal
    D}(\overline{\ell})$.

  (b) Define $\{c_j\}_{j\in{\mathcal F}}$ and $\{b_j\}_{j\in{\mathcal F}}$ in $[1,p]^{\mathcal F}$, with $c_j<p$ and $b_{j'}<p$ for some $j, j'$, such that$$\widetilde{\ell}\equiv-\sum_{j\in
  {{{\mathcal F}}}}b_jp^j,\quad\quad\quad\quad\overline{\ell}\equiv-\sum_{j\in
  {{{\mathcal F}}}}c_jp^j$$ modulo $\xi{\mathbb
  Z}$. Let $i\in {{{\mathcal F}}}$ such that $b_i=p$.

(b1) If $i\in ({\mathcal E}(\widetilde{\ell})\cup{\mathcal D}(\overline{\ell}))\cap J^{c,1}$ then
$c_i=p$.

(b2) If $i\in {\mathcal D}(\widetilde{\ell})\cap{\mathcal
  E}(\overline{\ell})\cap (J^{c,1}\cup J^+)$ then
$c_i=p-1$.  
\end{lem}

{\sc Proof:} This is not difficult.\hfill$\Box$\\
  
For $S\subset{{{\mathcal F}}}$ and $x\in {{{\mathcal F}}}$ write $\chi^S(x)=1$ if $x\in S$, but $\chi^S(x)=0$ otherwise.

Let $\ell_1,\ell_2\in[0,q-2]$. For $x\in{\mathbb Z}$ let $\check{\nu}_{\ell_1}^{\ell_2}(x)\in{\mathbb Z}$
denote the largest representative of the class of
${\nu}_{\ell_1}^{\ell_2}(x)$ modulo $f{\mathbb Z}$ with
$\check{\nu}_{\ell_1}^{\ell_2}(x)\le x$. For $a<b$
define in ${\mathcal F}$ the subset$$I_{\ell_1}^{\ell_2}[a,b)=[\check{\nu}_{\ell_1}^{\ell_2}(a),\check{\nu}_{\ell_1}^{\ell_2}(b)-1]\quad\mbox{
    if }\check{\nu}_{\ell_1}^{\ell_2}(a)<\check{\nu}_{\ell_1}^{\ell_2}(b),$$but
  $I_{\ell_1}^{\ell_2}[a,b)=\emptyset$ otherwise. Define its subsets$$I_{\ell_1}^{\ell_2}[a,b)_{\vdash}=\left\{\begin{array}{l@{\quad:\quad}l}[0,\check{\nu}_{\ell_1}^{\ell_2}(b)-1]&0\in  I_{\ell_1}^{\ell_2}[a,b)\\I_{\ell_1}^{\ell_2}[a,b)&0\notin I_{\ell_1}^{\ell_2}[a,b)\mbox{ and }(\check{\nu}_{\ell_1}^{\ell_2}(b),b)\ne(0,1)\\\emptyset&\mbox{else }\end{array}\right.$$and $I_{\ell_1}^{\ell_2}[a,b)_{\dashv}=I_{\ell_1}^{\ell_2}[a,b)-I_{\ell_1}^{\ell_2}[a,b)_{\vdash }$. For integers $k\le i<f+k-1$ define in ${\mathcal F}$ the subset$$[k,i]_{\vdash}=\left\{\begin{array}{l@{\quad:\quad}l}[0,i]&0\in [k,i]\\\emptyset&i=-1\\{[k,i]}&\mbox{else }\end{array}\right.$$and $[k,i]_{\dashv}=[k,i]-[k,i]_{\vdash}$.
  
\begin{lem}\label{Josef0} Let $k\le i<f+k-1$ and $\ell_1,\ell_2\in[0,q-2]$
  with $\ell_1\succ_{[k,i]}\ell_2$. Let $\{a_j\}_{j\in{\mathcal F}}$, $\{b_j\}_{j\in{\mathcal F}}$ in $[1,p]^{{\mathcal F}}$ such that $\ell_1\equiv-\sum_{j\in
  {{{\mathcal F}}}}a_jp^j$ and ${\ell}_2\equiv-\sum_{j\in
  {{{\mathcal F}}}}b_jp^j$ modulo $\xi{\mathbb
  Z}$, with $a_j<p$ and $b_{j'}<p$ for some $j, j'$.\footnote{Our conventions imply the understanding $a_j=a_{j+f}$ and  $b_j=b_{j+f}$ for all $j\in{\mathbb Z}$. Similarly, in the following, the arguments and the values of the function $\nu_{\ell_1}^{\ell_2}$ are to be understood modulo $f$.}
  
 (a) We have $\nu_{\ell_1}^{\ell_2}(k)\in\{k,k-1\}$. For $j\in[i+2,k-1]$ we have:
  
$\bullet$ If $\nu_{\ell_1}^{\ell_2}(j)\notin\{j, j-1\}$ then $j=i+2$.

$\bullet$ If $\nu_{\ell_1}^{\ell_2}(j)=j$ then $a_j=b_j$ and $\nu_{\ell_1}^{\ell_2}(j+1)=j+1$.

$\bullet$ If $\nu_{\ell_1}^{\ell_2}(j)\ne j$ then (i) or (ii) holds true:

(i) We have $b_j=p$ and $a_j=1$ and $\nu_{\ell_1}^{\ell_2}(j+1)=j$.

(ii) We have $b_{j}=a_j-1$ and $\nu_{\ell_1}^{\ell_2}(j+1)=j+1$.
 
(b) Let $J$ and $I$ be subsets of ${{{\mathcal F}}}$. For all $s\in[i+2,k-1]$, all $j\in I_{\ell_1}^{\ell_2}[s,s+1)$ assume \begin{gather}\chi^J(s)=\chi^I(j).\label{neukaffee0}\end{gather}For all $j\in I_{\ell_1}^{\ell_2}[k,i+2)$, all $s\in[k,i]$ assume\begin{gather}\chi^J(i+1)=\chi^I(j)=\chi^{J^c}(s).\label{neukaffee1}\end{gather}Then\begin{gather}\sum_{j\in J}a_jp^j-\sum_{j\in J^c}a_jp^j\equiv \sum_{j\in I}b_jp^j-\sum_{j\in I^c}b_jp^j\quad\mbox{ modulo }\xi{\mathbb
    Z}.\label{karfrei0}\end{gather}

(c) Let $J$ and $I$ be subsets of ${{{\mathcal F}}}$. For all $s\in[i+2,k-1]$, all $j'\in I_{\ell_1}^{\ell_2}[s,s+1)_{\dashv}$, all $j\in
  I_{\ell_1}^{\ell_2}[s,s+1)_{\vdash}$
assume\begin{gather}\chi^{I^c}(j')=\chi^{I}(j)=\chi^{J}(s).\label{neukaffee2}\end{gather}For
all $s\in [k,i]_{\vdash}$, all $s'\in [k,i]_{\dashv}$, all $j\in I_{\ell_1}^{\ell_2}[k,i+2)_{\vdash}$, all $j'\in I_{\ell_1}^{\ell_2}[k,i+2)_{\dashv}$ assume\begin{gather}\chi^J(i+1)=\chi^{J^c}(s)=\chi^{J}(s')=\chi^{I^c}(j')=\chi^{I}(j).\label{neukaffee3}\end{gather}Then\begin{gather}\sum_{j\in J}a_jp^j-\sum_{j\in J^c}a_jp^j\equiv \sum_{j\in I}b_jp^j-\sum_{j\in I^c}b_jp^j\quad\mbox{ modulo }(q+1){\mathbb
    Z}\label{karfrei1}\end{gather}
    
\end{lem}

{\sc Proof:} To see (a) is straightforward. To see (b) and (c) we will use (a). Recall
first that for $j\in[i+2,k-1]$ with $\nu_{\ell_1}^{\ell_2}(j)\ne j$ we have
$b_{\nu_{\ell_1}^{\ell_2}(j)}=p$ and $b_{s}=p-1$ for all
$\nu_{\ell_1}^{\ell_2}(j)<s<j$. With this we deduce$$a_sp^s\equiv\sum_{j\in I_{\ell_1}^{\ell_2}[s,s+1)}b_jp^j\quad\mbox{ modulo }\xi{\mathbb
    Z}$$for each $s\in [i+2,k-1]$, as well as$$a_{i+1}p^{i+1}-\sum_{j\in
    [k,i]}a_jp^j\equiv\sum_{j\in I_{\ell_1}^{\ell_2}[k,i+2)}b_jp^j\quad\mbox{ modulo }\xi{\mathbb
    Z}.$$Summing up yields formula (\ref{karfrei0}). Similarly, we find$$a_sp^s\equiv\sum_{j\in I_{\ell_1}^{\ell_2}[s,s+1)_{\vdash}}b_jp^j-\sum_{j\in I_{\ell_1}^{\ell_2}[s,s+1)_{\dashv}}b_jp^j\quad\mbox{ modulo }(q+1){\mathbb
    Z}$$ for each $s\in [i+2,k-1]$, as well as$$a_{i+1}p^{i+1}-\sum_{j\in
    [k,i]_{\vdash}}a_jp^j+\sum_{j\in
    [k,i]_{\dashv}}a_jp^j\equiv\sum_{j\in I_{\ell_1}^{\ell_2}[k,i+2)_{\vdash}}b_jp^j-\sum_{j\in I_{\ell_1}^{\ell_2}[k,i+2)_{\dashv}}b_jp^j\quad\mbox{ modulo }(q+1){\mathbb
    Z}.$$Summing up yields formula (\ref{karfrei1}).\hfill$\Box$\\

  {\bf Remark:} Formulae (\ref{neukaffee0}) and (\ref{neukaffee1}) together are equivalent with\begin{gather}J=\left\{\begin{array}{l@{\quad:\quad}l}J'\cup \{i+1\}&I_{\ell_1}^{\ell_2}[k,i+2)\subset I\\J'\cup [k,i]&\mbox{ else }\end{array}\right.\label{viodurch0}\end{gather}where$$J'=\{s\in[i+2,k-1]\quad|\quad I_{\ell_1}^{\ell_2}[s,s+1)\subset I\}.$$It helps to observe that if $\nu_{\ell_1}^{\ell_2}={\rm id}$ (which is the generic case), $J$ is simply the symmetric difference of $I$ and $[k,i]$, i.e. $J=I\cup [k,i]-(I\cap [k,i])$. Similarly, formulae (\ref{neukaffee2}) and (\ref{neukaffee3}) together are equivalent with\begin{gather}J=\left\{\begin{array}{l@{\quad:\quad}l}J'\cup \{i+1\}\cup[k,i]_{\dashv}&I_{\ell_1}^{\ell_2}[k,i+2)_{\vdash}\subset I\\J'\cup [k,i]_{\vdash}&\mbox{ else }\end{array}\right.\label{viodurch1}\end{gather}where$$J'=\{s\in[i+2,k-1]\quad|\quad I_{\ell_1}^{\ell_2}[s,s+1)_{\vdash}\subset I\}.$$Again, if $\nu_{\ell_1}^{\ell_2}={\rm id}$ this simplifies, since then $I_{\ell_1}^{\ell_2}[x,y)_{\vdash}=[0,y-1]$ and $I_{\ell_1}^{\ell_2}[x,y)_{\dashv}=[x,-1]$ (the latter understood to be empty if $x=0$) if $0\in[x,y-1]$, but $I_{\ell_1}^{\ell_2}[x,y)_{\vdash}=[x,y-1]$ and $I_{\ell_1}^{\ell_2}[x,y)_{\dashv}=\emptyset$ otherwise.
    
\begin{pro}\label{Josef2} Let $({\ell}_1,{\mathcal
    I}_1)\succ_J({\ell}_3,{\mathcal I}_3)$ and $k\in J^-$. There are $J_2$, $\ell_2$ and ${\mathcal I}_2$ such that, with $J_1=[k,k_J^+-1]$, we have$$({\ell}_1,{\mathcal I}_1)\succ_{J_1}({\ell}_2,{\mathcal
    I}_2)\succ_{J_2}({\ell}_3,{\mathcal I}_3).$$
\end{pro}

{\sc Proof:} Put $i=k_J^+-1$, so that $J_1=[k,i]$. Define $\ell_2$ by $\ell_1\succ_{J_1}\ell_2$ and let $\{a_j,
b_j\}_{j\in{{{\mathcal F}}}}$ in $[1,p]$ such that $\ell_1\equiv-\sum_{j\in
  {{{\mathcal F}}}}a_jp^j$ and ${\ell}_2\equiv-\sum_{j\in
  {{{\mathcal F}}}}b_jp^j$ modulo $\xi{\mathbb
  Z}$, with $a_j<p$ and $b_{j'}<p$ for some $j, j'$. Define$${\mathcal I}_2=\nu_{\ell_1}^{\ell_2}({\mathcal
  I}_1-\{i+1\}).$$We then
claim\begin{gather}({\ell}_1,{\mathcal
  I}_1)\succ_{{J}_1}({\ell}_2,{\mathcal
  I}_2).\label{freund2}\end{gather}Indeed, for the requested $J_1^+\cup J_1^-\subset {\mathcal I}_1$ observe $J_1^+\cup J_1^-=\{k,i+1\}$ and $i+1\in J^+\subset {\mathcal I}_1$ and $k\in J^-\subset {\mathcal I}_1$. The requested ${\mathcal I}_2\subset\nu_{\ell_1}^{\ell_2}({\mathcal
  I}_1\cap J_1^{c,1})$ follows from ${\mathcal I}_1-\{i+1\}={\mathcal
  I}_1\cap J_1^{c,1}$. Define$$J_2=\coprod_{j\in J^-\atop j\ne k}I_{\ell_1}^{\ell_2}[j,j_J^+).$$We claim \begin{gather}({\ell}_2,{\mathcal
  I}_2)\succ_{{J}_2}({\ell}_3,{\mathcal I}_3).\label{Josef3}\end{gather}By construction we have$$J^-_2\subset\{\nu_{\ell_1}^{\ell_2}(j) | j\in J^--\{k\}\},\quad\quad\quad J^+_2\subset\{\nu_{\ell_1}^{\ell_2}(j_J^+)|j\in J^--\{k\}\}.$$As $i+1\notin J$ this implies the
requested $J_2^-\cup J_2^+\subset{\mathcal I}_2$. It remains to see\begin{gather}{\mathcal
  I}_3\subset\nu^{{\ell}_3}_{{\ell}_2}({\mathcal I}_2\cap
{J}_2^{c,1}),\label{neuboruspsv0}\\\sum_{j\in J_2}b_jp^j-\sum_{j\in
  J_2^c}b_jp^j\equiv \sum_{j\in J}a_jp^j-\sum_{j\in
  J^c}a_jp^j\quad\mbox{ modulo }\xi{\mathbb
  Z}.\label{boruspsv1}\end{gather}Formula (\ref{boruspsv1}) follows from Lemma
\ref{Josef0}. For
formula (\ref{neuboruspsv0}) it is enough, since ${\mathcal I}_3\subset\nu^{{\ell}_3}_{{\ell}_1}({\mathcal I}_1\cap {J}^{c,1})$ by assumption, to see \begin{gather}\nu_{\ell_1}^{\ell_3}(j)=\nu_{\ell_2}^{\ell_3}\nu_{\ell_1}^{\ell_2}(j)\quad\quad\mbox{
    for all }j\in {J}^{c,1}.\label{00oewvnach}\end{gather}For $j\in
{\mathcal D}(\ell_1)\cap{\mathcal D}(\ell_2)$ as well as $ j\in
{\mathcal E}(\ell_1)\cap{\mathcal E}(\ell_2)$, formula (\ref{00oewvnach}) is
immediate. Next, it is not hard to see that ${J}^{c,1}\cap {\mathcal E}(\ell_1)\cap{\mathcal
  D}(\ell_2)\subset {J}_2^{c,1}$. Thus, Lemma \ref{mario1} says that for $j\in
{J}^{c,1}\cap {\mathcal E}(\ell_1)\cap{\mathcal
  D}(\ell_2)$ we have $j\in {\mathcal
  D}(\ell_3)$, and again formula (\ref{00oewvnach}) is
immediate.

It remains to look at $j\in {J}^{c,1}\cap {\mathcal D}(\ell_1)\cap{\mathcal
  E}(\ell_2)$. We have $\nu_{\ell_1}^{\ell_2}(j)\in {J}_2^{c,1}\cup J_2^+$. Write
$\ell_3=-\sum_{j\in{{{\mathcal F}}}}c_ip^i$ with $c_i\in[1, p]$, not all of them
$=p$.

$\bullet$ If $\nu_{\ell_1}^{\ell_2}(j)\in{\mathcal E}(\ell_2)\cup {\mathcal D}(\ell_3)$
then $\nu_{\ell_1}^{\ell_2}(j)\in {J}_2^{c,1}$. It therefore follows
from Lemma \ref{mario1} that
$c_{\nu_{\ell_1}^{\ell_2}(j)}=b_{\nu_{\ell_1}^{\ell_2}(j)}$, as well as $b_{s}=c_{s}$ for all $\nu_{\ell_1}^{\ell_2}(j)<s<j$. 

$\bullet$ If $\nu_{\ell_1}^{\ell_2}(j)\in{\mathcal D}(\ell_2)\cap {\mathcal
  E}(\ell_3)$ then it follows from Lemma \ref{mario1} that
$c_{\nu_{\ell_1}^{\ell_2}(j)}=p-1$, as well as $b_{s}=c_{s}$ for all $\nu_{\ell_1}^{\ell_2}(j)<s<j$.

In either case, the implication is that
$j$ satisfies formula (\ref{00oewvnach}).\hfill$\Box$\\ 

\begin{pro}\label{0jupp} For any chain $({\ell}_1,{\mathcal I}_1)\succ_{J_1}({\ell}_2,{\mathcal I}_2)\succ_{J_2}({\ell}_3,{\mathcal I}_3)$ in ${\mathbb D}$ with $|J_1^-|=1$ there is a $J$ with $({\ell}_1,{\mathcal I}_1)\succ_{J}({\ell}_3,{\mathcal I}_3)$. \end{pro}

{\sc Proof:} Write $J_1=[k,i]$ with $k\le i<f+k-1$; thus,
$J_1^{c,1}=[i+2,k]$. Our assumptions imply$$J^-_2\cup J^+_2\subset{\mathcal
  I}_2\subset \nu_{\ell_1}^{\ell_2}({\mathcal I}_1\cap J_1^{c,1}).$$We first
claim \begin{gather}I_{\ell_1}^{\ell_2}[s,s+1)\cap
    J_2\in\{\emptyset,I_{\ell_1}^{\ell_2}[s,s+1)\}\quad\mbox{ for each }s\in[i+2,k-1],\label{osterfrei0}\\I_{\ell_1}^{\ell_2}[k,i+2)\cap
    J_2\in\{\emptyset,I_{\ell_1}^{\ell_2}[k,i+2)\}.\label{osterfrei1}\end{gather}Formula (\ref{osterfrei0}) follows from $J^-_2\cup J^+_2\subset{\rm
  im}(\nu_{\ell_1}^{\ell_2})$. Formula (\ref{osterfrei1}) follows from the
more precise $J^-_2\cup J^+_2\subset\nu_{\ell_1}^{\ell_2}([i+2,k])$. Given formulae
(\ref{osterfrei0}) and (\ref{osterfrei1}) we find a unique $J\subset{{{\mathcal F}}}$
which satisfies formulae (\ref{neukaffee0}) and (\ref{neukaffee1}) with
$I=J_2$, cf. formula (\ref{viodurch0}). By Lemma \ref{Josef0} it satisfies$$\sum_{j\in J_2}b_jp^j-\sum_{j\in
  J_2^c}b_jp^j\equiv \sum_{j\in J}a_jp^j-\sum_{j\in
  J^c}a_jp^j\quad\mbox{ modulo }\xi{\mathbb Z}.$$To see that $J$ works as
desired it remains to see \begin{gather}J^-\cup J^+\subset{\mathcal I}_1,\label{0juppmaerzsam1}\\{\mathcal I}_3\subset\nu_{\ell_1}^{\ell_3}({\mathcal I}_1\cap J^{c,1}).\label{0juppmaerzsam2}\end{gather}
    
Formula (\ref{0juppmaerzsam1}) follows from the assumptions $\{k,i+1\}=J^+_1\cup
J^-_1\subset{\mathcal I}_1$ and $J^-_2\cup J^+_2\subset{\mathcal
  I}_2\subset \nu_{\ell_1}^{\ell_2}({\mathcal I}_1)$, by the very definition
of $J$.

To see formula (\ref{0juppmaerzsam2}) it will be enough, as ${\mathcal I}_3\subset
\nu_{\ell_2}^{\ell_3}({\mathcal I}_2\cap J_2^{c,1})$ and ${\mathcal I}_2\subset
\nu_{\ell_1}^{\ell_2}({\mathcal I}_1\cap J_1^{c,1})$ by assumption, to
see$$\nu_{\ell_2}^{\ell_3}(\nu_{\ell_1}^{\ell_2}({\mathcal I}_1\cap
J_1^{c,1})\cap J_2^{c,1})\subset \nu_{\ell_1}^{\ell_3}({\mathcal I}_1\cap
J^{c,1}),$$and as $J_1^{c,1}=[i+2,k]$ it will be enough to
see\begin{gather}\nu_{\ell_2}^{\ell_3}\nu_{\ell_1}^{\ell_2}(j)=\nu_{\ell_1}^{\ell_3}(j)\quad\mbox{
  for all }j\in[i+2,k]\mbox{ with }\nu_{\ell_1}^{\ell_2}(j)\in
J_2^{c,1},\label{0juppjosefstag}\\ [i+2,k]\cap
(\nu_{\ell_1}^{\ell_2})^{-1}(J_2^{c,1})\subset J^{c,1}.\label{ost0juppjosefstag}\end{gather}For $j\in{\mathcal D}(\ell_1)\cap{\mathcal D}(\ell_2)$ as well as for
$j\in{\mathcal E}(\ell_1)\cap{\mathcal E}(\ell_2)$, formula (\ref{0juppjosefstag}) is
clear. For $j\in {\mathcal E}(\ell_1)\cap{\mathcal D}(\ell_2)$ we have
$j=\nu_{\ell_1}^{\ell_2}(j)$ and by our assumption this belongs to
$J_2^{c,1}$, therefore Lemma \ref{mario1} shows $j\in{\mathcal D}(\ell_3)$,
and this again implies formula (\ref{0juppjosefstag}).

It remains to consider the case $j\in {\mathcal D}(\ell_1)\cap{\mathcal
  E}(\ell_2)$. We then have $j\ne\nu_{\ell_1}^{\ell_2}(j)\in {J}_2^{c,1}$. Write
$\ell_3=-\sum_{j\in{{{\mathcal F}}}}c_ip^i$ with $c_i\in[1, p]$, not all of them
$=p$.

$\bullet$ If $\nu_{\ell_1}^{\ell_2}(j)\in{\mathcal E}(\ell_2)\cup {\mathcal D}(\ell_3)$
then it follows
from Lemma \ref{mario1} that $c_{\nu_{\ell_1}^{\ell_2}(j)}=b_{\nu_{\ell_1}^{\ell_2}(j)}$, as well as $b_{s}=c_{s}$ for all $\nu_{\ell_1}^{\ell_2}(j)<s<j$.

$\bullet$ If $\nu_{\ell_1}^{\ell_2}(j)\in{\mathcal D}(\ell_2)\cap {\mathcal
  E}(\ell_3)$ then it follows from Lemma \ref{mario1} that
$c_{\nu_{\ell_1}^{\ell_2}(j)}=p-1$, as well as $b_{s}=c_{s}$ for all $\nu_{\ell_1}^{\ell_2}(j)<s<j$.

In either case, the implication is that
$j$ satisfies formula (\ref{0juppjosefstag}).

Formula (\ref{ost0juppjosefstag}) follows from the definition of $J$.\hfill$\Box$\\

  \begin{satz}\label{tobbu} (a) Given $(\widetilde{\ell},\widetilde{\mathcal I})\succ_J(\overline{\ell},\overline{\mathcal I})$ in ${\mathbb D}$, there are $m\ge1$ and $({\mathcal I}_n,\ell_n)$ for $0\le n\le m$, and $J_n$ with $|J_n^-|=1$ for $1\le n\le m$, such that $(\widetilde{\ell},\widetilde{\mathcal I})=({\ell}_0,{\mathcal I}_0)$ and $(\overline{\ell},\overline{\mathcal I})=({\ell}_m,{\mathcal I}_m)$ and $({\ell}_{n-1},{\mathcal I}_{n-1})\succ_{J_n}({\ell}_{n},{\mathcal I}_n)$ for all $1\le n\le m$.
  
  (b) Given $({\ell}_1,{\mathcal I}_1)\succ_{J_1}({\ell}_2,{\mathcal I}_2)\succ_{J_2}({\ell}_3,{\mathcal I}_3)$ in ${\mathbb D}$ there is a $J$ with $({\ell}_1,{\mathcal I}_1)\succ_{J}({\ell}_3,{\mathcal I}_3)$.
\end{satz}

{\sc Proof:} Statement (a) follows from Proposition \ref{Josef2} by induction. To see statement (b) we apply statement (a) to $({\ell}_1,{\mathcal I}_1)\succ_{J_1}({\ell}_2,{\mathcal I}_2)$ and then use an induction argument to reduce to the case where $|J_1^-|=1$; this case was handled in Proposition \ref{0jupp}.\hfill$\Box$\\

{\bf Definition:} (a) For $h\in{\mathbb Z}$ and $\ell\in[0,q-2]$ and $J\subset{{{\mathcal F}}}$ we declare the relation $\ell\succ^+_Jh$ to hold if there are $1\le a_j\le p$ for all $j\in{\mathcal F}$, not all of them $=p$, such that$${\ell}\equiv-\sum_{j\in{{{\mathcal F}}}}a_jp^j\mbox{ modulo }\xi{\mathbb Z},\quad\quad\quad\quad h\equiv\sum_{j\in J}a_jp^j-\sum_{j\in J^c}a_jp^j\mbox{ modulo }(q+1){\mathbb Z}.$$
  
(b) For $h\in{\mathbb Z}$ denote by $[[h]]$ its class in ${\mathbb Z}/(q+1){\mathbb Z}$. For $(\ell,{\mathcal I})\in{\mathbb D}$ and $J\subset{{{\mathcal F}}}$ we declare the relation $(\ell,{\mathcal I})\succ_J[[h]]$ to hold if and only if $${\mathcal I}^c\cap\{0\}\subset J^+\cup J^-\subset{\mathcal I}\cup\{0\}\quad\quad\mbox{ and }\quad\quad\ell\succ^+_Jh.$$

(c) On $\breve{\mathbb D}={\mathbb D}\coprod {\mathbb Z}/(q+1){\mathbb Z}$ define the relations $\succ_J$ for $J\subset{\mathcal F}$ by taking the union of the relations $\succ_J$ on ${\mathbb D}$ and all the $(\ell,{\mathcal I})\succ_J[[h]]$ with $(\ell,{\mathcal I})\in{\mathbb D}$ and $[[h]]\in {\mathbb Z}/(q+1){\mathbb Z}$. 

\begin{satz}\label{ostersa} Let $h\in{\mathbb Z}$ and $(\ell,{\mathcal I})\in{\mathbb D}$, let $J\subset{{{\mathcal F}}}$. We have $(\ell,{\mathcal I})\succ_J[[h]]$ if and only if there are some $0\le i<f$, some $\overline{\ell}$ and $I$ such that $(\ell,{\mathcal I})\succ_I(\overline{\ell},\{i\})$ and $(\overline{\ell},\{i\}) \succ_{K}[[h]]$, where either $K=\{j\,|\,0\le j<i\}$ or $K=\{j\,|\,i\le j<f\}$.

\end{satz}

{\sc Proof:} {\it Claim 1:} For any chain $({\ell}_1,{\mathcal I}_1)\succ_{J_1}({\ell}_2,{\mathcal I}_2)\succ_{J_2}[[h]]$ with $|J_1^-|=1$ there is a $J$ with $({\ell}_1,{\mathcal I}_1)\succ_{J}[[h]]$.

{\it Claim 2:} Let $({\ell}_1,{\mathcal
    I}_1)\succ_J[[h]]$ and $k\in J^-$. There are $J_2$, $\ell_2$ and ${\mathcal I}_2$ such that, with $J_1=[k,k_J^+-1]$, we have $({\ell}_1,{\mathcal I}_1)\succ_{J_1}({\ell}_2,{\mathcal
    I}_2)\succ_{J_2}[[h]]$.

{\it Conclusion:} Given claims 1 and 2, proceed as in the proof of Theorem \ref{tobbu}.

{\it Proof of claim 1:} This is similar to the proof of Proposition \ref{0jupp}. Write $J_1=[k,i]$ with $k\le i<f+k-1$; thus,
$J_1^{c,1}=[i+2,k]$. Our assumptions imply\begin{gather}{\mathcal
  I}_2^c\cap\{0\}\subset J^-_2\cup J^+_2\subset{\mathcal
  I}_2\cup\{0\}\quad\mbox{ and }\quad{\mathcal
  I}_2\subset \nu_{\ell_1}^{\ell_2}({\mathcal I}_1\cap J_1^{c,1}).\label{nguyen}\end{gather}We first
claim \begin{gather}(I_{\ell_1}^{\ell_2}[s,s+1)_{\dashv}\cap
    J_2)\cup(I_{\ell_1}^{\ell_2}[s,s+1)_{\vdash}\cap
    J_2)\in\{I_{\ell_1}^{\ell_2}[s,s+1)_{\dashv},I_{\ell_1}^{\ell_2}[s,s+1)_{\vdash}\}\label{freiosterfrei0}\end{gather}for each $s\in[i+2,k-1]$, and similarly\begin{gather}(I_{\ell_1}^{\ell_2}[k,i+2)_{\dashv}\cap
    J_2)\cup(I_{\ell_1}^{\ell_2}[k,i+2)_{\vdash}\cap
      J_2)\in\{I_{\ell_1}^{\ell_2}[k,i+2)_{\dashv},I_{\ell_1}^{\ell_2}[k,i+2)_{\vdash}\}.\label{freiosterfrei1}\end{gather}Indeed, first notice that formulae (\ref{osterfrei0}) and
(\ref{osterfrei1}) hold again in the present context. If both
$I_{\ell_1}^{\ell_2}[s,s+1)_{\dashv}$ and $I_{\ell_1}^{\ell_2}[s,s+1)_{\vdash}$ are
    non-empty, then there is no $x\in J_1^{c,1}$ with
    $\nu_{\ell_1}^{\ell_2}(x)=0$. Thus,
    $0\in\nu_{\ell_1}^{\ell_2}(J_1^{c,1})^c$, hence $0\in {\mathcal  I}_2$ by assumption, hence $0\in  J^-_2\cup J^+_2$ by assumption. But this
    (together with formula (\ref{osterfrei0})) implies formula
    (\ref{freiosterfrei0}). Formula (\ref{freiosterfrei1}) is proven in the
    same way.

Given formulae
(\ref{freiosterfrei0}) and (\ref{freiosterfrei1}) we find a unique $J\subset{{{\mathcal F}}}$
which satisfies formulae (\ref{neukaffee2}) and (\ref{neukaffee3}) with
$I=J_2$, cf. formula (\ref{viodurch1}). By Lemma \ref{Josef0} it satisfies$$\sum_{j\in J}a_jp^j-\sum_{j\in J^c}a_jp^j\equiv \sum_{j\in I}b_jp^j-\sum_{j\in I^c}b_jp^j\quad\mbox{ modulo }(q+1){\mathbb
    Z}.$$To see that $J$ works as
  desired it remains to see\begin{gather}{\mathcal I}_1^c\cap\{0\}\subset J^-\cup J^+\subset{\mathcal I}_1\cup\{0\}.\label{osterdi0}\end{gather}Similarly as in the proof of Proposition \ref{0jupp} this follows now from $\{k,i+1\}=J_1^+\cup J_1^{-}\subset{\mathcal I}_1$ and formula (\ref{nguyen}).

  {\it Proof of claim 2:} This is similar to the proof of Proposition \ref{Josef2}. Put $i=k_J^+-1$, so that $J_1=[k,i]$. Define $\ell_2$ by $\ell_1\succ_{J_1}\ell_2$ and let $\{a_j,
b_j\}_{j\in{{{\mathcal F}}}}$ in $[1,p]$ such that $\ell_1\equiv-\sum_{j\in
  {{{\mathcal F}}}}a_jp^j$ and ${\ell}_2\equiv-\sum_{j\in
  {{{\mathcal F}}}}b_jp^j$ modulo $\xi{\mathbb
  Z}$, with $a_j<p$ and $b_{j'}<p$ for some $j, j'$. Define$${\mathcal I}_2=\nu_{\ell_1}^{\ell_2}({\mathcal
  I}_1-\{i+1\}).$$We then
claim\begin{gather}({\ell}_1,{\mathcal
  I}_1)\succ_{{J}_1}({\ell}_2,{\mathcal
  I}_2).\label{OOfreund2}\end{gather}Indeed, for the requested $J_1^+\cup J_1^-\subset {\mathcal I}_1$ observe $J_1^+\cup J_1^-=\{k,i+1\}$ and $i+1\in J^+\subset {\mathcal I}_1$ and $k\in J^-\subset {\mathcal I}_1$. The requested ${\mathcal I}_2\subset\nu_{\ell_1}^{\ell_2}({\mathcal
  I}_1\cap J_1^{c,1})$ follows from ${\mathcal I}_1-\{i+1\}={\mathcal
  I}_1\cap J_1^{c,1}$. Define$$J_2=I_{\ell_1}^{\ell_2}[k,i+2)_{\dashv}\coprod\coprod_{j\in J^-\atop j\ne k}I_{\ell_1}^{\ell_2}[j,j_J^+)_{\vdash}.$$We claim \begin{gather}({\ell}_2,{\mathcal
  I}_2)\succ_{{J}_2}[[h]].\label{OOJosef3}\end{gather}Indeed, Lemma \ref{Josef0} with $I=J_2$ yields $$\sum_{j\in J}a_jp^j-\sum_{j\in J^c}a_jp^j\equiv \sum_{j\in J_2}b_jp^j-\sum_{j\in J_2^c}b_jp^j\quad\mbox{ modulo }(q+1){\mathbb
    Z},$$and similarly as before we also see$${\mathcal I}_2^c\cap\{0\}\subset J_2^-\cup J_2^+\subset{\mathcal I}_2 \cup\{0\}.$$\hfill$\Box$\\

\section{The combinatorics of degenerating families of $(\varphi,\Gamma)$-modules}

\label{colook}

{\bf Definition:} (a) On ${\mathbb E}$ (cf. subsection \ref{stackstrataset}) we define the binary relation $\rhd$ by
declaring$$(\widetilde{\ell},\widetilde{u},\widetilde{\mathcal
  I})\rhd(\overline{\ell},\overline{u},\overline{\mathcal I})$$if and only if (at least) one of the following scenarios is fulfilled:

(I) $(\widetilde{\ell},\widetilde{u})=(\overline{\ell},\overline{u})$ and $\overline{\mathcal I}\subset
\widetilde{\mathcal I}$.

(II) $\widetilde{\ell}\equiv \widetilde{u}-\overline{u}\equiv -\overline{\ell}$ modulo $\xi{\mathbb Z}$ and $\overline{\mathcal I}=\emptyset$.

(III) There is some $i\ge 2f$ with $\widetilde{\ell}+2p^i\in\xi{\mathbb Z}$ and $\overline{\ell}=0$ and
$\overline{u}\equiv p^i+\widetilde{u}$ modulo $\xi{\mathbb Z}$ and $\overline{\mathcal I}=\emptyset$ and $\widetilde{\mathcal
  I}=\{\Pi(i)\}$.

(IV) There are $2f\le i_1<i_2<i_1+f$ such that the conditions (\ref{1rook2}), (\ref{2rook2}), (\ref{ftklaus}), (\ref{armand}) and (\ref{armand1}) are fulfilled. 

(b) Denote the class of $h\in{\mathbb Z}$ in ${\mathbb Z}/(q^2-1){\mathbb Z}$ (and hence in $\breve{{\mathbb E}}$, cf. subsection \ref{stackstrataset}) by $[h]$. Extend the binary relation $\rhd$ from ${\mathbb E}$ to $\breve{{\mathbb E}}$ as follows: We declare $$[h] \rhd [qh]$$for all $h\in{\mathbb Z}$. We further declare$$(\widetilde{\ell},\widetilde{u},\widetilde{\mathcal
  I})\rhd [h]$$if and only if the following hold. Let $1\le a_i\le p$ for $0\le i<f$ be such that $\widetilde{\ell}\equiv-\sum_{i=0}^{f-1}a_ip^i$ modulo $\xi{\mathbb Z}$, assuming $a_i\ne p$ for at least one $i$. Write $a_i=a_{i+f}$ for all $i\in{\mathbb Z}$. We demand $\widetilde{\ell}\ne0$ and that there is some $i\ge0$ with $\widetilde{\mathcal I}=\{\Pi(i)\}$, with $\widetilde{\ell}+2p^i\notin \xi{\mathbb Z}$, with\begin{gather}h\equiv -\sum_{j=i}^{i+f-1}a_jp^j\quad\mbox{ modulo }(q+1){\mathbb Z}\label{asch2}\end{gather}and with\begin{gather}\frac{h+\sum_{j=i}^{i+f-1}a_jp^j}{q+1}\equiv -\widetilde{u}\quad\mbox{ modulo }\xi{\mathbb Z}.\label{asch1}\end{gather} 
  
{\bf Remark:} Given $i\ge0$, define in ${{{\mathcal F}}}$ the subset\begin{gather}J=\left\{\begin{array}{l@{\quad:\quad}l}\Pi([(2s+1)f,f+i-1])&2sf\le i<(2s+1)f\\\Pi([i, 2sf-1])&(2s-1)f\le i<2sf\end{array}\right..\label{kruthilf}\end{gather}with the appropriate $s\in{\mathbb Z}$. Formula (\ref{asch2}) is equivalent with\begin{gather}h\equiv\sum_{j\in J}a_jp^j-\sum_{j\in J^c}a_jp^j\quad\mbox{ modulo }(q+1){\mathbb Z}.\label{stoch}\end{gather}
  
{\bf Remark:} Following Lemma \ref{fernandodrei} below we explain how the relation $\rhd$ (and its transitive closure $\rhd\rhd$) relates to the relations $\succ_J$ of section \ref{combi}.
  
\begin{kor}\label{palmmon} Let $(\widetilde{\ell},\widetilde{u},\widetilde{\mathcal
  I})\in{\mathbb E}$ and a $(\varphi,\Gamma)$-module $\overline{\bf D}$ be such that one of the following occurs:

(A) $\overline{\bf D}\in{\mathcal V}(\overline{\ell},\overline{u},\overline{\mathcal I})$ and $(\widetilde{\ell},\widetilde{u},\widetilde{\mathcal
  I})\rhd(\overline{\ell},\overline{u},\overline{\mathcal I})$ for some $(\overline{\ell},\overline{u},\overline{\mathcal I})\in {\mathbb E}$ with $\overline{\ell}\ne0$.

(B) $\overline{\bf D}\in{\mathcal V}(\overline{\ell},\overline{u},\overline{\mathcal I})$ and $(\widetilde{\ell},\widetilde{u},\widetilde{\mathcal
  I})\rhd(\overline{\ell},\overline{u},\overline{\mathcal I})$ for some $(\overline{\ell},\overline{u},\overline{\mathcal I})\in {\mathbb E}$ with $\overline{\ell}=0$.

(C) $\overline{\bf D}\in{\mathcal V}(h)$ and $(\widetilde{\ell},\widetilde{u},\widetilde{\mathcal
  I})\rhd [h]$ for some $h\in{\mathbb Z}$.

Then there exists a $(\varphi,\Gamma)$-module ${\bf D}$ over $B$, with $B=k[\tau]$ in cases (A) and (C), and with $B=k[\tau][\frac{1}{\tau^2-\epsilon}]$ for some $\epsilon\in k^{\times}$ in case (B), with the following property: We have ${\bf D}\otimes_{B}k\cong \overline{\bf D}$ via $\tau\mapsto 0\in k$, and ${\bf D}\otimes_{B}k\in {\mathcal V}(\widetilde{\ell},\widetilde{u},\widetilde{\mathcal
  I})$ via $\tau\mapsto a$ for each $a\in k^{\times}$ (resp. each $a\in k^{\times}$ with $a^2\ne\epsilon$ in case (B)).
\end{kor}

  {\sc Proof:} First consider (A) and (B). In case (I)
  (in the definition of $(\widetilde{\ell},\widetilde{u},\widetilde{\mathcal
    I})\rhd(\overline{\ell},\overline{u},\overline{\mathcal I})$) the claim is trivial. In case (II) we may assume, after composing with an instance of case (I), that $\widetilde{\mathcal
    I}=\emptyset$. Then it corresponds to switching the two characters appearing in a two dimensional reducible split Galois representation. In case (III) apply Corollary \ref{cusanus00}\footnote{This is a subcase of (B), but here we may even take $B=k[\tau]$, not only $B=k[\tau][\frac{1}{\tau^2-\epsilon}]$ for some $\epsilon\in k^{\times}$}, in case (IV) apply Corollary \ref{cusanus}.
  
Now consider (C). Formulae (\ref{ashwed}), (\ref{asch2}) and (\ref{asch1}) together show$$h\equiv\widetilde{\ell}_{[i]}-\sigma^{\widetilde{\ell}}(i)p^i\xi-(q+1)\widetilde{u}\quad\mbox{ modulo }(q^2-1){\mathbb Z}.$$Thus, we may again apply Corollary \ref{cusanus00}.\hfill$\Box$\\

{\bf Remark:} Let $(\widetilde{\ell},\widetilde{u},\widetilde{\mathcal
  I})\in {\mathbb E}$. Recall (section \ref{stackstrataset}) that $(\widetilde{\ell},\widetilde{u},\widetilde{\mathcal
  I})$ corresponds to a family ${\bf D}(\widetilde{\ell},\widetilde{u},\widetilde{\mathcal I})$ of $(\varphi,\Gamma)$-modules over $({\mathbb G}_m)^2\times{\mathbb A}(\widetilde{\ell},\widetilde{u},\widetilde{\mathcal
  I})$ for some affine space ${\mathbb A}(\widetilde{\ell},\widetilde{u},\widetilde{\mathcal
  I})$ of dimension $|\widetilde{\mathcal I}|$, parametrizing certain ${\mathcal G}_F$-representations which, when restricted to ${\mathcal I}_F$, are extensions of $\omega^{\widetilde{\ell}-\widetilde{u}}$ by $\omega^{-\widetilde{u}}$. Recall also that $\widetilde{\mathcal I}$ is in bijection with a spanning set of lines in ${\mathbb A}(\widetilde{\ell},\widetilde{u},\widetilde{\mathcal
  I})$. We recapitulate:

${\bullet}$ (Degeneration) For each $\Pi(i_1)\in \widetilde{\mathcal I}$, 'going to infinity' along the line in ${\mathbb A}(\widetilde{\ell},\widetilde{u},\widetilde{\mathcal
  I})$ corresponding to $\Pi(i_1)$ gives rise to some $e\in \breve{{\mathbb E}}$ with $(\widetilde{\ell},\widetilde{u},\widetilde{\mathcal
  I})\rhd e$, in such a way that the said family of ${\mathcal G}_F$-representations degenerates into ${\mathcal G}_F$-representations belonging to $e$. Namely, if $|\widetilde{\mathcal
  I}|=1$ then define $e=[h]$ with $i$ as in the definition of $(\widetilde{\ell},\widetilde{u},\widetilde{\mathcal
  I})\rhd [h]$ such that $\Pi(i_1)=\Pi(i)$. If however $|\widetilde{\mathcal
  I}|\ge 2$ then let $2f\le i_1<i_2<i_1+f$ (with $i_1$ lifting the given $\Pi(i_1)$) such that $\Pi(i_2)\in \widetilde{\mathcal
  I}$ but $\Pi(i)\notin \widetilde{\mathcal
  I}$ for all $i_2<i<i_1+f$. This yields a well defined $e=(\overline{\ell},\overline{u},\overline{\mathcal
  I})$ as before.

${\bullet}$ (Deformation) By Corollary \ref{palmmon} we can also go the other way: For given $(\widetilde{\ell},\widetilde{u},\widetilde{\mathcal
  I})\rhd e$ in $\breve{\mathbb E}$, each ${\mathcal G}_F$-representation belonging to $e$ may be deformed into a family which generically belongs to $(\widetilde{\ell},\widetilde{u},\widetilde{\mathcal
  I})$.\\

{\bf Definition:} (a) Denote by $\rhd\rhd$ the transitive closure of the relation $\rhd$ on $\breve{{\mathbb E}}$.

(b) For integers $u$ and $0\le\ell<\xi$ define the tuples $a_{\bullet}=(a_i)_{0\le i<f}\in [1,p]^f$ and $d_{\bullet}=(d_i)_{0\le i<f}\in [0,p-1]^f-\{(p-1,\ldots,p-1)\}$ by asking the modulo $\xi{\mathbb Z}$-congruences$$\ell\equiv-\sum_{i=0}^{f-1}a_ip^i\quad\mbox{ and }\quad u\equiv-\sum_{i=0}^{f-1}(a_i+d_i)p^i$$where if $\ell=\frac{(p-2)\xi}{p-1}$ choose $a_{\bullet}=(1,\ldots,1)$. Then put$$V(\ell,u):=V_{d_{\bullet},a_{\bullet}}:=\bigotimes_{i=0}^{f-1}({\rm det}^{d_i}\otimes_{{\mathbb F}_q}{\rm Sym}^{a_i-1}({\mathbb F}_q^2))\otimes_{{\mathbb F}_q,\omega^i}{\overline{\mathbb F}_p}.$$This is an irreducible representation of ${\rm GL}_2({\mathbb F}_q)$; as such, it uniquely determines the pair $(a_{\bullet},d_{\bullet})$, and hence also $\ell$ and $u$ modulo $\xi{\mathbb Z}$.\footnote{The case $\ell=0$ corresponds to $a_{\bullet}=(p-1,\ldots,p-1)$ and $\omega^{\ell}$ being the trivial character. The case $\ell=\frac{(p-2)\xi}{p-1}$ corresponds to $a_{\bullet}=(1,\ldots,1)$ and $\omega^{\ell}$ being the cyclotomic character.} For integers $u$ again now define $d_{\bullet}=(d_i)_{0\le i<f}\in [0,p-1]^f-\{(p-1,\ldots,p-1)\}$ by asking the modulo $\xi{\mathbb Z}$-congruence$$u\equiv-\sum_{i=0}^{f-1}(p+d_i)p^i\equiv-\sum_{i=0}^{f-1}(1+d_i)p^i$$and put$$V(\frac{(p-2)\xi}{p-1},u)_{\rm St}:=V_{d_{\bullet},(p,\ldots,p)}:=\bigotimes_{i=0}^{f-1}({\rm det}^{d_i}\otimes_{{\mathbb F}_q}{\rm Sym}^{p-1}({\mathbb F}_q^2))\otimes_{{\mathbb F}_q,\omega^i}{\overline{\mathbb F}_p}.$$

(c) Let $W$ be a two dimensional ${\mathcal G}_F$-representation with corresponding $(\varphi,\Gamma)$-module ${\bf D}$. Let $e=[h]$ resp. $e=(\overline{\ell},\overline{u},\overline{\mathcal I})$ be such that ${\bf D}\in {\mathcal V}(h)$ resp. ${\bf D}\in{\mathcal V}^0(\overline{\ell},\overline{u},\overline{\mathcal I})$ (cf. formula (\ref{fernandovertex})). Define$${\mathbb W}'_{\rhd}(W)=\{V({\ell},u)\,\,|\,\,({\ell},u,{\mathcal F})\rhd\rhd e\}$$and then ${\mathbb W}_{\rhd}(W)={\mathbb W}'_{\rhd}(W)$ if there is no $V(\frac{(p-2)\xi}{p-1},u)\in {\mathbb W}'_{\rhd}(W)$, but ${\mathbb W}_{\rhd}(W)={\mathbb W}'_{\rhd}(W)\cup \{V(\frac{(p-2)\xi}{p-1},u)_{\rm St}\}$ if $V(\frac{(p-2)\xi}{p-1},u)\in {\mathbb W}'_{\rhd}(W)$ (there is at most one such). 

(d) Let $W$ and ${\bf D}$ be as in (c). Consider the stack ${\mathcal X}_{2,{\rm red}}$ described in \cite{eg}, or rather its Lubin Tate version described in \cite{pham}.\footnote{By the main result of \cite{pham}, these two stacks are in fact isomorphic.} It parametrizes \'{e}tale $(\varphi,\Gamma)$-modules over $\overline{\mathbb F}_p$ of rank two. Its irreducible components are labelled by the irreducible representations of ${\rm GL}_2({\mathbb F}_q)$ over $\overline{\mathbb F}_p$. Define ${\mathbb W}_{\rm geom}(W)$ to be the set of labels of those irreducible components ${\mathcal C}$ of ${\mathcal X}_{2,{\rm red}}$ such that ${\bf D}$ corresponds to a $k$-point of ${\mathcal C}$.

(e) For $W$ as in (c) let ${\mathbb W}_{BDJ}(W)$ denote the set of Serre weights attached to $W$ in \cite{bdj}.

\begin{satz}\label{freimai} ${\mathbb W}_{\rm geom}(W)={\mathbb W}_{\rhd}(W)$.\end{satz}

{\sc Proof:} Its proof shows that Corollary \ref{palmmon} also admits family versions, stating the existence of a deformation of a given subfamily of ${\bf D}(\overline{\ell},\overline{u},\overline{\mathcal I})$, over an irreducible subscheme of $({\mathbb G}_m)^2\times {\mathbb A}(\overline{\ell},\overline{u},\overline{\mathcal I})$, into a family of $(\varphi,\Gamma)$-modules which generically is a subfamily of ${\bf D}(\widetilde{\ell},\widetilde{u},\widetilde{\mathcal I})$ (over an irreducible subscheme of $({\mathbb G}_m)^2\times {\mathbb A}(\widetilde{\ell},\widetilde{u},\widetilde{\mathcal I})$). This becomes clear from looking at the proof of Corollary \ref{cusanus} (which underlies Corollary \ref{palmmon}): The coefficients in the generators of $\nabla$ --- namely, the coefficients of $C$, $D$, or equivalently of $\widetilde{A}$, $\widetilde{B}$, in subsection \ref{osterfreitag} (or in subsection \ref{keydefo}: the scalars $\alpha_i$) --- must now be allowed to vary (rather than being fixed in $k$) over a suitable dense open subspace of an affine $k$-space.

Now let $V\in{\mathbb W}_{\rhd}(W)$. There is a chain $$(\ell_0,u_0,{\mathcal I}_0)\rhd (\ell_1,u_1,{\mathcal I}_1)\rhd\ldots\rhd (\ell_r,u_r,{\mathcal I}_r)\rhd e$$(some $r$) with $V=V(\ell_0,u_0)$, with ${\mathcal I}_0={\mathcal F}$ and with $e$ attached to $W$ as in the definition of $\rhd$. We claim that there are affine $k$-schemes $X_i$, closed subschemes $Y_i\subset X_i$ with dense complement, and $(\varphi,\Gamma)$-modules ${\bf D}_i$ over the $X_i$, for $0\le i\le r+1$, with the following properties: $X_{r+1}={\rm Spec}(k)$ and $Y_{r+1}=\emptyset$; for $1\le i\le r+1$ there are isomorphisms $X_i-Y_i\cong Y_{i-1}$ together with isomorphisms between the respective $(\varphi,\Gamma)$-modules $${\bf D}_i|_{X_i-Y_i}\cong {\bf D}_{i-1}|_{Y_{i-1}}$$and ${\bf D}_i|_{X_i-Y_i}$ for $0\le i\le r$ is a subfamily of the family ${\bf D}(\ell_i,u_i,{\mathcal I}_i)$ of subsection \ref{stackstrataset}.

Indeed, this follows from repeated application the family version of Corollary \ref{palmmon} mentioned above.

Given this, the definition of ${\mathcal X}_{2,{\rm red}}$ implies $V\in{\mathbb W}_{\rm geom}(W)$, since ${\bf D}_0|_{X_0-Y_0}$ is a subfamily of ${\bf D}(\ell_0,u_0,{\mathcal F})$, hence ${\bf D}_0$ lives over the component of ${\mathcal X}_{2,{\rm red}}$ labelled by $V=V(\ell_0,u_0)$.

We have shown the inclusion ${\mathbb W}_{\rhd}(W)\subset{\mathbb W}_{\rm geom}(W)$. The reverse inclusion now simply follows from $|{\mathbb W}_{\rhd}(W)|=|{\mathbb W}_{\rm geom}(W)|$, for which see Corollary \ref{movorchrihi} below.

Alternatively, if one wants to avoid usage of the results from \cite{dediro}, \cite{ecgm}, \cite{gls} (which go into Corollary \ref{movorchrihi}), one might infer the reverse inclusion from seeing that the preceding constructions capture {\it all} $(\varphi,\Gamma)$-module deformations/degenerations. But indeed, they capture the degeneration of the $(\widetilde{\ell},\widetilde{u},\widetilde{\mathcal
  I})$-family, parametrized by the product of lines indexed by $\widetilde{\mathcal
  I}$, by 'going to infinity' along {\it any} of these lines. \hfill$\Box$\\

The wit in the definition of ${\mathbb W}_{\rhd}(W)$ lies in Theorem \ref{marienmai} below. Define the map\begin{gather}\breve{\mathbb E}\longrightarrow\breve{\mathbb D},\quad e\mapsto[[e]]\label{fast}\end{gather}by putting $[[({\ell},{u},{\mathcal
  I})]]= 
({\ell},{\mathcal
  I})\in {\mathbb D}$ for $({\ell},{u},{\mathcal
  I})\in {\mathbb E}$ and putting $[[[h]]]=[[h]]\in {\mathbb Z}/(q+1){\mathbb Z}$ for $[h]\in{\mathbb Z}/(q^2-1){\mathbb Z}$.

 \begin{lem}\label{fernandodrei} (a) For $\widetilde{e}\in{\mathbb E}$ and $\overline{e}\in\breve{\mathbb E}$ with $\widetilde{e}\rhd\overline{e}$ we have $[[\widetilde{e}]]\succ_J[[\overline{e}]]$ for some $J\subset{\mathcal F}$ with $|J^-|\le 1$.

   (b) For $\widetilde{e}\in{\mathbb E}$ and $\underline{\overline{e}}\in\breve{\mathbb D}$ with $[[\widetilde{e}]]\succ_J\underline{\overline{e}}$ for some $J\subset{\mathcal F}$ with $|J^-|\le 1$ there is a unique $\overline{e}\in\breve{\mathbb E}$ with $\widetilde{e}\rhd\overline{e}$ and $[[\overline{e}]]=\underline{\overline{e}}$.
   
 For $\underline{\widetilde{e}}\in{\mathbb D}$ and ${\overline{e}}\in\breve{\mathbb E}$ with $\underline{\widetilde{e}}\succ_J[[\overline{e}]]$ for some $J\subset{\mathcal F}$ with $|J^-|\le 1$ there is a unique ${\widetilde{e}}\in{\mathbb E}$ with $\widetilde{e}\rhd\overline{e}$ and $[[{\widetilde{e}}]]=\underline{\widetilde{e}}$. \end{lem}

 {\sc Proof:} (a) For $\overline{e}\in{\mathbb E}$ the statement is trivial. For $\overline{e}\in{\mathbb Z}/(q^2-1){\mathbb Z}$ take the
 $J$ defined by formula (\ref{kruthilf}).

 (b) If $(\widetilde{\ell},\widetilde{\mathcal
  I})\succ_J(\overline{\ell},\overline{\mathcal
  I})$ for $J=\emptyset$ (resp. for $J={\mathcal F}$, resp. for $J=\{i\}$ in case $\overline{\ell}=0$, resp. for $J=[i_2,i_1+f-1]$ in case $\overline{\ell}\ne 0$) then, given $\widetilde{u}$ (resp. $\overline{u}$), we find $\overline{u}$ (resp. $\widetilde{u}$) so that $(\widetilde{\ell},\widetilde{u},\widetilde{\mathcal
  I})\rhd(\overline{\ell},\overline{u},\overline{\mathcal
  I})$ by an instance of case (I) (resp. of case (II), resp. of case (III),
 resp. of case (IV)) in the definition of $\rhd$. Similarly for $(\widetilde{\ell},\widetilde{\mathcal
   I})\succ_J[[h]]$.\hfill$\Box$\\

{\bf Remark:} Lemma \ref{fernandodrei} says that $\rhd$ on $\breve{\mathbb E}$ is faithfully reflected by the relations $\succ_J$ on $\breve{\mathbb D}$ with $|J^-|\le 1$. It follows that $\rhd\rhd$ on $\breve{\mathbb E}$ is faithfully reflected by the transitive closure of the relations $\succ_J$ on $\breve{\mathbb D}$ with $|J^-|\le 1$. By Theorem \ref{tobbu} and Theorem \ref{ostersa} this transitive closure is given by the relations $\succ_J$ on $\breve{\mathbb D}$ for arbitrary $J$.

\begin{satz}\label{marienmai} (a) Given $(\overline{\ell},\overline{u},\overline{\mathcal I})\in{\mathbb E}$ and $\widetilde{\ell}\in[0,q-2]$ we have $(\widetilde{\ell},\widetilde{u},{\mathcal F})\rhd\rhd (\overline{\ell},\overline{u},\overline{\mathcal I})$ for some $\widetilde{u}$ if and only if $\widetilde{\ell}\succ_J\overline{\ell}$ and $\overline{\mathcal I}\subset\nu_{\widetilde{\ell}}^{\overline{\ell}}(J^{c,1})$ for some $J\subset{\mathcal F}$. Moreover, $\widetilde{u}$ modulo $\xi{\mathbb Z}$ is then uniquely determined.

  (b) Given $[h]\in{\mathbb Z}/(q^2-1){\mathbb Z}$ and $\widetilde{\ell}\in[0,q-2]$ we have $(\widetilde{\ell},\widetilde{u},{\mathcal F})\rhd\rhd[h]$ for some $\widetilde{u}$ if and only if $\widetilde{\ell}\succ_J^+h$ for some $J\subset{\mathcal F}$. Moreover, $\widetilde{u}$ modulo $\xi{\mathbb Z}$ is then uniquely determined.
\end{satz}

{\sc Proof:} (a) It follows from Lemma \ref{fernandodrei} and Theorem \ref{tobbu} that $(\widetilde{\ell},\widetilde{u},{\mathcal F})\rhd\rhd (\overline{\ell},\overline{u},\overline{\mathcal I})$ becomes equivalent with $(\widetilde{\ell},{\mathcal F})\succ_J(\overline{\ell},\overline{\mathcal I})$ for some $J$, together with a uniquely solvable equation for $\widetilde{u}$ modulo $\xi{\mathbb Z}$ (which in particular implies $\widetilde{\ell}-2\widetilde{u}\equiv\overline{\ell}-2\overline{u}$ modulo $\xi{\mathbb Z}$). But $(\widetilde{\ell},{\mathcal F})\succ_J(\overline{\ell},\overline{\mathcal I})$ is just the condition stated. 

(b) The argument is the same, this time using also Theorem \ref{ostersa}.\hfill$\Box$\\

Fix $\widetilde{\ell}, \overline{\ell}\in[0,q-2]$ and $J\subset{{{\mathcal F}}}$ such that $\widetilde{\ell}\succ_J\overline{\ell}$. Write $\widetilde{\ell}\equiv-\sum_{i=0}^{f-1}a_ip^i$ modulo $\xi{\mathbb Z}$, with $a_i\in[1,p]$, not all of them $=p$. Write $a_i=a_{i+f}$ for all $i\in{\mathbb Z}$. Following \cite{dediro} par. 7.2 we say that $J^c$ is maximal (with respect to $\widetilde{\ell}$) if there is no pair $i<j$ with $(a_i,\ldots,a_j)=(p,p-1,p-1,\ldots,p-1,1)$ and $[i,j-1]\subset J$ and $j\in J^c$, and if moreover $J\ne {{{\mathcal F}}}$ in the case where $a_i=p-1$ for all $i$.\footnote{In \cite{dediro} one also needs  $J\ne {{{\mathcal F}}}$ if $p=2$ and $a_i=2$ for all $i$, but here we assume $p>2$ anyway.} Consider the map $$\delta=\delta_{\overline{\ell}}:{{{\mathcal F}}}\to {{{\mathcal F}}},\quad\quad j\mapsto \Pi(i_{\mathcal D}^{\overline{\ell}}(j+1)-1).$$Write $\delta(j)=\delta(\Pi(j))$. Define $\mu(J^c)\subset{{{\mathcal F}}}$ as follows. If $\delta(J^c)\subset J^c$, then $\mu(J^c)=J^c$. Otherwise choose some $i_1\in\delta(J^c)$ with $i_1\notin J^c$. Let $j_1$ be the least
integer such that $j_1 > i_1$ and $\Pi(j_1)\in J^c$ and $\delta(j_1)=i_1$. Write $J^c= \{\Pi(j_1),\ldots,\Pi(j_r)\}$ with $j_1 < j_2 <\cdots < j_r < j_1 + f$, and define $i_{\kappa}$ for $\kappa= 2,\ldots, r$ inductively as follows: $i_{\kappa} =\delta(j_k)$ if $i_{\kappa-1} <\delta(j_{\kappa})$, but $i_{\kappa} =j_{\kappa}$ otherwise. Set $\mu(J^c)= \{\Pi(i_1),\ldots, \Pi(i_r)\}$.

  \begin{pro}\label{goldin} If $J^{c}$ is maximal\footnote{A slightly (but in general truly) weaker condition than maximality would be sufficient here.} then $\nu_{\widetilde{\ell}}^{\overline{\ell}}|_{J^{c,1}}$ is injective, and $\Pi(\nu_{\widetilde{\ell}}^{\overline{\ell}}(J^{c,1})-1)=\mu(J^{c})$.
    \end{pro}

    {\sc Proof:} As before we drop $\Pi$ from our notation. We put $\mu(j_{\kappa})=i_{\kappa}$ and we claim that $\nu_{\widetilde{\ell}}^{\overline{\ell}}(x+1)=\mu(x)+1$ for all $x\in J^c$, which then implies the Proposition.
    
   We will rely on a consideration of the $b_i$, where we write $\overline{\ell}\equiv-\sum_{i=0}^{f-1}b_ip^i$ modulo $\xi{\mathbb Z}$, with $b_i\in[1,p]$, not all of them $=p$.
    
    Let first $x\in J^c$ such that $x-1\in J$. If $\delta(x)\ne x$ then maximality of $J^c$ and a look at the possible $b_i$'s shows that $x+1\in{\mathcal D}(\widetilde{\ell})$, hence $\mu(x)=\delta(x)$. On the other hand, another look at the possible $b_i$'s (and again $\delta(x)\ne x$ and $x\in J^c$ and $x-1\in J$) also shows that $x+1\in{\mathcal D}(\widetilde{\ell})$, hence $\nu_{\widetilde{\ell}}^{\overline{\ell}}(x+1)=i_{\mathcal D}^{\overline{\ell}}(x+1)=\delta(x)+1$.
    
    If however $\delta(x)=x$ then $\nu_{\widetilde{\ell}}^{\overline{\ell}}(x+1)=x+1=\delta(x)+1$.
    
    Now let $x\in J^c$ such that also $x-1\in J^c$. If $b_x=p$ a small computation shows $a_x=1$, hence $\nu_{\widetilde{\ell}}^{\overline{\ell}}(x+1)=x=\delta(x)+1$. Similarly, if $b_x\ne p$ a small computation shows $\nu_{\widetilde{\ell}}^{\overline{\ell}}(x+1)=x+1=\delta(x)+1$.\hfill$\Box$\\
    
Now, besides $\widetilde{\ell}, \overline{\ell}$ fix $\widetilde{u}$, $\overline{u}$ in $[0,q-2]$ and $\alpha,\beta$ in $k^{\times}$. Let $W$ be the two dimensional ${\mathcal G}_F$-representation over $\overline{\mathbb F}_p$ corresponding to an element in ${\rm Ext}^1({\bf E}(\beta,\overline{u}-\overline{\ell}),{\bf E}(\alpha,\overline{u}))\otimes_k\overline{\mathbb F}_p$. Put $\chi=\omega^{-\overline{\ell}}\mu_{\beta\alpha^{-1}}:{\mathcal G}_F\to \overline{\mathbb F}_p^{\times}$ (in particular $\chi|_{{\mathcal I}_F}=\omega^{-\overline{\ell}}$), thus $W$ also defines an element in $H^1({\mathcal G}_F,\overline{\mathbb F}_p(\chi))$. In section 7 of \cite{dediro} an $\overline{\mathbb F}_p$-basis of $H^1({\mathcal G}_F,\overline{\mathbb F}_p(\chi))$ is described, indexed by ${\mathcal F}\cong{\mathbb Z}/f{\mathbb Z}$, together with an additional basis element if $\chi$ is trivial resp. cyclotomic. This is used, if $J^c$ is maximal, to define the subspace $L_{V(\widetilde{\ell},\widetilde{u})}=L^{\rm AH}_{V(\widetilde{\ell},\widetilde{u})}$ as being spanned by the basis elements with index in $\mu(J^c)$, possibly together with the said additional basis element if $\chi$ is trivial resp. cyclotomic.

    \begin{pro}\label{neupro} (\cite{gls}) Keeping $\widetilde{\ell}$, $\overline{\ell}$, $\widetilde{u}$, $\overline{u}$, $\alpha$, $\beta$, $W$, we may (if necessary) replace $J$ (while maintaining $\widetilde{\ell}\prec_J\overline{\ell}$) so that $J^c$ is maximal.
\end{pro}

{\sc Proof:} See the discussion before Proposition 8.8 in \cite{gls}, and its reformulation at the beginning of subsection 7.2 in \cite{dediro}.\hfill$\Box$\\

\begin{satz}\label{erstfreimai0} (\cite{ecgm}) $V(\widetilde{\ell},\widetilde{u})$ belongs to ${\mathbb W}_{BDJ}(W)$ if and only if the extension class of $W$ in  $H^1({\mathcal G}_F,\overline{\mathbb F}_p(\chi))$ lies in $L_{V(\widetilde{\ell},\widetilde{u})}$.
\end{satz}

 {\sc Proof:} This is Conjecture 7.2 in \cite{dediro} and Theorem 3.5.3 in \cite{ecgm}.\hfill$\Box$\\

\begin{kor}\label{movorchrihi} $|{\mathbb W}_{BDJ}(W)|=|{\mathbb W}_{\rhd}(W)|$.
\end{kor}

{\sc Proof:} The characterization of ${\mathbb W}_{BDJ}(W)$ given by Theorem \ref{erstfreimai0} is the precise parallel to the characterization of ${\mathbb W}_{\rhd}(W)$ resulting from Theorem \ref{marienmai}. Namely, containment in $L_{V(\widetilde{\ell},\widetilde{u})}$ (Theorem \ref{erstfreimai0}) inside $H^1$ corresponds to containment in the subspace inside ${\rm Ext}^1$ spanned by the lines with index in $\nu_{\widetilde{\ell}}^{\overline{\ell}}(J^{c,1})$ (Theorem \ref{marienmai}). By Proposition \ref{goldin} the respective subspaces have the same dimension.\hfill$\Box$\\
 
\begin{satz}\label{erstfreimai1} ${\mathbb W}_{BDJ}(W)={\mathbb W}_{\rm geom}(W)$.
\end{satz}

 {\sc Proof:} This was proven first in Theorem 1.2 of \cite{cegs}. Our Theorem \ref{freimai} in conjunction with Theorem \ref{marienmai} and Theorem \ref{erstfreimai0} reproves this. (But contrary to the present paper, \cite{cegs} does not require the finite extension $F/{\mathbb Q}_p$ to be unramified.)\hfill$\Box$\\

 {\bf Remark:} By Theorems \ref{freimai}, \ref{marienmai}, \ref{erstfreimai0} and \ref {erstfreimai1} the natural isomorphism\begin{gather} H^1({\mathcal G}_F,\overline{\mathbb F}_p(\chi))\cong{\rm Ext}^1({\bf E}(\beta,\overline{u}-\overline{\ell}),{\bf E}(\alpha,\overline{u}))\otimes_k\overline{\mathbb F}_p\label{yufasha}\end{gather}sends the subspace $L_{V(\widetilde{\ell},\widetilde{u})}$ to the subspace $E_{\nu_{\widetilde{\ell}}^{\overline{\ell}}(J^{c,1})}^{\alpha,\beta}\otimes_k\overline{\mathbb F}_p$. By Proposition \ref{goldin} this means that the subspace of $H^1({\mathcal G}_F,\overline{\mathbb F}_p(\chi))$ generated by the basis vectors with index in $\Pi(\nu_{\widetilde{\ell}}^{\overline{\ell}}(J^{c,1})-1)$ corresponds to the subspace of ${\rm Ext}^1({\bf E}(\beta,\overline{u}-\overline{\ell}),{\bf E}(\alpha,\overline{u}))\otimes_k\overline{\mathbb F}_p$ generated by the basis vectors with index in $\nu_{\widetilde{\ell}}^{\overline{\ell}}(J^{c,1})$, i.e. the indexings differ by a $\pm1$-shift.\\ 

 {\bf Remark:} For an irreducible two dimensional ${\mathcal G}_F$-representation $W$ it was similarly proven in \cite{cegs} that ${\mathbb W}_{BDJ}(W)={\mathbb W}_{\rm geom}(W)$. Again, Theorem \ref{freimai} in conjunction with Theorem \ref{marienmai} reproves this (but contrary to the present paper, \cite{cegs} does not require the finite extension $F/{\mathbb Q}_p$ to be unramified).\\

 {\bf Remark:} In sections \ref{okt3} and \ref{zoomgergelyfengwei} we described one-parameter degenerations of families of $(\varphi,\Gamma)$-modules along {\it specifically chosen} paths. For these families we worked out the combinatorial data in $\breve{\mathbb E}$ of its generic as well as its special members (this is encapsulated in the relation $\rhd$ on $\breve{\mathbb E}$, cf. Corollary \ref{palmmon}). The general theory seems to entail that there should be one-parameter degenerations of families of $(\varphi,\Gamma)$-modules along other paths with {\it different} limits (namely, the limiting objects which could possibly occur should be encapsulated in the relation $\rhd\rhd$ on $\breve{\mathbb E}$). However, we find it hard to work out explicitly such more general one-parameter families.

 \section{Multiplicities}

 \label{giovanni}

By construction, for a two dimensional ${\mathcal G}_F$-representation $W$ (irreducible or not), the elements of ${\mathbb W}_{\rhd}(W)$ come along with a canonically defined {\it multiplicity}. In other words, ${\mathbb W}_{\rhd}(W)$, by construction, can be seen as a {\it multiset}.

Due to Lemma \ref{fernandodrei}, the definition and analysis of these multiplicities can be reduced to the combinatorics of the set $\breve{\mathbb D}$ in section \ref{combi}. Namely, we propose:\\

 {\bf Definition:} For $d\in\breve{\mathbb D}$ and $\ell\in[0,q-2]$, the {\it multiplicity of $\ell$ in $d$} is the number $m({\ell}|d)$ of subsets ${\mathcal I}$ of ${{{\mathcal F}}}$ for which there is a $J\subset {{\mathcal F}}$ with $(\ell,{\mathcal I})\succ_Jd$ and so that there is no proper subset ${\mathcal I}'\subsetneq{\mathcal I}$ with $(\ell,{\mathcal I}')\succ_Jd$.

Equivalently, $m(\ell|d)$ is the number of subsets $J\subset {{{\mathcal F}}}$ with $(\ell,{\mathcal I})\succ_Jd$ for some ${\mathcal I}$.\\

{\bf Remarks:} (a) Proposition \ref{neupro} accounts for multiplicities possibly $>1$.

(b) If $d=(\overline{\ell},\overline{\mathcal I})$, asking for $m(\ell|d)$ is asking for the number of all $J$ with $\overline{\mathcal I}\subset\nu_{{\ell}}^{\overline{\ell}}(J^{c,1})$, where ${\ell}=\ell(J)$ is such that ${\ell}\succ_J\overline{\ell}$. Is $J$ uniquely determined by $\nu_{{\ell}}^{\overline{\ell}}(J^{c,1})$ ?

(c) Problem: Characterize the sets $\nu_{{\ell}}^{\overline{\ell}}(J^{c,1})$ for a given $(\overline{\ell},\overline{\mathcal I})\in {\mathbb D}$.\\

{\bf Remark:} Recall the strata ${\mathcal V}(\ell,u,{\mathcal I})$ of ${\mathcal X}_{2,{\rm red}}$ described in subsection \ref{stackstrataset}. Let $W$ be a two dimensional ${\mathcal G}_F$-representation with corresponding $(\varphi,\Gamma)$-module ${\bf D}$. Let $e=[h]$ resp. $e=(\overline{\ell},\overline{u},\overline{\mathcal I})$ be such that ${\bf D}\in {\mathcal V}(h)$ resp. ${\bf D}\in{\mathcal V}^0(\overline{\ell},\overline{u},\overline{\mathcal I})$ (cf. formula (\ref{fernandovertex})). Denote by $d\in\breve{\mathbb D}$ the image element of $e\in \breve{\mathbb E}$ under the map (\ref{fast}). Let $\ell,u\in[0,q-2]$ be such that $({\ell},{u},{\mathcal F})\rhd\rhd d$. Then $m(\ell|d)$ should be the number of minimal subsets ${\mathcal I}\subset{\mathcal F}$ such that the family of $(\varphi,\Gamma)$-modules carried by ${\mathcal V}(\ell,u,{\mathcal I})$ degenerates into ${\bf D}$.

 \end{document}